\documentclass[a4paper,12pt]{article} 
\usepackage{amssymb}
\usepackage{amscd}
\usepackage{amsmath} 
\usepackage[dvipdfmx]{graphicx}
\usepackage{latexsym} 
\usepackage{enumerate}

\usepackage{theorem}
\theoremstyle{plain}
\theorembodyfont{\itshape}
\newtheorem{theorem}{Theorem}[section]
\newtheorem{proposition}[theorem]{Proposition}
\newtheorem{lemma}[theorem]{Lemma}
\newtheorem{corollary}[theorem]{Corollary} 
\newtheorem{result}[theorem]{Result} 
\theorembodyfont{\rmfamily}

\newtheorem{remark}[theorem]{Remark}
\newtheorem{problem}[theorem]{Problem}

\title{\bf On Some Hypergeometric Summations\thanks{ 
formula; closed-form expression; asymptotic analysis; Diophantine 
approximation; contiguous relation.}} 
\author{Katsunori Iwasaki\thanks{Department of Mathematics, 
Hokkaido University, Kita 10, Nishi 8, Kita-ku, Sapporo 060-0810 Japan. 
E-mail: {\tt iwasaki@math.sci.hokudai.ac.jp}}}
\date{January 6, 2015} 
\begin{document}
\maketitle
\begin{abstract}
We develop a theoretical study of non-terminating hypergeometric 
summations with one free parameter.  
Composing various methods in complex and asymptotic analysis,  
geometry and arithmetic of certain transcendental curves and 
rational approximations of irrational numbers, we are able to 
obtain some necessary conditions of arithmetic flavor for 
a given hypergeometric sum to admit a gamma product formula. 
This kind of research seems to be new even in the most 
classical case of the Gauss hypergeometric series.
\end{abstract} 
\section{Introduction} \label{sec:intro}
Let ${}_2F_1(\alpha,\beta;\gamma;z)$ be the Gauss hypergeometric series:  
\begin{equation*} \label{eqn:ghg}
{}_2F_1(\alpha,\beta;\gamma;z) := \sum_{k=0}^{\infty} 
\dfrac{(\alpha)_k (\beta)_k}{(\gamma)_k(1)_k} \, z^k,   
\qquad (\alpha)_k := \alpha(\alpha+1)\cdots(\alpha+k-1).        
\end{equation*}
Given a parameter $\lambda = (p,q,r;a,b;x) \in \mathbb{C}^5 \times 
\mathbb{D}$, where $\mathbb{D} = \{x \in \mathbb{C} \,:\, |x| < 1\}$ is 
the open unit disk in the complex plane $\mathbb{C}$, we consider an 
entire meromorphic function of $w$:  
\begin{equation} \label{eqn:f}  
f(w) = f(w;\lambda) := {}_2F_1(pw+a, \, qw+b; \, rw; \, x).        
\end{equation} 
\par
We say that $f(w)$ admits a {\sl gamma product formula} (GPF),  
if there exist a rational function $S(w) \in \mathbb{C}(w)$; a nonzero 
complex number $d \in \mathbb{C}^{\times}$; two nonnegative integers 
$m, n \in \mathbb{Z}_{\ge0}$; $m$ complex numbers 
$u_1, \dots, u_m \in \mathbb{C}$; and $n$ complex numbers 
$v_1, \dots, v_n \in \mathbb{C}$ such that  
\begin{equation} \label{eqn:gpf}
f(w) = S(w) \cdot d^w \cdot 
\dfrac{\varGamma(w+u_1)\cdots\varGamma(w+u_m)}{\varGamma(w+v_1)\cdots\varGamma(w+v_n)},     
\end{equation}
where $\varGamma(w)$ is the Euler gamma function. 
We shall discuss the following.    
\begin{problem} \label{prob:gpf} 
Find a parameter $\lambda = (p,q,r;a,b;x)$ for which $f(w)$ admits a GPF. 
\end{problem}
\par
A bit of search in the classical literature gives us three solutions 
in Table \ref{tab:erdelyi}, where $x = e^{\pm \pi i/3}$ is allowed to 
lie on the unit circle in the third solution.  
\begin{table}[t]
\begin{align*}
{}_2F_1\left(-\frac{w}{2}+\frac{3}{4},-\frac{w}{2}+\frac{5}{4};w;-\frac{1}{3}\right) 
&= \left(\frac{8}{9}\right)^{w-3/2} 
\frac{\varGamma(4/3)\varGamma(w)}{\varGamma(3/2)\varGamma(w-1/6)}, \\
{}_2F_1\left(3 w-\frac{5}{6},3 w-\frac{1}{3};3 w;\frac{1}{9}\right) 
&= \frac{2^{17/18}}{3^{5/6}} \left(\frac{3^3}{2^4}\right)^w 
\frac{\varGamma(w+5/36)\varGamma(w+23/36)}{\varGamma(w+2/9)\varGamma(w+5/9)}, \\
{}_2F_1\left(w,3 w-1;2 w;e^{\pm \pi i/3}\right) 
&= e^{\mp \pi i/6} \left(2^2 \cdot 3^{-3/2} \cdot e^{\pm \pi i/2}\right)^w 
\frac{\varGamma(1/2)\varGamma(w+1/2)}{\varGamma(2/3)\varGamma(w+1/3)}. 
\end{align*}
\caption{Three solutions to Problem \ref{prob:gpf} found in  
Erd\'{e}lyi \cite{Erdelyi}.}
\label{tab:erdelyi}
\end{table}
They are readily derived from formulas (53), (54) and (55)+(56) of   
Erd\'{e}lyi \cite[Chapter II, \S 2.8]{Erdelyi} by affine changes of 
a variable and the duplication formula for the gamma function. 
A more extensive search in the literature would give us more solutions. 
Problem \ref{prob:gpf} naturally extends to the generalized hypergeometric 
series ${}_{p+1}F_p$.  
Our ultimate goal is to enumerate all solutions to the problem, 
but it is far beyond the scope of this article even in the most 
classical case of ${}_2F_1$, being simply too hard to settle at present.  
Our aim here is more moderate, that is, to embark on a theoretical 
study of this very classical problem, but in a direction that 
has hitherto attracted scant attention.        
\par
The subject of hypergeometric evaluations has a long history.  
Recently it saw an important progress with the development of 
Wilf-Zeilberger methods \cite{WZ} and Zeilberger's 
algorithms \cite{Zeilberger1,Zeilberger2}, which enabled 
systematic proofs (and sometimes discoveries) of 
a lot of combinatorial identities (see also \cite{Koepf2,PWZ}). 
Gessel \cite{Gessel} and Koepf \cite{Koepf} applied these 
techniques to terminating hypergeometric sums to give mechanical proofs 
of classical evaluations in Bailey's book \cite{Bailey} as well as 
Gosper, Gessel and Stanton's ``strange" evaluations \cite{GS}; see 
also Apagodu and Zeilberger \cite{AZ} and Ekhad \cite{Ekhad} 
for more evaluations. 
In some cases, formulas for terminating series remain true for 
non-terminating ones, due to Carlson's theorem in function 
theory \cite[p. 39]{Bailey}. 
These methods might be useful for our purpose in one direction, that is, 
toward finding and proving as many solutions as possible, but {\sl not} 
in the other direction, that is, toward establishing necessary 
conditions, as strong as possible, for a given $\lambda = (p,q,r;a,b;x)$ 
to be a solution. 
It is expected that $\lambda$ must be subject to some severe 
constraints, perhaps of arithmetic flavor. 
In the latter direction we need more transcendental methods 
based on complex and asymptotic analysis, geometry and arithmetic 
of certain transcendental curves, as well as on rational 
approximations of irrational numbers. 
These are exactly what we want to explore in this article.  
\par
Problem \ref{prob:gpf} has a close relative.   
To state it we say that $f(w)$ is {\sl of closed form} if 
\begin{equation} \label{eqn:ocf}
\dfrac{f(w+1)}{f(w)} =: R(w) \in \mathbb{C}(w) \,\,:\,\, 
\mbox{a rational function of $w$}. 
\end{equation}
\begin{problem} \label{prob:ocf} 
Find a parameter $\lambda = (p,q,r;a,b;x)$ for which 
$f(w)$ is of closed form.  
\end{problem}
\par
Any solution to Problem \ref{prob:gpf} leads to a solution to 
Problem \ref{prob:ocf}. 
Indeed, by the recurrence formula $\varGamma(w+1) = w \, \varGamma(w)$ for the 
gamma function, if $f(w)$ admits the product formula 
\eqref{eqn:gpf} then it fulfills condition \eqref{eqn:ocf} with 
rational function 
\begin{equation} \label{eqn:canonical}
R(w) = \dfrac{S(w+1)}{S(w)} \cdot d \cdot 
\dfrac{(w+u_1)\cdots(w+u_m)}{(w+v_1)\cdots(w+v_n)}.  
\end{equation}
\par
Now we can pose a converse question. 
To formulate it, notice that any rational function $R(w) \in \mathbb{C}(w)$ 
can be written in the form \eqref{eqn:canonical} such that  
$S(w) \in \mathbb{C}(w)$ is a rational function, $u_i - v_j$ is not an 
integer for every $i = 1, \dots,m$ and $j = 1, \dots, n$, and moreover if 
$R(w) \in\mathbb{R}(w)$ then $S(w)$, $d$, $P(w) := (w+u_1)\cdots(w+u_m)$ and 
$Q(w) := (w+v_1)\cdots(w+v_n)$ should be real (see Lemma \ref{lem:rational1}).  
Such a representation \eqref{eqn:canonical} is said to be {\sl canonical}. 
Note that in formula \eqref{eqn:canonical} one can multiply $S(w)$ by a 
nonzero (real) constant (if $R(w)$ is real) without changing the form.       
\begin{problem} \label{prob:back}
Does a solution $f(w)$ to Problem \ref{prob:ocf} lead back to a 
solution to Problem \ref{prob:gpf}?   
Namely, suppose that $f(w)$ satisfies condition \eqref{eqn:ocf} with a 
canonical representation \eqref{eqn:canonical} and multiply $S(w)$ by a 
suitable nonzero constant.   
Then does $f(w)$ admit GPF \eqref{eqn:gpf} ? 
\end{problem} 
We shall give an affirmative solution when  $\lambda$ lies in a certain 
real region (see Theorem \ref{thm:gpf}). 
\par
There is a method of finding (partial) solutions to 
Problem \ref{prob:ocf}, which we call {\sl the method of contiguous 
relations}. 
It works only when $p, q, r \in \mathbb{Z}$ and relies on the fifteen 
contiguous relations of Gauss (see e.g. Andrews et. al. \cite[\S2.5]{AAR}).    
Composing a series of contiguous relations yields      
\begin{equation} \label{eqn:contig-F} 
\begin{split}
{}_2F_1(\alpha + p, \beta + q; \gamma + r; z)  
&= r(\alpha,\beta;\gamma;z) \, {}_2F_1(\alpha, \beta; \gamma; z) \\
&\phantom{=} + q(\alpha,\beta;\gamma;z) \, {}_2F_1(\alpha+1, \beta+1; \gamma+1; z), 
\end{split}
\end{equation}
where $q(\alpha,\beta;\gamma;z)$ and $r(\alpha,\beta;\gamma;z)$ are rational 
functions of $(\alpha,\beta;\gamma;z)$ depending uniquely upon $(p,q,r)$. 
Recently Vidunas \cite{Vidunas1} and Ebisu \cite{Ebisu1} discussed some 
computational aspects of formula \eqref{eqn:contig-F} and showed how to 
compute $q(\alpha,\beta;\gamma;z)$ and $r(\alpha,\beta;\gamma;z)$ rapidly 
and efficiently. 
\par
Given a parameter $\lambda = (p,q,r;a,b;x)$ with $p, q, r \in \mathbb{Z}$, we put     
\begin{equation} \label{eqn:f-tilde}  
\tilde{f}(w) = \tilde{f}(w;\lambda) := {}_2F_1(pw+a+1, \, qw+b+1; \, rw+1; \, x).         
\end{equation} 
Substituting $(\alpha,\beta;\gamma;z) = (pw+a, qw+b; rw;x)$ into formula 
\eqref{eqn:contig-F} we have  
\begin{equation} \label{eqn:contig-f}
f(w+1) = R(w) \, f(w) + Q(w) \, \tilde{f}(w),  
\end{equation}
where $Q(w) = Q(w;\lambda)$ and $R(w) = R(w;\lambda)$ are rational 
functions of $w$ depending on $\lambda$.  
If $\lambda$ happens to be such a parameter that $Q(w)$ vanishes 
identically in $\mathbb{C}(w)$, then three-term relation \eqref{eqn:contig-f} 
reduces to the two-term one \eqref{eqn:ocf} so that $\lambda$ happens 
to be a solution to Problem \ref{prob:ocf}. 
The method of contiguous relations is developed by Ebisu \cite{Ebisu3,Ebisu2} 
(mostly for terminating series) and it will be amplified    
for non-terminating ones in this article. 
We say that a solution to Problem \ref{prob:ocf} {\sl comes from 
contiguous relations} if it is obtained by this method. 
\begin{problem} \label{prob:contiguous1}
When does a solution to Problem \ref{prob:ocf} come from 
contiguous relations? 
\end{problem} 
\par
In a certain real parameter region we will be able to show that all solutions 
{\sl essentially} come from contiguous relations. 
For the precise statement of this result, including what we mean by {\sl essentially}, 
we refer to Theorem \ref{thm:contiguous1} and a comment right after it.  
\par
\begin{figure}[t]
\begin{center}
\includegraphics[width=60mm,clip]{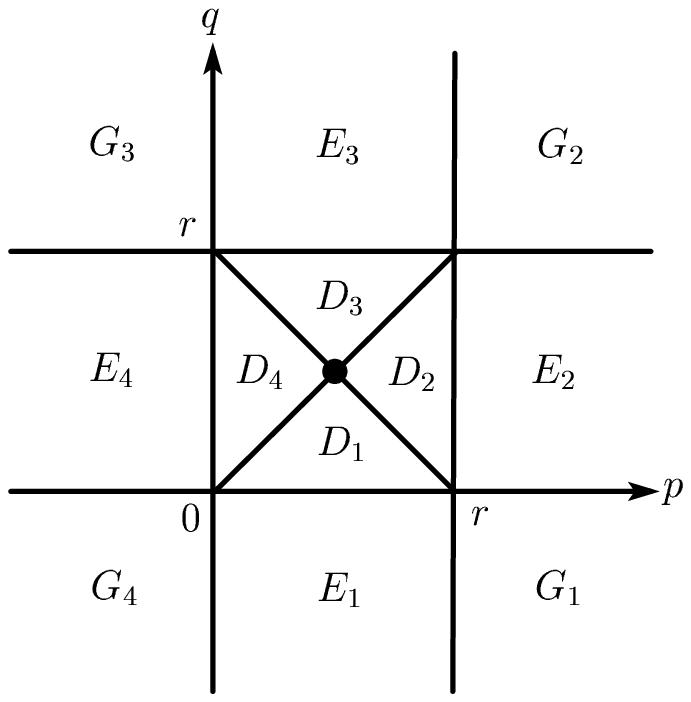}
\end{center}
\caption{The real $(p,q)$-plane for a fixed $r > 0$.} 
\label{fig:pq-plane} 
\end{figure}

Problems \ref{prob:gpf} and \ref{prob:ocf} are difficult for a 
general complex parameter $\lambda = (p,q,r;a,b;x)$. 
So we content ourselves to suppose that $\lambda$ be {\sl real} and 
furthermore to restrict it into a real region  
\begin{equation} \label{eqn:cross} 
p, \, q, \, r \in \mathbb{R}, \quad 0 < p < r \quad \mbox{or} \quad  
0 < q < r; \qquad a, \, b \in \mathbb{R};  \qquad -1 < x < 1,    
\end{equation}     
whose $(p,q)$-component forms the cross-shaped domain $D \cup E$ 
in Figure \ref{fig:pq-plane}, where $D := D_1 \cup D_2 \cup D_3 \cup 
D_4$ and $E := E_1 \cup E_2 \cup E_3 \cup E_4$. 
Of particular interest among \eqref{eqn:cross} is a subregion   
\begin{equation} \label{eqn:square} 
p, \, q, \, r \in \mathbb{R}, \quad 0 < p < r \quad \mbox{and} \quad  
0 < q < r; \qquad a, \, b \in \mathbb{R};  \qquad -1 < x < 1,    
\end{equation}     
whose $(p,q)$-component comprises the central square $D$ in 
Figure~\ref{fig:pq-plane}. 
The main results of this article will be stated either in region 
\eqref{eqn:cross} or \eqref{eqn:square} and other parameter regions 
will be left elsewhere.  
\par
The hypergeometric series enjoys well-known 
$\mathbb{Z}_2 \ltimes (\mathbb{Z}_2 \times \mathbb{Z}_2)$ symmetries: 
\begin{subequations} \label{eqn:transf}
\begin{align}
{}_2F_1(\alpha, \, \beta; \, \gamma; \, z)  
&= {}_2F_1(\beta, \, \alpha; \, \gamma; \, z),  \label{eqn:transf0} \\ 
&= (1-z)^{\gamma-\alpha-\beta} {}_2F_1(\gamma-\alpha, \, \gamma-\beta; \, \gamma; \, z),  
\label{eqn:transf1} \\
&= (1-z)^{-\alpha} {}_2F_1(\alpha, \, \gamma-\beta; \, \gamma; \, z/(z-1)),  
\label{eqn:transf2} \\
&= (1-z)^{-\beta} {}_2F_1(\gamma-\alpha, \, \beta; \, \gamma; \, z/(z-1)),  
\label{eqn:transf3}    
\end{align}
\end{subequations}
where transformation \eqref{eqn:transf0} is obvious, \eqref{eqn:transf1} 
is due to Euler, whereas \eqref{eqn:transf2} and \eqref{eqn:transf3} are 
due to Pfaff (see \cite[Theorem 2.2.5]{AAR}). 
They induce symmetries of solutions to Problem \ref{prob:gpf} or 
\ref{prob:ocf}:  
\begin{subequations} \label{eqn:sym}
\begin{align}
\lambda = (p,q,r; a,b; x)  
& \mapsto (q,p,r; b,a; x), \label{eqn:sym0} \\
& \mapsto (r-p, \, r-q, \, r; \, -a, \, -b; \, x), \label{eqn:sym1} \\
& \mapsto (p, \, r-q, \, r; \, a, \, -b; \,x/(x-1)),  \label{eqn:sym2} \\
& \mapsto (r-p, \, q, \, r; \, -a, \, b; \, x/(x-1)). \label{eqn:sym3}   
\end{align}
\end{subequations}
\begin{lemma} \label{lem:reductions}
By symmetries \eqref{eqn:sym} the parameter region 
\eqref{eqn:cross} can be reduced to a subregion  
\begin{equation} \label{eqn:pencil} 
p, \, q, \, r \in \mathbb{R}, \quad 0 < p < r \,\,\, \mbox{and} \,\,\,   
q \le p \,\,\, \mbox{and} \,\,\, p+q \le r; \qquad a, \, b \in \mathbb{R}; 
\qquad 0 \le x < 1,    
\end{equation}     
whose $(p,q)$-component forms the ``pencil-like" domain $D_1 \cup E_1$ 
in Figure $\ref{fig:pq-plane}$. 
In a similar manner the parameter region \eqref{eqn:square} can be 
reduced to a subregion 
\begin{equation} \label{eqn:apollo} 
p, \, q, \, r \in \mathbb{R}, \quad 0 < q \le p \,\,\, \mbox{and} \,\,\, 
p+q \le r; \qquad a, \, b \in \mathbb{R};  \qquad 0 \le x < 1,    
\end{equation}     
whose $(p,q)$-component corresponds to the isosceles right triangle 
$D_1$ in Figure $\ref{fig:pq-plane}$. 
\end{lemma}
{\it Proof}. First we claim that condition \eqref{eqn:cross} can be 
reduced to 
\begin{equation} \label{eqn:h-strip} 
p, \, q, \, r \in \mathbb{R}, \quad 0 < p < r \,\,\, \mbox{and} \,\,\, q < r; 
\qquad a, \, b \in \mathbb{R};  \qquad 0 \le x < 1,    
\end{equation}
whose $(p,q)$-component corresponds to the region $D \cup E_1$ in 
Figure \ref{fig:pq-plane}. 
Indeed, we may assume $0 < p < r$ in condition \eqref{eqn:cross}, 
for otherwise transformation \eqref{eqn:sym0} takes $0 < q < r$ to  $0 < p < r$. 
Then $\lambda = (p,q,r; a,b; x)$ is brought to region 
\eqref{eqn:h-strip} by transformation \eqref{eqn:sym1} if 
$q \ge r$ and $0 \le x <1$; by \eqref{eqn:sym2} if $q \ge r$ and $-1 < x < 0$; 
and by \eqref{eqn:sym3} if $q < r$ and $-1 < x < 0$, respectively. 
Here the map $q \mapsto r-q$ exchanges the conditions  
$q \ge r$ and $q \le 0$ $(< r)$, whereas $x \mapsto x/(x-1)$ 
maps the negative interval $-1 < x < 0$ to the positive one $0 < x < 1/2$.  
Next we show that condition \eqref{eqn:h-strip} can be reduced to 
condition \eqref{eqn:pencil}. 
If $(p,q) \in E_1$ then we are already done. 
When $(p,q) \in D$, we use transformations \eqref{eqn:sym0} and/or 
\eqref{eqn:sym1} which keep $x$ invariant. 
If $(p,q) \in D_2 \cup D_3$, apply \eqref{eqn:sym1} 
to make $(p,q) \in D_1 \cup D_4$; and if $(p,q) \in D_4$, apply 
\eqref{eqn:sym0} to have $(p,q) \in D_1$. 
These procedures also reduce condition \eqref{eqn:square} to 
condition \eqref{eqn:apollo}.  \hfill $\Box$ \par\medskip
In what follows we shall always exclude the case $x=0$, for which 
everything is trivial. 
\section{Main Results} \label{sec:results} 
In this section we present the main results together with an outline 
of this article. 
\subsection{Elementary Solutions} \label{subsec:elementary}
A solution to Problem \ref{prob:gpf} or \ref{prob:ocf} is said to be 
{\sl elementary} if the corresponding function $f(w)$ in definition  
\eqref{eqn:f} has at most finitely many poles in the entire complex 
$w$-plane $\mathbb{C}_w$. 
\begin{theorem} \label{thm:elementary}
In region \eqref{eqn:cross} Problem $\ref{prob:ocf}$ has only two types 
of elementary solutions. 
Up to symmetries the first type of elementary solutions are given by  
\begin{equation} \label{eqn:elementary1}
q = 0, \qquad b \in \mathbb{Z}_{\le 0},    
\end{equation}
in which case $f(w)$ itself is a rational function because the series 
\eqref{eqn:f} that defines it is terminating.  
These solutions are attached to the boundary of two strips 
$E_1 \cup D \cup E_3$ and $E_2 \cup D \cup E_4$ in 
Figure $\ref{fig:pq-plane}$.  
The second type of elementary solutions are given by  
\begin{equation} \label{eqn:elementary2}
p = q = \frac{r}{2} > 0, \qquad 
a = i,  \quad b = j -\frac{1}{2}, \qquad i, \, j \in \mathbb{Z},      
\end{equation}
where $a$ and $b$ are permutable by symmetry. 
They are attached to the ``core" of the pencil $D_1 \cup E_1$, namely, 
to the bullet $\bullet$ in Figure $\ref{fig:pq-plane}$.     
Via the scaling transformation $w \mapsto w/r$, 
solutions \eqref{eqn:elementary2} correspond to a contiguous family of 
degenerate Gauss hypergeometric functions     
\begin{equation} \label{eqn:dihedral} 
f(w/r) = {}_2F_1 \left(\frac{w}{2}+i, \, \frac{w-1}{2}+j, \, w; \, x \right) = 
S_{ij}(w;x) \cdot \left(\frac{1+\sqrt{1-x}}{2} \right)^{1-w},   
\end{equation}
where $S_{00}(w;x) = 1$ and $S_{ij}(w;x)$ is a rational function of 
$w$ defined by  
\[
S_{ij}(w;x) = (1-x)^{-\frac{i+j}{2}} 
F_3\Big(i+j,\, j-i;\, 1-i-j,\, 1+i-j;\, w;\, 
-\frac{1-\sqrt{1-x}}{2 \sqrt{1-x}}, \, \frac{1-\sqrt{1-x}}{2} \, \Big),    
\] 
with $F_3(\alpha_1, \alpha_2; \beta_1, \beta_2; \gamma; u, v)$ being 
Appell's hypergeometric series of two variables: 
\[
F_3(\alpha_1, \alpha_2; \beta_1, \beta_2; \gamma; u, v) 
:= \sum_{i=0}^{\infty} \sum_{j=0}^{\infty} 
\frac{(\alpha_1)_i(\alpha_2)_j(\beta_1)_i(\beta_2)_j}{(\gamma)_{i+j} \, i! \, j!} \, 
u^i v^j.    
\]
\end{theorem}
\par
Indeed, $S_{ij}(w;x)$ is rational in $w$ because the Appell  
series that defines it is terminating for every $i,j \in \mathbb{Z}$.  
Functions \eqref{eqn:dihedral} are instances of Gauss hypergeometric functions 
with dihedral monodromy groups studied by Vidunas \cite{Vidunas2}. 
Theorem \ref{thm:elementary} will be proved in Proposition \ref{prop:f-pole}.   
\subsection{Some Examples of Non-Elementary Solutions} \label{subsec:non-elementary}
\begin{table}[tt]
\[
\begin{array}{c|c|c|c|c|c|l}
\hline
  &   &   &     &     &     &                    \\[-4mm]
r & p & q & x   & a   & b   & f(w+1)/f(w) = R(w) \\[1mm]
\cline{1-7}
  &   &   &              &     &              &                                                                                 \\[-4mm]  
4 & 1 & 1 & \dfrac{8}{9} & 0   & \dfrac{1}{4} & \dfrac{4}{3} \cdot \dfrac{(w+2/4)(w+3/4)}{1\cdot 1 \cdot(w+2/3) \cdot (w+7/12)} \\[3mm]
\cline{5-7}
  &   &   &              &     &              &                                                                         \\[-4mm]   
&   &   &     & \dfrac{1}{2}   & \dfrac{1}{4} & \dfrac{4}{3} \cdot \dfrac{w(w+3/4)}{1\cdot 1\cdot (w+7/12) \cdot (w+1/6)} \\[3mm]
\cline{5-7}
  &   &   &              &     &              &                                                                          \\[-4mm]   
  &   &   &     & 0   & \dfrac{1}{2} & \dfrac{4}{3} \cdot \dfrac{(w+1/4)(w+3/4)}{1\cdot 1\cdot (w+1/3)(w+2/3) \cdot 1}   \\[3mm]  
\cline{1-7}
  &   &   &              &     &              &                                                                                 
\\[-4mm]   
6 & 1 & 1 & \dfrac{4}{5} & 0  & \dfrac{1}{2} & \dfrac{3^6}{5^4} \cdot \dfrac{(w+1/6)(w+2/6)(w+4/6)(w+5/6)}{1 \cdot 1\cdot (w+1/5)(w+4/5) 
\cdot  (w+3/10)(w+7/10)} \\[3mm]
\cline{5-7}
  &   &   &              &     &              &                                                                                 
\\[-4mm]   
  &   &   &     & \dfrac{2}{3} & \dfrac{1}{6} & \dfrac{3^6}{5^4} \cdot \dfrac{w(w+2/6)(w+3/6)(w+5/6)}{1 \cdot 1\cdot (w+17/30)(w+23/30) 
\cdot (w+1/15)(w+4/15)} \\[3mm]
\cline{2-7}
  &   &   &              &     &              &                                                                                 
\\[-4mm]   
  & 2 & 2 & \dfrac{3}{4}(3-\sqrt{3}) & 0 & \dfrac{1}{3} & \dfrac{3\sqrt{3}}{2} \cdot \dfrac{(w+2/6)(w+5/6)}{1 \cdot 1\cdot 
(w+3/4)\cdot (w+5/12)} \\[3mm]  
\cline{2-7}
  &   &   &              &     &              &                                                                                 
\\[-4mm]   
  & 3 & 1 & 4(\sqrt{5}-2) & 0 & \dfrac{1}{6} & \dfrac{27}{125}(5+2\sqrt{5}) \cdot \dfrac{(w+3/6)(w+5/6)}{1 \cdot 1\cdot 1\cdot 
(w+17/30)(w+23/30)} \\[3mm]
\cline{5-7}
  &   &   &              &     &              &                                                                                 
\\[-4mm]   
  &   &   &               & 0 & \dfrac{1}{2} & \dfrac{27}{125}(5+2\sqrt{5}) \cdot \dfrac{(w+1/6)(w+5/6)}{1\cdot 1\cdot 1\cdot 
(w+3/10)(w+7/10)} \\[3mm]
\cline{1-7}
  &   &   &              &     &              &                                                                                 
\\[-4mm]   
8 & 4 & 2 & 4(3\sqrt{2}-4) & 0 & \dfrac{1}{4} & \dfrac{4}{27}(17+12\sqrt{2}) \cdot \dfrac{(w+3/8)(w+7/8)}{1 \cdot 1 \cdot 
1\cdot (w+11/24)(w+19/24)}  \\[3mm]
\hline  
\end{array}
\]
\caption{Some non-elementary solutions to Problem \ref{prob:ocf} in region 
\eqref{eqn:apollo}.} 
\label{tab:solutions}
\end{table}
Table \ref{tab:solutions} exhibits some non-elementary solutions to 
Problem \ref{prob:ocf} in region \eqref{eqn:square}; they are presented 
in region \eqref{eqn:apollo} upon reduced by symmetries \eqref{eqn:sym}.  
These examples are obtained by the method of contiguous relations. 
They also lead to solutions to Problem \ref{prob:gpf} due to 
Theorem \ref{thm:gpf} below. 
There are also non-elementary solutions to Problem \ref{prob:gpf} 
not lying in region \eqref{eqn:square}. 
For example, 
\[
\begin{split}
{}_2F_1(w,-w+1/2;3w;1/2) &= \dfrac{\csc(\pi/8)}{\sqrt{6}} \cdot  
\left(\frac{3^3}{2^5}\right)^w \cdot 
\frac{\varGamma(w+1/3) \varGamma(w+2/3)}{\varGamma(w+3/8)\varGamma(w+5/8)}, \\
{}_2F_1(w+1/3,-w+7/6;3w;1/2) &= \dfrac{2^{5/6}\csc(3\pi/8)}{\sqrt{3}} \cdot 
\left(\frac{3^3}{2^5}\right)^w \cdot 
\frac{\varGamma(w) \varGamma(w+2/3)}{\varGamma(w+5/24)\varGamma(w+11/24)} 
\end{split}
\]
are solutions in region \eqref{eqn:cross}, whose $(p,q)$-component lies 
in domain $E_1$ of Figure \ref{fig:pq-plane}.   
\begin{problem} \label{prob:necessary} 
Find necessary conditions on $\lambda = (p,q,r;a,b,;x)$ as well as on 
the related quantities $d$; $m$, $n$; $u_1, \dots, u_m$ and $v_1, \dots, v_n$, 
in order for $\lambda$ to be a non-elementary solution to 
Problem \ref{prob:gpf} or \ref{prob:ocf}. 
Of particular interests to us are {\sl arithmetic properties} of 
these quantities.  
\end{problem}
\par
Of course our ultimate goal is to establish necessary and sufficient 
conditions in the entire parameter region, but it is far beyond the 
scope of this article as mentioned in \S1.  
We should content ourselves to establish some necessary conditions 
upon restricting $\lambda$ to subregion \eqref{eqn:cross} or \eqref{eqn:square}.  
``Arithmetic properties" in Problem \ref{prob:necessary} refer to such 
questions as follows. 
Observe that $p$, $q$ and $r$ in Table \ref{tab:solutions} are {\sl integers};   
they are trivially so because they come from contiguous relations. 
But to what extent is this true for a {\sl general} solution? 
According to Table \ref{tab:solutions}, the numbers $a$ and $b$ are 
{\sl rational}, while $x$ and $d$ are {\sl algebraic};  
moreover, $u_1, \dots, u_m$ and $v_1, \dots, v_n$ are 
also {\sl rational}, with $m = n$. 
Do these observations remain true in general?  
If so, what kinds of integers occur as the numerators and denominators 
of those rational numbers?             
\subsection{Arithmetic Properties} 
\label{subsec:integer}
A first thing that can be said about Problem \ref{prob:necessary} in 
region \eqref{eqn:cross} is the following. 
\begin{theorem} \label{thm:gpf} 
Let $\lambda =(p,q,r;a,b;x)$ be a solution to Problem 
$\ref{prob:ocf}$ in region \eqref{eqn:cross} and write $R(w)$  
in a canonical form \eqref{eqn:canonical}. 
Then the gamma product formula \eqref{eqn:gpf} is valid with $m = n$. 
This means that in region \eqref{eqn:cross} any solutions to 
Problem $\ref{prob:ocf}$ leads back to a solution to 
Problem $\ref{prob:gpf}$.  
If the solution is non-elementary then $r$ must be integer with 
$1 \le m = n \le r$ and there exist $m$ integers $s_1, \dots, s_m$, 
mutually distinct modulo $r$, such that
\begin{equation} \label{eqn:ui} 
u_i = s_i/r \qquad (i = 1, \dots, m).
\end{equation}
\end{theorem}
\par
An essence for the proof of this theorem lies in asymptotic analysis. 
It consists of investigating the asymptotic behavior of 
$f(w)$ (see Proposition \ref{prop:asympt}) and that of a gamma 
product expression as in the right-hand side of formula \eqref{eqn:gpf} 
(in \S\ref{sec:gamma-p}) and then comparing both results  
(see Propositions \ref{prop:f=g} and \ref{prop:pole2}).     
Theorem \ref{thm:gpf} will be established at the end of \S\ref{sec:rational}   
right after Proposition \ref{prop:pole2}.  
It turned out that $r$ must be an integer for any non-elementary solution.  
How about $p$ and $q$?  
\begin{table}[tt]
\[
\begin{array}{c|c|c|c|c|c|l}
\hline
  &   &   &     &     &     &                    \\[-4mm]
r & p & q & x   & a   & b   & f(w+1)/f(w) = R(w) \\[1mm]
\cline{1-7}
  &   &   &     &     &     &                    \\[-4mm] 
3 & \dfrac{1}{2} & \dfrac{1}{2} & \dfrac{4}{5} & 0  & \dfrac{1}{2} & \dfrac{3^3}{5^2} \cdot 
\dfrac{(w+1/3)(w+2/3)}{(w+2/5)(w+3/5)} \\[3mm]
\cline{5-7}
  &   &   &     &     &     &                    \\[-4mm] 
  &   &   &     & \dfrac{2}{3} & \dfrac{1}{6} & \dfrac{3^3}{5^2} \cdot 
\dfrac{w(w+2/3)}{(w+2/15)(w+8/15)} \\[3mm]
\hline 
\end{array}
\]
\caption{Two non-elementary solutions of type (B).} 
\label{tab:solutions2}
\end{table}
So far we have an answer to this question only when $\lambda$ 
is in subregion \eqref{eqn:square}. 
Put 
\begin{equation} \label{eqn:ai-bi}
\alpha_i := p u_i - a, \qquad \beta_i := q u_i - b \qquad 
(i = 1, \dots, m). 
\end{equation} 
\begin{theorem} \label{thm:integer} 
For any non-elementary solution $\lambda = (p,q,r;a,b;x)$ in 
region \eqref{eqn:square}, in addition to $r \in \mathbb{N} := \mathbb{Z}_{>0}$ 
mentioned in Theorem $\ref{thm:gpf}$, either $(\mathrm{A})$ or $(\mathrm{B})$ 
below must be satisfied:   
\begin{enumerate}
\item[$(\mathrm{A})$] \ $p$, $q \in \mathbb{Z}$; \, $p+q+r$ is even; \, 
$\alpha_i$, $\beta_i \in \mathbb{R} \setminus \mathbb{Z}$ \, $(i = 1, \dots, m)$. 
\item[$(\mathrm{B})$] \ $p$, $q \in 1/2 + \mathbb{Z}$; \, 
$\alpha_i$, $\beta_i \in \mathbb{R} \setminus \mathbb{Z}$; \,   
$\alpha_i - (-1)^{p+q+r} \, \beta_i \in 1/2 + \mathbb{Z}$ \, $(i = 1, \dots, m)$.  
\end{enumerate}
The dilation constant $d$ in the gamma product formula \eqref{eqn:gpf} must be  
\begin{equation} \label{eqn:d}
d = \dfrac{r^r}{\sqrt{p^p \cdot q^q \cdot (r-p)^{r-p} \cdot (r-q)^{r-q} \cdot 
x^r \cdot (1-x)^{p+q-r}}}.   
\end{equation}
\end{theorem}
\begin{remark} \label{rem:integer} 
In Theorem \ref{thm:integer}, if $\lambda = (p,q,r;a,b;x)$ is a solution 
of type (B) to Problem \ref{prob:ocf}, then $\hat{\lambda} := (2p,2q,2r;a,b;x)$ 
becomes a solution of type (A) to the problem, with the corresponding 
rational function $\hat{R}(w) = R(2w) \cdot R(2w+1)$. 
Indeed, since $f(w;\hat{\lambda}) = f(2w;\lambda)$, 
\[
\dfrac{f(w+1;\hat{\lambda})}{f(w;\hat{\lambda})} = 
\dfrac{f(2w+2;\lambda)}{f(2w;\lambda)} = 
\dfrac{f(2w+1;\lambda)}{f(2w;\lambda)} \cdot 
\dfrac{f(2w+2;\lambda)}{f(2w+1;\lambda)} = 
R(2w) \cdot R(2w+1) \in \mathbb{C}(w). 
\]
We call $\hat{\lambda}$ the {\sl duplication} of $\lambda$. 
Table \ref{tab:solutions2} exhibits two solutions of type (B), 
the duplications of which are just those solutions with 
$(p,q,r) = (1,1,6)$ in Table \ref{tab:solutions}. 
\end{remark}
\par
The proof of Theorem \ref{thm:integer} relies on an asymptotic analysis 
(details of which are developed in \S\ref{sec:a-r}) and a key lemma 
below (Lemma \ref{lem:sin-sin}).  
For a real number $t \in \mathbb{R}$, let $\{t\} := t-[t] \in [0, \, 1)$ 
denote the fractional part of $t$, where $[t] \in \mathbb{Z}$ stands for the 
largest integer not exceeding $t$. 
\begin{figure}[t]
\begin{minipage}{0.45\hsize}
\vspace{-4mm}
\begin{center}
\includegraphics[width=60mm,clip]{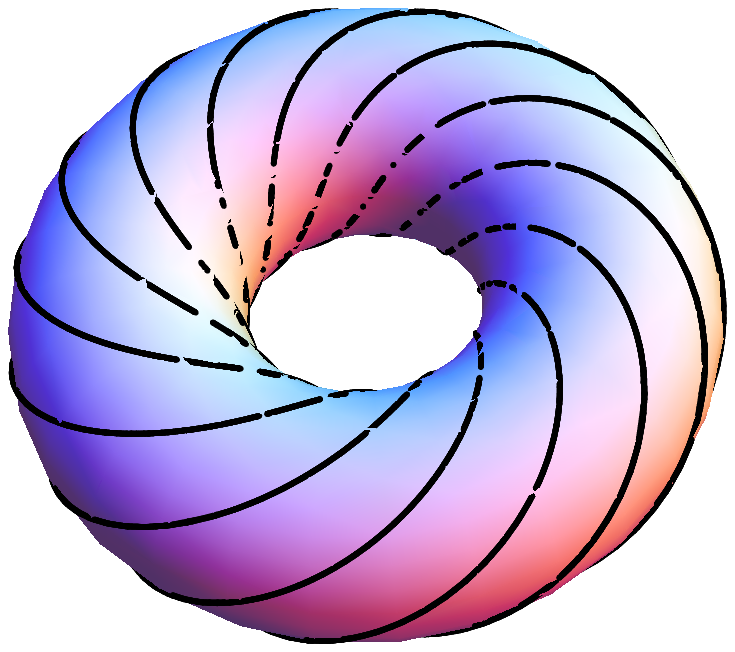} 
\end{center}
\vspace{-11mm} 
\caption{A torus line.}
\label{fig:torus}
\end{minipage}
\begin{minipage}{0.45\hsize}
\begin{center}
\includegraphics[width=45mm,clip]{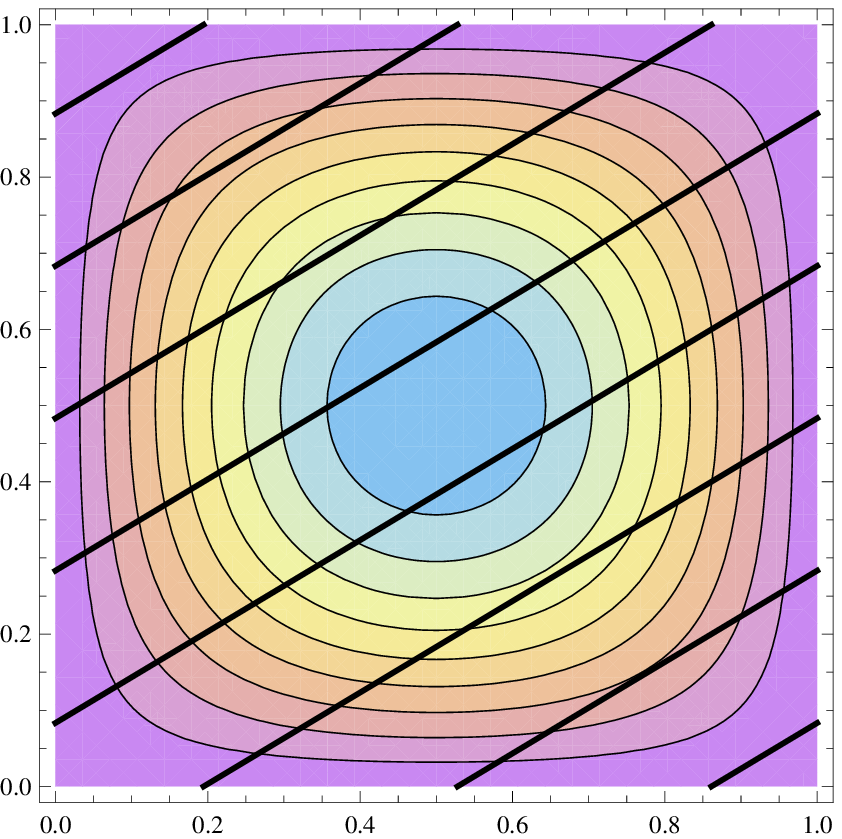}  
\end{center}
\vspace{-3mm} 
\caption{Level curves of $\sin\pi\{x\} \cdot\sin\pi\{y\}$.}
\label{fig:sin-sin-line}
\end{minipage}
\end{figure}
\begin{lemma}[Sine-Sine] \label{lem:sin-sin} 
Let $r \in \mathbb{N}$, $p$, $q \in \mathbb{R}^{\times}$ and 
$\alpha$, $\beta \in \mathbb{R}$. 
Suppose that   
\begin{enumerate}
\item there exist positive constants $C > 0$ and $h > 0$ such that 
\begin{equation} \label{eqn:sin-sin} 
C \cdot h^j \cdot \sin\pi\{p j + \alpha \}\cdot \sin\pi\{q j + \beta \} 
\to 1 \quad \mbox{as} \quad \mathbb{N} \ni j \to \infty,     
\end{equation}
\item the parities of integers $[p j + \alpha]+[q j + \beta] + r j$ 
are stationary as $\mathbb{N} \ni j \to \infty$.   
\end{enumerate}
Then $h = 1$ and either condition $(\mathrm{A})$ or $(\mathrm{B})$ below 
must be satisfied. 
\begin{enumerate}
\item[$(\mathrm{A})$] \ $p$, $q \in \mathbb{Z}$; \, $p+q+r$ is even; \, 
$\alpha$, $\beta \in \mathbb{R} \setminus \mathbb{Z}$.  
\item[$(\mathrm{B})$] \ $p$, $q \in 1/2 + \mathbb{Z}$; \, 
$\alpha$, $\beta \in \mathbb{R} \setminus \mathbb{Z}$, \,   
$\alpha - (-1)^{p+q+r} \, \beta \in 1/2 + \mathbb{Z}$.   
\end{enumerate} 
\end{lemma}
\par
It turns out that condition \eqref{eqn:sin-sin} is equivalent 
to a seemingly more restrictive condition    
\begin{equation} \label{eqn:sin-sin2} 
C \cdot \sin\pi\{p j + \alpha \}\cdot \sin\pi\{q j + \beta \} = 1 \qquad 
\mbox{for all} \quad j \in \mathbb{Z}.   
\end{equation}
The proof of Lemma \ref{lem:sin-sin} is divided into two parts. 
In the first part (in \S\ref{sec:kronecker}), Kronecker's theorem on 
Diophantine approximations is used to show the rationality of 
$p$ and $q$; see Figure \ref{fig:torus} to get a feel for the 
discussions around here. 
In the second part (in \S\ref{sec:level} and \S\ref{sec:parity}), 
a certain geometry and arithmetic (motivated by condition 
\eqref{eqn:sin-sin2}) for the level curves of the function 
$\sin\pi\{x\} \cdot\sin\pi\{y\}$ on the torus $(\mathbb{R}/\mathbb{Z}) 
\times (\mathbb{R}/\mathbb{Z})$ (see Figure \ref{fig:sin-sin-line}) is 
used to reduce the possibilities of $(p, q; \alpha, \beta)$ into 
two types (A) and (B).
Theorem \ref{thm:integer} will be established at the end 
of \S\ref{sec:parity} (see Proposition \ref{prop:parity2}).  
\subsection{Coming from Contiguous Relations} \label{subsec:contiguous}
As mentioned in \S\ref{sec:intro} we can find (partial) solutions 
to Problem \ref{prob:ocf} by the method of contiguous relations and we may 
ask how often this class of solutions occur among all solutions to the problem 
(see Problem \ref{prob:contiguous1}).   
So far it is only in region  \eqref{eqn:square} that we have an answer to 
this question.  
Recall that in this region Problems \ref{prob:ocf} and \ref{prob:gpf} 
are equivalent by Theorem \ref{thm:gpf}. 
\begin{theorem} \label{thm:contiguous1} 
All non-elementary solutions of type $(\mathrm{A})$ to Problem $\ref{prob:ocf}$ 
in region \eqref{eqn:square} come from contiguous relations, where  
solution types are those mentioned in Theorem $\ref{thm:integer}$. 
\end{theorem}
\par
If $\lambda$ is a non-elementary solution of type (B), then by 
Remark \ref{rem:integer} its duplication $\hat{\lambda}$ is a 
non-elementary solution of type (A) and hence comes from contiguous 
relations by Theorem \ref{thm:contiguous1}. 
In this sense we can say that all non-elementary solutions in 
region \eqref{eqn:square} {\sl essentially} come from contiguous relations. 
Theorem \ref{thm:contiguous1} will be proved in Proposition \ref{prop:ocr}.  
\begin{problem} \label{prob:contiguous2} 
Let $\lambda = (p,q,r;a,b;x)$ be a parameter in region 
\eqref{eqn:square}, where $p$, $q$ and $r$ are a priori supposed   
to be integers. 
Find a necessary and sufficient condition in {\sl algebraic} terms  
in order for $\lambda$ to be a solution to Problem \ref{prob:ocf} that 
comes from contiguous relations.  
\end{problem}
\par
Reduced by symmetries \eqref{eqn:sym} this problem may be discussed in 
region \eqref{eqn:apollo}. 
Given a nonnegative integer $k$, let $\langle \, \varphi(z) \, \rangle_k 
:= \sum_{j=0}^k c_j z^j$ denote the {\sl truncation} at degree $k$ of 
a power series $\varphi(z) = \sum_{j=0}^{\infty} c_j z^j$.   
In what follows a truncation will always be taken with respect to 
variable $z$. 
We introduce a ``truncated hypergeometric product" defined by   
\begin{equation} \label{eqn:truncate}
\begin{split}
\Phi(w;\lambda) 
&:= (r w)_{r-1} \, \langle \, {}_2F_1(\mbox{\boldmath $\alpha$}^*(w); z) \cdot 
{}_2F_1(\mbox{\boldmath $v$}-\mbox{\boldmath $\alpha$}^*(w+1);z) \, \rangle_{r-q-1}\big|_{z = x} \\
&\phantom{:}= (r w)_{r-1} \, \langle \, (1-z)^{r-p-q} \cdot 
{}_2F_1(\mbox{\boldmath $\alpha$}(w); z) \cdot {}_2F_1(\mbox{\boldmath $v$}-\mbox{\boldmath $\alpha$}(w+1);z) \, 
\rangle_{r-q-1}\big|_{z = x},
\end{split}    
\end{equation}
where $\mbox{\boldmath $\alpha$}(w) := (p w+a, \, q w+b; \, r w)$, $\mbox{\boldmath $\alpha$}^*(w) := ((r-p)w-a, \, 
(r-q)w-b; \, r w)$, $\mbox{\boldmath $v$} := (1,1;2)$ and the second equality in definition 
\eqref{eqn:truncate} is due to Euler's transformation \eqref{eqn:transf1}. 
An inspection of definition \eqref{eqn:truncate} shows that 
$\Phi(w;\lambda)$ is a polynomial of $(w;\lambda) = (w;p,q,r;a,b;x)$ 
over $\mathbb{Q}$ with degree at most $2r-q-2$ in $w$. 
It is not immediate from definition \eqref{eqn:truncate} but can be 
seen that $\Phi(w;\lambda)$ is more strictly of degree at most $r-1$ 
in $w$ (see Lemma \ref{lem:ebisu2}). 
So we can write  
\[
\Phi(w;\lambda) = \sum_{k=0}^{r-1} \Phi_k(\lambda) w^k. 
\]   
\begin{theorem} \label{thm:contiguous2} 
A parameter $\lambda = (p,q,r;a,b;x)$ in region \eqref{eqn:apollo} 
with $p,q,r \in \mathbb{Z}$ is a solution to Problem $\ref{prob:ocf}$ that 
comes from contiguous relations if and only if $\Phi(w;\lambda)$ 
vanishes identically as a polynomial of $w$, that is, if and only if 
$\lambda$ is a simultaneous root of the following $r$ equations:   
\begin{equation} \label{eqn:Phi}
\Phi_k(\lambda) = 0 \qquad (k = 0, 1, \dots, r-1).  
\end{equation}
\end{theorem}
\par
If a triple $(p,q,r)$ is given then condition \eqref{eqn:Phi} yields  
an {\sl overdetermined system of algebraic equations over} $\mathbb{Q}$ for  
an unknown $(a,b;x)$. 
For a various value of $(p,q,r)$, ask if system \eqref{eqn:Phi} 
admits at least one root and, if so, solve it to obtain an example  
of solution to Problem \ref{prob:contiguous1}.  
Actually all solutions in Table \ref{tab:solutions} were found 
in this manner with the aid of computer, where the use of Theorem 
\ref{thm:contiguous3} below was also helpful. 
Theorem \ref{thm:contiguous2} will be proved in 
Proposition \ref{prop:ebisu}. 
It will turn out that \eqref{eqn:Phi} leads to an algebraic system 
involving certain {\sl terminating} hypergeometric summations, 
which is in some sense more explicit than \eqref{eqn:Phi} itself 
(see Proposition \ref{prop:terminate1}). 
\par
In view of Theorem \ref{thm:contiguous1}, the method of contiguous 
relations brings us a further understanding of non-elementary 
solutions of type (A) in region \eqref{eqn:apollo} as in 
Theorem \ref{thm:contiguous3} below.  
To state it, let $p$, $q$ and $r$ be integers such that 
$0 < q \le p < r$ and $p + q \le r$, and set 
\begin{equation} \label{eqn:Delta} 
\varDelta = \varDelta(z) := (p-q)^2 z^2 - 2 \{(p+q)r-2 p q\} z + r^2. 
\end{equation}
Note that $\varDelta(0) = r^2 > 0$ and $\varDelta(1) = (r-p-q)^2 \ge 0$. 
If $p \neq q$ then $\varDelta(z)$ is a quadratic polynomial with 
axis of symmetry $z = 1 + (p-q)^{-2} \{p(r-p)+q(r-q)\} > 1$.  
If $p = q$ then $\varDelta(z)$ is a linear polynomial with slope 
$-4 p(r-p) < 0$. 
In either case, $\varDelta = \varDelta(z)$ is strictly decreasing 
and positive in $0 \le z < 1$, where we take the 
branch of $\sqrt{\varDelta}$ so that $\sqrt{\varDelta} > 0$.  
Put 
\begin{equation} \label{eqn:Cpm}
Z_{\pm}(z) := 
\big\{r+(p-q)z \pm \sqrt{\varDelta} \big\}^p 
\big\{r-(p-q)z \pm \sqrt{\varDelta} \big\}^q
\big\{(2r-p-q)z-r \mp \sqrt{\varDelta} \big\}^{r-p-q}.   
\end{equation}
Then there exist polynomials $X(z)$ and $Y(z)$ with integer 
coefficients such that 
\begin{equation} \label{eqn:CpmAB}
Z_{\pm}(z) = X(z) \pm Y(z) \, \sqrt{\varDelta}.
\end{equation}
Observe that $Z_{+}(0) = (-1)^{r-p-q}\, (2r)^r$ and $Z_{-}(0) = 0$;  
thus $Y(0) = (-1)^{r-p-q}(2r)^{r-1} \neq 0$ and so the polynomial 
$Y(z)$ is nonzero. 
It is not hard to see that the degree of $Y(z)$ is at most $r-1$ if 
$p \neq q$ and at most $r-p-1$ if $p = q$. 
The following will be shown in Lemma \ref{lem:roots}: \\[2mm]      
\quad $\bullet$ If $r-p-q = 0$ then $Y(z)$ has no root in $0 \le z < 1$. \\[1mm]   
\quad $\bullet$ If $r-p-q$ is a positive even integer then 
$Y(z)$ has at least one root in $0 < z < 1$.  
\par\vspace{1mm}
To state our result we introduce another truncated hypergeometric 
product defined by 
\begin{equation} \label{eqn:P}
\begin{split}
P(w) &:= (r w)_r \, \langle \, {}_2F_1(\mbox{\boldmath $\alpha$}^*(w); z) \cdot 
{}_2F_1(\mbox{\boldmath $1$}-\mbox{\boldmath $\alpha$}^*(w+1);z) \, \rangle_{r-q-1}\big|_{z = x} \\
&\phantom{:}= (r w)_r \, \langle \, (1-z)^{r-p-q-1} \cdot {}_2F_1(\mbox{\boldmath $\alpha$}(w); z)  
\cdot {}_2F_1(\mbox{\boldmath $e$}_3-\mbox{\boldmath $\alpha$}(w+1);z) \, \rangle_{r-q-1}\big|_{z = x},  
\end{split}      
\end{equation}
where $\mbox{\boldmath $\alpha$}(w)$ and $\mbox{\boldmath $\alpha$}^*(w)$ 
are the same as in definition \eqref{eqn:truncate}, $\mbox{\boldmath $1$} := (1,1;1)$, 
$\mbox{\boldmath $e$}_3 := (0,0;1)$ and the second equality in \eqref{eqn:P} 
follows from Euler's transformation \eqref{eqn:transf1}. 
An inspection of definition \eqref{eqn:P} shows that $P(w)$ is a 
polynomial of degree at most $2r-q-1$ in $w$.   
It turns out that if $r-p-q \ge 1$ then $P(w)$ is more strictly   
of degree at most $r$ in $w$ (see Lemma \ref{lem:ebisu2}).  
\par
Regarding the number $x$ and the rational function 
$R(w)$ in formula \eqref{eqn:ocf} we have the following. 
\begin{theorem} \label{thm:contiguous3}
For any non-elementary solution $\lambda = (p,q,r;a,b;x)$ of type 
$(\mathrm{A})$ to Problem $\ref{prob:ocf}$ in region \eqref{eqn:apollo}, 
the following statements must be true: 
\begin{enumerate}
\item The number $x$ is a root of $Y(z)$ in the interval $0 < z < 1$. 
In particular $x$ is an algebraic number of degree at most $r-1$  
if $p \neq q$ and at most $r-p-1$ if $p = q$, respectively.  
\item The integer $r-p-q$ is positive and even.            
\item $P(w)$ is exactly of degree $r$ and the rational function 
$R(w)$ in formula \eqref{eqn:ocf} is given by   
\begin{equation} \label{eqn:R(w)}
R(w) = (1-x)^{r-p-q-1} \cdot \dfrac{(r w)_r}{P(w)}.     
\end{equation}
Moreover, in $\mathbb{C}[w]$ the polynomial $P(w)$ admits 
a division relation 
\begin{equation} \label{eqn:division}
P(w)\, \big|\, (p w+a+1)_{p-1}(q w + b+1)_{q-1}((r-p)w-a)_{r-p}((r-q)w-b)_{r-q}. 
\end{equation}
\end{enumerate}
\end{theorem}
\par
Assertions (1) and (2) of Theorem \ref{thm:contiguous3} will be 
proved in Proposition \ref{prop:roots}, while assertion (3) will 
be established in Proposition \ref{prop:ebisu}, respectively.  
\begin{remark} \label{rem:contiguous} 
A few comments on Theorem \ref{thm:contiguous3} should be in order 
at this stage.  
\begin{enumerate} 
\item In assertion (1) the condition that $x$ should be a root of 
$Y(z)$ is equivalent to the equation $\Phi_{r-1}(\lambda) = 0$ in 
system \eqref{eqn:Phi} (see Lemma \ref{lem:ebisu2}).  
The degree bound for $x$ there is by no means optimal.  
In fact, for every solution $\lambda$ known to the author, 
$Y(z)$ is reducible and $x$ is either rational or quadratic.  
Since $Y(z)$ depends only on $(p,q;r)$, so does the root $x$ and 
hence the dilation constant $d$ in formula \eqref{eqn:d}, that is, 
$d$ is independent of $(a,b)$. 
\item Exactly $r$ among all the $2r-2$ factors on the right side of 
division relation \eqref{eqn:division} appear as factors of $P(w)$.  
It is yet to be decided which $r$ should be chosen.  
This question seems quite hard in general, but it has something 
to do with certain {\sl terminating} hypergeometric sums (see 
Proposition \ref{prop:terminate2}). 
At least one can say that $P(w)$ contains each of 
\begin{equation} \label{eqn:factor}
\textstyle 
w + \frac{r-p-1-a}{r-p}; \quad w + \frac{r-q-1-b}{r-q}; \quad 
w + \frac{a+p-1}{p} \quad \mbox{(if $p \ge 2$)}; \quad 
w + \frac{b+q-1}{q} \quad \mbox{(if $q \ge 2$)},   
\end{equation}
as a factor (see the end of \S\ref{subsec:terminate}).  
In any case, division relation \eqref{eqn:division} provides us with 
much, though not full, information about the numbers 
$u_1, \dots, u_n$ in formula \eqref{eqn:gpf}. 
In particular they are real numbers (rational numbers if so are $a$ and $b$).     
\item It often occurs in formula \eqref{eqn:R(w)} that the 
numerator $(rw)_r$ and denominator $P(w)$ have some factors in 
common that can be canceled to get a reduced representation.       
Or rather the author knows no example in region \eqref{eqn:square} 
for which such a cancellation does not occur. 
\end{enumerate}
\end{remark}
\section{Stationary Phase Method for Euler's Integral} \label{sec:euler}
In this section $f(w)$ is just the function defined by formula 
\eqref{eqn:f}, that is, it may or may not be a solution to Problem 
\ref{prob:gpf} or \ref{prob:ocf}.    
We consider it in the parameter region \eqref{eqn:cross}, which 
can be reduced to subregion \eqref{eqn:pencil} by 
Lemma \ref{lem:reductions}. 
Rather we may and shall work in the intermediate region 
\eqref{eqn:h-strip}.   
Thus under condition \eqref{eqn:h-strip} we study the asymptotic 
behavior of $f(w)$ as $w \to \infty$ on a right half-plane.  
Euler's integral representation for the hypergeometric function 
allows us to write $f(w) = \psi(w) f_1(w)$, where $\psi(w)$ and 
$f_1(w)$ are given by 
\begin{align} 
\psi(w) &= \frac{\varGamma(r w)}{\varGamma(p w+a) \varGamma((r-p)w-a)}, \label{eqn:psi}  \\  
f_1(w) &= \int_0^1 t^{p w+a-1}(1-t)^{(r-p)w-a-1}(1-x t)^{-q w-b} dt.  
\label{eqn:f1}
\end{align}
The improper integral in formula \eqref{eqn:f1} converges if 
$p \, \mathrm{Re}(w) + a > 0$ and $(r-p) \, \mathrm{Re}(w) -a > 0$. 
Due to assumption \eqref{eqn:h-strip} this condition 
is satisfied on the right half-plane $\mathrm{Re}(w) \ge R_1$, if  
\begin{equation} \label{eqn:R1}
R_1 > \max\{-a/p, \, a/(r-p) \} \, (\ge 0). 
\end{equation}     
\par
The gamma factor $\psi(w)$ can be estimated by Stirling's formula, 
which states that $\varGamma(t) \sim \sqrt{2\pi} \, e^{-t} \, 
t^{t-1/2}$ as $t \to \infty$ uniformly on every proper subsector of 
the sector $|\arg(t)| < \pi$, where $* \sim **$ indicates that the 
ratio of $*$ and $**$ tends to $1$ as $t \to \infty$ in the region 
considered.  
It is convenient to note a slightly generalized version of 
Stirling's formula: for any $\alpha > 0$ and $\beta \in \mathbb{C}$, 
\begin{equation} \label{eqn:stirling}
\varGamma(\alpha t + \beta) \sim \sqrt{2\pi} \, \alpha^{\beta-1/2} \,  
(\alpha/e)^{\alpha t} \, t^{\alpha t + \beta -1/2} \qquad 
\mbox{as} \quad t \to \infty,      
\end{equation}  
which is valid on the same sector as above and is easily derived  
from the original formula.   
\begin{lemma} \label{lem:gamma1} 
The function $\psi(w)$ is holomorphic and admits a uniform estimate  
\[
\psi(w) \sim  A_1 \cdot B_1^w \cdot w^{1/2},  
\]
on the right half-plane $\mathrm{Re}(w) > 0$, where $A_1$ and $B_1$ are given by  
\[ 
A_1 = \frac{1}{\sqrt{2\pi}} \cdot 
\frac{(r-p)^{a+1/2}}{p^{a-1/2} \, r^{1/2}}, 
\qquad 
B_1 = \frac{r^r}{p^p \, (r-p)^{r-p}}. 
\]
\end{lemma}
{\it Proof}. 
The poles of $\psi(w)$ are contained in the arithmetic progression 
$\{ -j/r \}_{j = 0}^{\infty}$ and so $\psi(w)$ is 
holomorphic on $\mathrm{Re}(w) > 0$.    
By Stirling's formula \eqref{eqn:stirling} we have 
\begin{align*}
\psi(w) &\sim \frac{\sqrt{2\pi} \, r^{-1/2} \, (r/e)^{r w} \, w^{r w-1/2}}{
\sqrt{2\pi} \, p^{a-1/2} \, (p/e)^{p w} \, w^{p w+a-1/2} \cdot 
\sqrt{2\pi} \, (r-p)^{-a-1/2} \, ((r-p)/e)^{(r-p)w} \, w^{(r-p)w-a-1/2}} \\
&=  \frac{1}{\sqrt{2\pi}} \cdot   
\frac{(r-p)^{a+1/2}}{p^{a-1/2} \, r^{1/2}} \cdot
\left(\frac{r^r}{p^p \, (r-p)^{r-p}} \right)^w \cdot w^{1/2} 
= A_1 \cdot B_1^w \cdot w^{1/2},  
\end{align*}
as $w \to \infty$ uniformly on $\mathrm{Re}(w) > 0$. 
This proves the lemma.  \hfill $\Box$ \par\medskip 
The integral in formula \eqref{eqn:f1} can be written 
\begin{equation} \label{eqn:s-phase}
f_1(w) = \int_0^1 \Phi(t)^w \, \eta(t) \, dt = 
\int_0^1 e^{-w \phi(t)} \, \eta(t) \, dt,   
\end{equation}
where $\Phi(t)$, $\phi(t)$ and $\eta(t)$ are defined by 
\begin{equation} \label{eqn:Phi-eta} 
\begin{split}
\Phi(t) &= t^p (1-t)^{r-p} (1-xt)^{-q}, \qquad 
\phi(t) = - \log \Phi(t), \\ 
\eta(t) &= t^{a-1} (1-t)^{-a-1} (1-x t)^{-b}.
\end{split}   
\end{equation}
\par
We apply the stationary phase method to evaluate 
integral \eqref{eqn:s-phase}. 
Observe that 
\[
\phi'(t) = \frac{\phi_1(t)}{t(1-t)(1-x t)}, 
\qquad \phi_1(t) := -(r-q) x t^2 + \{(p-q) x + r \} t - p, 
\]
where $\phi_1(t)$ is a concave quadratic function thanks to assumption 
\eqref{eqn:h-strip}. 
Here recall that the trivial case $x = 0$ is excluded.   
Since $\phi_1(0) = - p < 0$ and $\phi_1(1) = (r-p)(1-x) > 0$, 
there is a unique root $0 < t_0 < 1$ of the quadratic 
equation $\phi_1(t) = 0$. 
Note that $\phi_1'(t_0) > 0$ and hence  
\[
\phi''(t_0) =  \frac{\phi_1'(t_0)}{t_0(1-t_0)(1-x t_0)} 
= \frac{-2(r-q) x t_0 +(p-q) x + r}{t_0(1-t_0)(1-x t_0)} > 0,   
\]
because $t_0$ lies strictly to the left of the axis of symmetry 
for the parabola $\phi_1(t)$.         
\begin{lemma} \label{lem:asympt} 
The function $f_1(w)$ is holomorphic and admits a uniform estimate  
\begin{equation} \label{eqn:asympt1}
f_1(w) \sim \sqrt{2\pi} \, \frac{\eta(t_0)}{\sqrt{\phi''(t_0)}} \,   
\Phi(t_0)^w \, w^{-1/2} \qquad \mbox{as} \quad w \to \infty,         
\end{equation}
on the right half-plane $\mathrm{Re}(w) \ge R_1$, where $R_1$ is any 
number satisfying condition \eqref{eqn:R1}.   
\end{lemma}
{\it Proof}. 
The function $f_1(w)$ is holomorphic on $\mathrm{Re}(w) \ge R_1$ by the 
convergence condition for the improper integral \eqref{eqn:f1} 
mentioned above. 
Asymptotic formula \eqref{eqn:asympt1} is obtained by the 
standard stationary phase method, so only an outline of its  
derivation will be included below.   
Suppose that $\arg w = 0$ for simplicity. 
Then the path of integration is just the real interval $0 < t < 1$ 
as taken in formula \eqref{eqn:s-phase}, where the phase function 
$\phi(t)$ attains its minimum at $t = t_0$ so that the vicinity 
of this point has the greatest contribution to the value of 
integral \eqref{eqn:s-phase}. 
Observing that 
\[
\phi(t) = \phi(t_0) + \frac{1}{2} \phi''(t_0) (t-t_0)^2 + O((t-t_0)^3), 
\quad \eta(t) = \eta(t_0) + O(t-t_0),  
\quad \mbox{as} \quad t \to t_0,  
\]
we have for any sufficiently small positive number $\varepsilon > 0$, 
\begin{align*}
f_1(w) &\sim \int_{t_0-\varepsilon}^{t_0+\varepsilon} 
e^{-w \{\phi(t_0) + \frac{1}{2} \phi''(t_0) (t-t_0)^2 \}} \, \eta(t) \, dt 
\sim \Phi(t_0)^w \int_{t_0-\varepsilon}^{t_0+\varepsilon} 
e^{- \frac{1}{2} w \phi''(t_0) (t-t_0)^2 } \, \eta(t_0) \, dt \\
&\sim \eta(t_0) \, \Phi(t_0)^w \int_{-\infty}^{\infty} 
e^{- \frac{1}{2} w \phi''(t_0) t^2 } \, dt = 
\dfrac{\eta(t_0)}{\sqrt{\frac{1}{2} w \phi''(t_0)}} \, \Phi(t_0)^w 
\int_{-\infty}^{\infty} e^{-t^2} dt,  
\end{align*}
from which formula \eqref{eqn:asympt1} follows, where we made   
a change of variable $t \sqrt{\frac{1}{2} w \phi''(t_0)} \mapsto t$ 
to obtain the last equality. 
\begin{figure}[t]
\begin{center}
\includegraphics[width=60mm,clip]{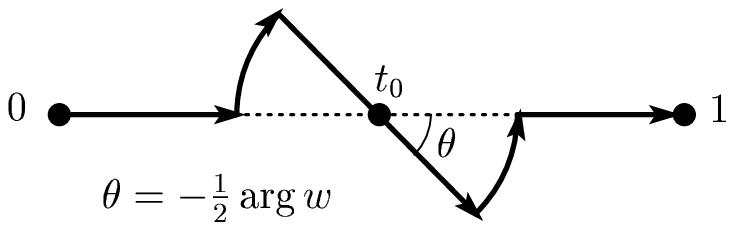}
\end{center}
\caption{A path of steepest descent in the $t$-plane.} 
\label{fig:s-phase}
\end{figure}
This argument carries over for a general complex variable $w$ on the 
right half-plane $\mathrm{Re}(w) \ge R_1$ if the path of integration is 
deformed as in Figure \ref{fig:s-phase}. \hfill $\Box$ \par\medskip
Lemmas \ref{lem:gamma1} and \ref{lem:asympt} are put together to 
yield the following. 
\begin{proposition} \label{prop:asympt}
The function $f(w)$ is holomorphic and admits a uniform estimate 
\begin{equation} \label{eqn:asympt2} 
f(w) \sim A \cdot B^{w} \qquad \mbox{as} \quad w \to \infty, 
\end{equation}
on the right half-plane $\mathrm{Re}(w) \ge R_1$, where $A$ and $B$ are given by 
\begin{equation} \label{eqn:AB} 
A = \frac{(r-p)^{a+1/2}}{p^{a-1/2} \, r^{1/2}} \cdot 
\frac{\eta(t_0)}{\sqrt{\phi''(t_0)}}, 
\qquad 
B = \frac{r^r}{p^p \, (r-p)^{r-p}} \cdot \Phi(t_0). 
\end{equation}
\end{proposition}
{\it Proof}. 
This proposition follows immediately from Lemmas \ref{lem:gamma1} 
and \ref{lem:asympt}. \hfill $\Box$ 
\section{Poles and Their Residues} \label{sec:p-r}
Also in this section $f(w)$ is just the function defined by 
formula \eqref{eqn:f}, which may or may not be a 
solution to Problem \ref{prob:gpf} or \ref{prob:ocf}, while 
condition \eqref{eqn:h-strip} is retained.     
We discuss the pole structure of the function $f(w)$. 
Any pole of $f(w)$ is simple and must lie in the arithmetic progression 
\begin{equation} \label{eqn:a-p0}
W := \{ w_j := -j/r \}_{j=0}^{\infty},   
\end{equation}
but $f(w)$ may be holomorphic at some points of \eqref{eqn:a-p0}.  
In order to know whether a given point $w_j$ is actually a pole or not, 
we need to calculate the residue of $f(w)$ at $w = w_j$.   
\begin{lemma} \label{lem:residue} 
The residue of $f(w)$ at $w=w_j$ admits a hypergeometric expression 
\begin{equation} \label{eqn:residue} 
\underset{\scriptstyle w = w_j}{\mathrm{Res}} \, f(w) = C_j \cdot 
{}_2F_1(a_j, \, b_j; \, j+2; \, x) 
\qquad (j = 0,1,2,\dots),  
\end{equation}
where $a_j := pw_j+j+a+1$, $b_j := qw_j+j+b+1$ and 
\begin{equation} \label{eqn:Cj}
C_j := \frac{(-1)^j}{r} \cdot
\frac{(pw_j+a)_{j+1} \, (qw_j+b)_{j+1}}{j! \, (j+1)!} \, x^{j+1}. 
\end{equation}
\end{lemma}
{\it Proof}. Let $j$ and $k$ be nonnegative integers. 
At the point $w = w_j$ the $k$-th summand of the hypergeometric 
series $f(w) = {}_2F_1(pw+a,qw+b;rw;x)$ has residue  
\[
\underset{\scriptstyle w = w_j}{\mathrm{Res}} \, 
\frac{(p w+a)_k (q w+b)_k}{(r w)_k \, k!} x^k = 
\begin{cases} 
0 & (k \le j), \\[2mm]
\dfrac{1}{r} \cdot 
\dfrac{(-1)^j}{j!} \cdot \dfrac{(p w_j+a)_k (q w_j+b)_k}{(k-j-1)! \, k!} 
\, x^k & (k \ge j+1). 
\end{cases} 
\]
Sum of these numbers over $k \ge j+1$ gives the residue of $f(w)$ at 
$w = w_j$.  
Putting $k = i+j+1$,    
\begin{align*}
\underset{\scriptstyle w = w_j}{\mathrm{Res}} \, f(w) 
&= \dfrac{x^{j+1}}{r} \cdot \dfrac{(-1)^j}{j!} 
\sum_{i=0}^{\infty} 
\dfrac{(p w_j+a)_{i+j+1} \, (q w_j+b)_{i+j+1}}{(i+j+1)! \, i!} \, x^i 
\nonumber \\
&= \frac{x^{j+1}}{r} \cdot \frac{(-1)^j}{j!} 
\sum_{i=0}^{\infty} \frac{(p w_j+a)_{j+1}(p w_j+a+j+1)_i \cdot 
(q w_j+b)_{j+1}(q w_j+b+j+1)_i}{(j+1)! \,(j+2)_i \, i!} \, x^i 
\nonumber \\
&= C_j \sum_{i=0}^{\infty} 
\frac{(a_j)_i \, (b_j)_i}{(j+2)_i \, i!} \, x^i = 
C_j \cdot {}_2F_1(a_j, \, b_j; \, j+2; \, x),  
\end{align*}
where $(t)_{i+j+1} = (t)_{j+1} \, (t+j+1)_i$ is used in 
the second equality. 
This proves formula \eqref{eqn:residue}. \hfill $\Box$ \par\medskip
For every sufficiently large integer $j$, Lemma \ref{lem:residue} 
reduces it to an elementary arithmetic to know whether $f(w)$ is 
holomorphic or has a pole at $w = w_j$.   
\begin{lemma} \label{lem:pole1} 
There exists a positive integer $j_0$ such that for any integer 
$j \ge j_0$ the function $f(w)$ is holomorphic at $w = w_j$ if and 
only if either condition \eqref{eqn:0-res1} or \eqref{eqn:0-res2} 
below holds:       
\begin{subequations} \label{eqn:0-res}
\begin{align} 
r(a+i) &= p j  \qquad \mbox{for some} \quad i \in \{0, \dots, j\},  
\label{eqn:0-res1} \\  
r(b+i) &= q j  \qquad \mbox{for some} \quad i \in \{0, \dots, j\}. 
\label{eqn:0-res2}
\end{align}
\end{subequations}
\end{lemma}
{\it Proof}. Observe that $a_j = \{(r-p)j+r(a+1) \}/r$ and 
$b_j = \{(r-q)j+r(b+1) \}/r$, where $r-p$ and $r-q$ are 
positive numbers due to assumption \eqref{eqn:h-strip}.  
Take an integer $j_0$ so that  
\[
j_0 > \max\left\{ - \frac{r(a+1)}{r-p}, \, -\frac{r(b+1)}{r-q}, 
\, 0 \right\}.   
\] 
Then $a_j$ and $b_j$ are positive for every $j \ge j_0$ 
so that $(a_j)_i$ and $(b_j)_i$ are also positive for 
every $i \ge 1$. 
Since $0 \le x < 1$ by assumption \eqref{eqn:h-strip}, 
we have ${}_2F_1(a_j, b_j; j+2;x) \ge {}_2F_1(a_j, b_j; j+2;0) 
= 1$. 
Thus formula \eqref{eqn:residue} tells us that 
$\underset{\scriptstyle w = w_j}{\mathrm{Res}} \, f(w) = 0$, that is, 
$f(w)$ is holomorphic at $w = w_j$ if and only if $C_j = 0$. 
In view of definition \eqref{eqn:Cj} this condition is 
equivalent to 
\[
(pw_j+a)_{j+1} \, (qw_j+b)_{j+1} = 0, 
\]
which in turn holds true exactly when either condition  
\eqref{eqn:0-res1} or \eqref{eqn:0-res2} is satisfied.    
\hfill $\Box$ \par\medskip
Lemma \ref{lem:pole1} poses the problem of finding nonnegative integers 
$j$ with property \eqref{eqn:0-res1}, which will be referred to as 
problem \eqref{eqn:0-res1}; this convention also applies to property 
\eqref{eqn:0-res2}. 
We wish to know whether each problem has infinitely many 
solutions and, if so, how the solutions look like. 
These questions will be answered after the following preliminary lemma.   
\begin{lemma} \label{lem:ap}
Let $\lambda$ and $\mu$ be coprime integers with $1 \le \mu < \lambda$  
and $\nu$ a real number. 
Let $J$ be the set of all nonnegative integers $j$ such that 
$\nu + \lambda i = \mu j$ for some $i \in \{0, \dots, j\}$.    
Then $J$ is nonempty if and only if $\nu$ is an integer. 
If this is the case then $J$ comprises an arithmetic progression 
$J = \{j_* + \lambda k \}_{k=0}^{\infty}$ with 
$j_* \equiv \nu \mu_* \mod \lambda$, where  
$\mu_* \in \mathbb{Z}$ is inverse to $\mu \mod \lambda$.    
\end{lemma}
{\it Proof}. 
If $J$ is nonempty and has an element $j \in J$ with the corresponding 
number $i \in \{0, \dots, j\}$, then evidently $\nu = \mu j - \lambda i$ 
must be an integer. 
Conversely suppose that $\nu$ is an integer. 
Since $\lambda$ and $\mu$ are coprime, there exist integers $i'$ and $j'$ 
such that $\nu = \mu j' - \lambda i'$, that is,  
$\nu + \lambda i' = \mu j'$.  
Consider $i = i' + \mu k$ and $j = j' + \lambda k$ for $k \in \mathbb{Z}$. 
Note that $j-i = j'-i'+(\lambda - \mu) k$. 
Since $\mu$ and $\lambda-\mu$ are positive, one has $0 \le i \le j$ 
and $\nu + \lambda i = \mu j$ for every sufficiently large $k$. 
Thus $j \in J$ and so $J$ is nonempty. 
Supposing that $J$ is nonempty, let $j_*$ be the smallest 
element of $J$ with the corresponding $i_*$. 
It is easy to see the inclusion 
$\{j_* + \lambda k \}_{k=0}^{\infty} \subset J$. 
Conversely, if $j$ is any element of $J$ with the corresponding $i$, 
then taking the difference from the smallest element yields 
$\lambda(i-i_*) = \mu(j-j_*)$. 
Since $\lambda$ and $\mu$ are coprime, there is a nonnegative 
integer $k$ such that $i-i_* = \mu k$ and $j-j_* = \lambda k$ so that 
$j = j_* + \lambda k$ belongs to the arithmetic progression 
$\{j_*+\lambda k\}_{k=0}^{\infty}$.  
It follows from $\nu + \lambda i_* = \mu j_*$ that 
$j_* \equiv (\mu j_*) \mu_* \equiv (\nu + \lambda i_*) \mu_* 
\equiv \nu \mu_* \mod \lambda$. \hfill $\Box$ 
\begin{lemma} \label{lem:inf-sol}
Problem \eqref{eqn:0-res1} has infinitely many solutions if and 
only if  
\begin{equation} \label{eqn:inf-sol1}
p/r = p_1/r_p \in \mathbb{Q}, \quad p_1, \, r_p \in \mathbb{Z}, \quad (p_1, r_p) = 1, 
\quad 1 \le p_1 < r_p; \quad a_1 := r_p \, a \in \mathbb{Z},  
\end{equation}
for some $p_1$ and $r_p$, in which case all solutions to 
\eqref{eqn:0-res1} comprise an arithmetic progression 
\begin{equation} \label{eqn:a-p3}
J_p = \{ j_p+ r_p k \}_{k = 0}^{\infty}, \qquad \mbox{with} \quad 
j_p \equiv a_2 := a_1 p_2 \mod r_p, 
\end{equation}
where $p_2$ is an $(\mbox{any})$ integer such that 
$p_1 p_2 \equiv 1 \mod r_p$.  
On the other hand, problem \eqref{eqn:0-res2} has infinitely many 
solutions if and only if either condition \eqref{eqn:elementary1} is 
satisfied, in which case all solutions to \eqref{eqn:0-res2} comprise 
an arithmetic progression $J_q = \{-b + k\}_{k=0}^{\infty};$ or 
otherwise, 
\begin{equation} \label{eqn:inf-sol2} 
q/r = q_1/r_q \in \mathbb{Q}, \quad q_1, \, r_q \in \mathbb{Z}, \quad (q_1, r_q) = 1, 
\quad 1 \le q_1 < r_q; \quad b_1 := r_q \, b \in \mathbb{Z},  
\end{equation}
for some $q_1$ and $r_q$, in which case all solutions to 
\eqref{eqn:0-res2} comprise an arithmetic progression 
\begin{equation} \label{eqn:a-p4}
J_q = \{ j_q + r_q k \}_{k=0}^{\infty}, \qquad \mbox{with} \quad 
j_q \equiv b_2 :=  b_1 q_2 \mod r_q,   
\end{equation} 
where $q_2$ is an $(\mbox{any})$ integer such that 
$q_1 q_2 \equiv 1 \mod r_q$.   
\end{lemma}
{\it Proof}. 
If problem \eqref{eqn:0-res1} has infinitely many solutions,  
then of course it has two solutions $j < j'$ with $i$ and $i'$ 
being the corresponding values of $i$.  
Then $r(a+i) = p j$ and $r(a+i') = p j'$, 
whose difference makes $r(i'-i) = p(j'-j)$. 
Thus $p/r = (i'-i)/(j'-j)$ must be a rational number. 
From $0 < p < r$ in assumption \eqref{eqn:h-strip} one has  
$1 \le p_1 < r_p$ in the reduced representation $p_1/r_p$ of the 
rational number $p/r$. 
In terms of $p_1$ and $r_p$, the condition \eqref{eqn:0-res1} is 
equivalent to  
\begin{equation} \label{eqn:0-res3}
a_1 + r_p i = p_1 j  \qquad \mbox{for some} \quad i \in \{0, \dots, j\}.
\end{equation}
Now that $p_1$ and $r_p$ are coprime integers such that $1 \le p_1 < r_p$, 
Lemma \ref{lem:ap} can be applied to problem \eqref{eqn:0-res3} to 
establish the assertion for problem \eqref{eqn:0-res1}. 
\par
We proceed to the assertion for problem \eqref{eqn:0-res2}.  
First we consider the case $q = 0$. 
In this case condition \eqref{eqn:0-res2} becomes 
$b+i = 0$ with $0 \le i \le j$, that is, $b = -i$ is a 
nonpositive integer and $j \ge i = -b$. 
Thus problem \eqref{eqn:0-res2} has infinitely many solutions 
precisely when condition \eqref{eqn:elementary1} is satisfied,  
in which case the integers $j \ge -b$ give all solutions. 
Next we consider the case where $q$ is nonzero.    
The proof is just the same as for problem \eqref{eqn:0-res1}   
except that we have to show $0 < q < r$ if problem \eqref{eqn:0-res2} 
has infinitely many solutions. 
Note that $q < r$ is evident from condition \eqref{eqn:h-strip}.   
To show $0 < q$, suppose the contrary that $q < 0$ and problem 
\eqref{eqn:0-res2} has infinitely many solutions $j \in J$ with 
the corresponding $i \in I$'s.  
Since $q$ and $r$ are nonzero and $r(b+i) = q j$, the 
correspondence $J \to I$, $j \mapsto i$ is a one-to-one mapping.  
The infinite set $J$ contains an infinite subset $J'$ 
such that $j > 0$ for every $j \in J'$. 
The corresponding subset $I' \subset I$ is also infinite. 
For every $i \in I'$ we have $r(b+i) = q j < 0$ and 
so $0 \le i < -b$.   
But this is absurd and hence we must have $0 < q < r$. 
\hfill $\Box$ \par\medskip
To describe how \eqref{eqn:a-p3} and \eqref{eqn:a-p4} intersect 
we require another preliminary lemma. 
\begin{lemma} \label{lem:intersect1} 
Consider two arithmetic progressions 
$J_1 = \{j_1 + \lambda_1 k \}_{k=0}^{\infty}$ and 
$J_2 = \{j_2 + \lambda_2 k \}_{k=0}^{\infty}$, where $\lambda_1$ and  
$\lambda_2$ are positive integers while $j_1$ and $j_2$ are integers. 
Then $J_1$ and $J_2$ intersect if and only if $j_1 \equiv j_2 \mod 
\lambda_3 := \gcd\{\lambda_1, \, \lambda_2\}$. 
In this case the intersection $J_1 \cap J_2$ comprises an 
arithmetic progression $J_3 = \{ j_3 + 
(\lambda_1 \lambda_2/\lambda_3) k \}_{k = 0}^{\infty}$, where 
$j_3 \equiv j_1 \equiv j_2 \mod \lambda_3$.   
\end{lemma} 
{\it Proof}. 
If $J_1$ and $J_2$ have an element in common, say, 
$j_1 + \lambda_1 k_1 = j_2 + \lambda_2 k_2$, then clearly 
$j_1 - j_2 = \lambda_2 k_2 - \lambda_1 k_1 \equiv 0 \mod \lambda_3$. 
Conversely suppose that $j_1 \equiv j_2 \mod \lambda_3$. 
Then there exist integers $k_1'$ and $k_2'$ such that 
$j_1 - j_2 = \lambda_2 k_2' - \lambda_1 k_1'$, that is, 
$j_1 + \lambda_1 k_1' = j_2 + \lambda_2 k_2'$. 
Consider the number $j := j_1 + \lambda_1 k_1 = j_2 + \lambda_2 k_2$,  
where $k_1 := k_1' + (\lambda_2/\lambda_3) l$ and 
$k_2 := k_2' + (\lambda_1/\lambda_3) l$ with $l \in \mathbb{Z}$.  
For every sufficiently large $l$, one has 
$k_1 \ge 0$ and $k_2 \ge 0$ so that $j \in J_1 \cap J_2$ and 
hence $J_1 \cap J_2$ is nonempty. 
Suppose now that $J_1 \cap J_2$ is nonempty and let $j_3 
= j_1 + \lambda_1 k_1' = j_2 + \lambda_2 k_2'$ be its 
smallest element.  
Then it is easy to see the inclusion $\{j_3 + 
(\lambda_1 \lambda_2/\lambda_3) k \}_{k = 0}^{\infty} \subset 
J_1 \cap J_2$. 
Conversely, if $j = j_1 + \lambda_1 k_1 = j_2 + \lambda_2 k_2$ 
is any element of $J_1 \cap J_2$, then its difference from the 
smallest element yields $\lambda_1(k_1-k_1') = \lambda_2(k_2-k_2') \ge 0$,  
or equivalently $(\lambda_1/\lambda_3)(k_1-k_1') = 
(\lambda_2/\lambda_3) (k_2-k_2') \ge 0$. 
Since $\lambda_1/\lambda_3$ and $\lambda_2/\lambda_3$ are coprime, 
there exists a nonnegative integer $k$ such that 
$k_1-k_1' = (\lambda_2/\lambda_3) k$ and $k_2-k_2' = (\lambda_1/\lambda_3) k$. 
Thus $j = j_1 + \lambda_1\{ k_1' + (\lambda_2/\lambda_3) k \} 
= j_3 + (\lambda_1 \lambda_2/\lambda_3) k$. 
This gives the reverse inclusion $J_1 \cap J_2 \subset 
\{j_3 + (\lambda_1 \lambda_2/\lambda_3) k \}_{k = 0}^{\infty}$. 
Finally $j_3 \equiv j_1 \equiv j_2 \mod \lambda_3$ is immediate 
from  $j_3 = j_1 + \lambda_1 k_1' = j_2 + \lambda_2 k_2'$. 
The proof is complete.  \hfill $\Box$ 
\begin{lemma} \label{lem:intersect2}
If conditions \eqref{eqn:inf-sol1} and \eqref{eqn:inf-sol2} 
are satisfied, then $J_p$ and $J_q$ intersect if and only if 
\begin{equation} \label{eqn:intersect1} 
a_2 \equiv b_2 \mod r_{pq} := \gcd\{r_p, \, r_q\},   
\end{equation}
in which case the intersection $J_p \cap J_q$ comprises  
an arithmetic progression  
\begin{equation} \label{eqn:intersect2}
J_{pq} := \{\, j_{pq} + (r_p r_q/r_{pq}) \, k \, \}_{k=0}^{\infty} 
\qquad \mbox{with} \quad 
j_{pq} \equiv a_2 \equiv b_2 \mod r_{pq}.  
\end{equation}  
\end{lemma}
{\it Proof}. 
This lemma is an immediate consequence of Lemmas \ref{lem:inf-sol} 
and \ref{lem:intersect1}. 
\hfill $\Box$ \par\medskip 
For a subset $A \subset \mathbb{R}$ bounded below or above its 
{\sl density} is defined by the limit    
\[
\delta(A) := \lim_{t \to \infty} \frac{1}{t} \, 
\# \{x \in A \,:\, \pm x \le t \},   
\]
where $\#$ denotes the cardinality of a set and plus (resp. minus) 
sign is chosen when $A$ is bounded below (resp. above).  
We consider only those subsets for which density is well defined.  
When density is considered for two or more subsets, they are 
simultaneously bounded below or above.  
Two sets $A$ and $B$ are said to be {\sl commensurable} if they 
share all but a finite number of elements, in which case 
we write $A \circeq B$.   
As basic properties of density we have    
\begin{align*}
\delta(A \cup B) &= \delta(A) + \delta(B) - \delta(A \cap B), \\   
\delta(A') &= |\mu|^{-1} \ \delta(A) \qquad \mbox{for}  \quad 
A' := \mu A + \nu \quad \mbox{with} \quad 
(\mu, \nu) \in \mathbb{R}^{\times} \times \mathbb{R}, \\
\delta(A) &= \delta(B) \qquad \mbox{if} \quad A \circeq B.    
\end{align*}
\par
Let $W_{\mathrm{hol}} \subset W$ be the set of all points 
in \eqref{eqn:a-p0} at which $f(w)$ is holomorphic and 
$J$ the set of all nonnegative integers $j$ that satisfy  
condition \eqref{eqn:0-res1} or \eqref{eqn:0-res2}. 
Lemma \ref{lem:pole1} implies that    
\begin{equation} \label{eqn:Whol}
W_{\mathrm{hol}} \circeq -J/r, \qquad 
\delta(W_{\mathrm{hol}}) = r \, \delta(J). 
\end{equation}  
\begin{lemma} \label{lem:density1} 
If condition \eqref{eqn:elementary1} is satisfied 
then $\delta(J) = 1$.  
Otherwise, $\delta(J)$ is given as in Table $\ref{tab:density}$, 
depending upon whether condition \eqref{eqn:inf-sol1}, 
\eqref{eqn:inf-sol2} or \eqref{eqn:intersect1} is true $(\mathrm{T})$ 
or false $(\mathrm{F})$.    
\begin{table}[t]
\begin{center} 
\begin{tabular}{ccccl}
\hline
            &                             &                             &                               &              \\[-4mm]
case \qquad & \eqref{eqn:inf-sol1} \qquad & \eqref{eqn:inf-sol2} \qquad & \eqref{eqn:intersect1} \qquad &  $\delta(J)$ \\[1mm]
\hline
            &                             &                             &                               &              \\[-4mm]  
I     & F     & F     &   ---    & $0$  \\[1mm]
\hline
            &                             &                             &                               &              \\[-4mm]
II    & F     & T     &   ---    & $1/r_q$  \\[1mm]
\hline
            &                             &                             &                               &              \\[-4mm]
III   & T     & F     &   ---    & $1/r_p$  \\[1mm]
\hline
            &                             &                             &                               &              \\[-4mm]
IV    & T     & T     &  F       & $1/r_p + 1/r_q$ \\[1mm]
\hline
            &                             &                             &                               &              \\[-4mm]
V     & T     & T     &  T       & $1/r_p + 1/r_q - r_{pq}/(r_p r_q)$  \\[1mm]
\hline
\end{tabular}
\end{center} 
\caption{Density $\delta(J)$ in the failure of    
condition \eqref{eqn:elementary1}.} \label{tab:density}
\end{table}
\end{lemma} 
{\it Proof}. 
Note that the whole arithmetic progression \eqref{eqn:a-p0} has 
density $r$. 
If condition \eqref{eqn:elementary1} holds true then our hypergeometric 
sum is terminating so that $f(w)$ is a rational function, clearly 
having only a finite number of poles. 
If $f(w)$ has only a finite number of poles then all but a finite 
number of points in the set \eqref{eqn:a-p0} belong to 
$W_{\mathrm{hol}}$ so that $\delta(W_{\mathrm{hol}}) = r$ and hence 
$\delta(J) = 1$ by formula \eqref{eqn:Whol}.     
Now suppose that condition \eqref{eqn:elementary1} is not satisfied. 
In case I of Table \ref{tab:density}, it is obvious that $J$ is empty 
and $\delta(J) = 0$.   
In case II one has $J = J_q$ and $\delta(J) = \delta(J_q) = 1/r_q$.   
Similarly in case III one has $J = J_p$ and $\delta(J) = 
\delta(J_p) = 1/r_p$. 
In cases IV and V one has $J = J_p \cup J_q$ and 
$\delta(J) = \delta(J_p \cup J_q) = \delta(J_p) + \delta(J_q) - 
\delta(J_p \cap J_q)$, so formulas \eqref{eqn:a-p3}, \eqref{eqn:a-p4} 
and \eqref{eqn:intersect2} yield $\delta(J) =1/r_p + 1/r_q$ in case 
IV and $\delta(J) = 1/r_p + 1/r_q - r_{pq}/(r_p r_q)$ in case V. 
\hfill $\Box$ 
\begin{proposition} \label{prop:f-pole}
The function $f(w)$ has only a finite number of poles if and only 
if either condition \eqref{eqn:elementary1} or \eqref{eqn:elementary2} 
is satisfied. 
If this is the case then $f(w)$ yields an elementary solution to 
Problem $\ref{prob:gpf}$ and hence to Problem $\ref{prob:ocf}$. 
Therefore Theorem $\ref{thm:elementary}$ follows.    
\end{proposition}
{\it Proof}. If condition \eqref{eqn:elementary1} holds then  
our hypergeometric sum is terminating so that $f(w)$ is a rational 
function, evidently having only a finite number of poles and giving 
an elementary solution to Problem \ref{prob:gpf}.  
Suppose that condition \eqref{eqn:elementary1} is not satisfied and 
so we are in one of the five cases in Table \ref{tab:density}. 
As is seen in the proof of Lemma \ref{lem:density1}, $f(w)$ has 
only a finite number of poles only when $\delta(J) = 1$.     
So let us consider when $\delta(J) = 1$ occurs.     
Obviously it cannot occur in cases I, II and III, because 
$r_p \ge 2$ and $r_q \ge 2$. 
A simple check shows that in case IV it is again impossible   
unless $r_p = r_q = 2$. 
This also implies that it is not possible in case V either, 
since the presence of positive term $r_{pq}/(r_p r_q) > 0$ forces 
$\delta(J) < 1$ even when $r_p = r_q = 2$. 
Finally, in case IV with $r_p = r_q = 2$ the validity of conditions 
\eqref{eqn:inf-sol1} and \eqref{eqn:inf-sol2} shows that $p_1 = q_1 = 1$, 
$p = q = r/2$ and $2 a$, $2 b \in \mathbb{Z}$, while the failure of 
condition \eqref{eqn:intersect1} means that 
$2a$ and $2b$ must have distinct parities.  
By symmetry \eqref{eqn:sym0} we may assume that 
$2a$ is even and $2b$ is odd, that is, $a = i$ and $b = j-1/2$ for 
some $i$, $j \in \mathbb{Z}$.  
This leads to condition \eqref{eqn:elementary2}. 
Conversely, if condition \eqref{eqn:elementary2} is satisfied 
then formula \eqref{eqn:dihedral} follows from an analysis of 
Vidunas \cite{Vidunas2}. 
This formula evidently shows that $f(w)$ yields an elementary 
solution to Problem \ref{prob:gpf}.   
\hfill $\Box$ 
\section{Gamma Product Formula} \label{sec:gamma-p}
We continue to work on the parameter region \eqref{eqn:h-strip}. 
Suppose that $f(w)$ is a solution to Problem \ref{prob:ocf} with 
$R(w)$ in condition \eqref{eqn:ocf} being of the form \eqref{eqn:canonical}, 
where representation \eqref{eqn:canonical} may or may not be in a canonical 
form, for example, it may be just the reduced expression of $R(w)$ with 
$S(w) \equiv 1$, while $S(w)$, $d$, $P(w) := (w+u_1)\cdots(w+u_m)$ 
and $Q(w) := (w+v_1)\cdots(w+v_n)$ are supposed to be real, since 
$R(w)$ is real. 
Consider the entire meromorphic function defined by 
\begin{equation} \label{eqn:g}
g(w) := S(w) \cdot d^w \cdot 
\frac{\varGamma(w+u_1) \cdots \varGamma(w+u_m)}{\varGamma(w+v_1) \cdots \varGamma(w+v_n)}.     
\end{equation}
Put $u := u_1+ \cdots+u_m$ and $v := v_1+\cdots+v_n$; they are real 
because $P(w)$ and $Q(w)$ are real.  
\begin{lemma} \label{lem:gamma2} 
There exists a constant $R_2$ such that on the right half-plane 
$\mathrm{Re}(w) \ge R_2$ the function $g(w)$ is holomorphic, nowhere 
vanishing, and admits a uniform estimate 
\[
g(w) \sim S_0 \cdot (2\pi)^{(m-n)/2} \left(d e^{n-m}\right)^w \cdot 
w^{-(m-n)/2+u-v+s_0} \cdot e^{(m-n) w \log w}, 
\]
where $S_0 \in \mathbb{R}^{\times}$ and $s_0 \in \mathbb{Z}$ are determined 
by the condition $S(w) \sim S_0 w^{s_0}$ as $w \to \infty$.   
\end{lemma}
{\it Proof}. 
Take a number $R_2 \in \mathbb{R}$ in such a manner that all the points  
$-u_1, \dots, -u_m; -v_1, \dots, -v_n$ as well as all the 
zeros and poles of $S(w)$ are strictly to the left of the vertical 
line $\mathrm{Re}(w) = R_2$.  
Then it is clear from the locations of its poles and zeros 
that $g(w)$ is holomorphic and non-vanishing on the half-plane 
$\mathrm{Re}(w) \ge R_2$.   
By Stirling's formula \eqref{eqn:stirling}, we have   
\begin{align*}
g(w) &= S(w) \cdot d^w \frac{\prod_{i=1}^m \varGamma(w+u_i)}{\prod_{j=1}^n \varGamma(w+v_j)} 
\sim S_0 w^{s_0} \cdot d^w \frac{\prod_{i=1}^m \sqrt{2\pi} \, e^{-w} \,   
w^{w+u_i-1/2}}{\prod_{j=1}^n \sqrt{2\pi} \, e^{-w} \, w^{w+v_j-1/2}} \\
&= S_0 w^{s_0} \cdot d^w \cdot (2\pi)^{(m-n)/2} \, e^{(n-m)w} \,  
w^{(m-n)(w-1/2) + u -v} \\
&=  S_0 \cdot (2\pi)^{(m-n)/2} \cdot (d e^{n-m})^w \cdot 
w^{-(m-n)/2+u-v+s_0} \cdot w^{(m-n) w} \\
&= S_0 \cdot (2\pi)^{(m-n)/2} \cdot (d e^{n-m})^w \cdot w^{-(m-n)/2+u-v+s_0} 
\cdot e^{(m-n) w \log w}.  
\end{align*}
uniformly on $\mathrm{Re}(w) \ge R_2$.  
This establishes the lemma.  \hfill $\Box$ \par\medskip
Observe that $g(w)$ satisfies the same recurrence relation  
\eqref{eqn:ocf} as the function $f(w)$. 
So it is natural to compare $f(w)$ with $g(w)$ or 
in other words to think of the ratio  
\begin{equation} \label{eqn:h}
h(w) := \frac{f(w)}{g(w)}.  
\end{equation}
It is clear that $h(w)$ is an entire meromorphic function that does 
not vanish identically.    
\begin{lemma} \label{lem:bound2} 
$h(w)$ is an entire holomorphic function which is periodic of 
period one.  
For any $R_3 > \max\{R_1, \, R_2, \, 1 \}$ there exists a constant 
$A_2 > 0$ such that  
\begin{equation} \label{eqn:bound2} 
|h(w)| \le A_2 \cdot K^{\mathrm{Re}(w)} 
\cdot |w|^{(n-m)\{\mathrm{Re}(w)-1/2\}+v-u-s_0} \cdot e^{-(n-m) \arg(w) \cdot \mathrm{Im}(w)}, 
\end{equation}
on $\mathrm{Re}(w) \ge R_3$, where $K := e^{m-n} B/d$ with $B$ being the  
positive constant in Proposition $\ref{prop:asympt}$.   
\end{lemma}
{\it Proof}. 
Since $f(w)$ and $g(w)$ satisfy the same recurrence relation 
\eqref{eqn:ocf}, their ratio $h(w)$ must be a periodic function of 
period one.  
From Proposition \ref{prop:asympt} and Lemma \ref{lem:gamma2} the 
function $h(w)$ has no poles on $\mathrm{Re}(w) \ge R_3$ and so holomorphic 
there.  
The periodicity then implies that $h(w)$ must be holomorphic 
on the entire complex plane. 
In view of $\left|e^{w \log w}\right| = 
e^{\mathrm{Re}(w \log w)} = e^{\mathrm{Re}(w) \cdot \log|w|-\mathrm{Im}(w) \cdot \arg(w)}  
= |w|^{\mathrm{Re}(w)} e^{-\mathrm{Im}(w) \cdot \arg(w)}$, Lemma \ref{lem:gamma2} 
implies that   
\[
|g(w)| \sim |S_0|\cdot (2\pi)^{(m-n)/2} \left(d e^{n-m}\right)^{\mathrm{Re}(w)} 
\cdot |w|^{(m-n)\{\mathrm{Re}(w)-1/2\}+u-v+s_0} \cdot 
e^{-(m-n) \mathrm{Im}(w) \cdot \arg(w)},  
\]
uniformly on $\mathrm{Re}(w) \ge R_3$. 
Since $g(w)$ has no zero there, there is a constant $A_3 > 0$ such that 
\[
|g(w)| \ge A_3 \cdot  \left(d e^{n-m}\right)^{\mathrm{Re}(w)} \cdot 
|w|^{(m-n)\{\mathrm{Re}(w)-1/2\}+u-v+s_0} \cdot e^{-(m-n) \mathrm{Im}(w) \cdot \arg(w)},  
\]
on $\mathrm{Re}(w) \ge R_3$. 
On the other hand, by Proposition \ref{prop:asympt} there exists a 
constant $A_4 > 0$ such that $|f(w)| \le A_4 \cdot B^{\mathrm{Re}(w)}$ on 
$\mathrm{Re}(w) \ge R_3$.  
Thus estimate \eqref{eqn:bound2} holds true with $A_2 := A_4/A_3$.  
\hfill $\Box$ 
\begin{lemma} \label{lem:const} 
We have $m=n$ in formulas \eqref{eqn:canonical} and \eqref{eqn:g}. 
Moreover $h(w)$ is a nonzero constant.  
\end{lemma}
{\it Proof}. 
First we show $m \le n$. 
Suppose the contrary $m > n$. 
Estimate \eqref{eqn:bound2} with real $w$ reads 
$|h(w)| \le A_2 \cdot K^w \cdot w^{-(m-n)(w-1/2)+v-u-s_0}$ for 
every $w \ge R_3$. 
Fix any $w \in \mathbb{R}$ and take a positive integer $k_0$ such that 
$w + k_0 \ge R_3$. 
Since $h(w)$ is periodic of period one, for any integer $k \ge k_0$, 
\begin{align*}
|h(w)| = |h(w+k)| & \le 
A_2 \cdot K^{w+k} \cdot (w+k)^{-(m-n)(w+k-1/2)+v-u-s_0} \\
&= A_2 \cdot K^w \cdot (1+w/k)^{\rho} \cdot 
(1+w/k)^{-(m-n)k} \cdot k^{\rho} \cdot (K/k^{m-n})^k \\
&\sim A_2 \cdot K^w \cdot e^{-(m-n)w} \cdot k^{\rho} \cdot (K/k^{m-n})^k 
\quad \mbox{as} \quad k \to +\infty, 
\end{align*}
where $\rho := -(m-n)(w-1/2)+v-u-s_0$. 
Since we are assuming that $m-n > 0$ there exists an integer 
$k_1 \ge k_0$ such that$0 < K/k_1^{m-n} \le 1/2$. 
Then there exists a constant $A_5 > A_2 \cdot K^w \cdot e^{-(m-n)w}$ 
such that $|h(w)| \le A_5 \cdot k^{\rho} \cdot 2^{-k}$ for every 
$k \ge k_1$.  
Letting $k \to +\infty$ we have $h(w) = 0$ for every 
$w \in \mathbb{R}$. 
By the unicity theorem for holomorphic functions, $h(w)$ must 
vanish identically. 
But this is absurd because $h(w)$ is nontrivial and thus 
we have proved $m \le n$. 
\par 
Next we show that $h(w)$ is a nonzero constant.  
We make use of estimate \eqref{eqn:bound2} on the strip 
$R_3 \le \mathrm{Re}(w) \le R_3+1$, where we recall $R_3 > 1$.   
On this strip $K^{\mathrm{Re}(w)}$ are bounded while 
$|w|^{(n-m)\{\mathrm{Re}(w)-1/2\}+v-u-s_0} \le |w|^{\mu} 
\le A_6 (1+|\mathrm{Im}(w)|^{\mu})$ for some constant $A_6$,  
where $\mu$ is a nonnegative number with 
$\mu \ge (n-m)(R_2-1/2)+v-u-s_0$.  
On the strip, if $|\mathrm{Im}(w)| \ge R_3+1$ then $|\arg(w)| \ge \pi/4$ 
and $\arg(w) \cdot \mathrm{Im}(w) \ge (\pi/4) |\mathrm{Im}(w)|$. 
So there is a constant $A_7$ such that 
\begin{equation} \label{eqn:bound3}
|h(w)| \le A_7 (1+|\mathrm{Im}(w)|^{\mu}) e^{-\pi(n-m)|\mathrm{Im}(w)|/4},  
\end{equation}
holds for any point $w$ on the strip $R_3 \le \mathrm{Re}(w) \le R_3+1$ 
such that $|\mathrm{Im}(w)| \ge R_3+1$.   
Estimate \eqref{eqn:bound3} remains true on the entire strip if 
$A_7$ is chosen sufficiently large. 
Moreover this estimate extends to the entire complex plane, since both 
sides of it are periodic functions of period one. 
In particular, in view of $m \le n$, estimate \eqref{eqn:bound3} yields  
$|h(w)| \le A_7 \cdot (1+|\mathrm{Im}(w)|^{\mu}) \le A_7 \cdot 
(1+|w|^{\mu})$ for every $w \in \mathbb{C}$. 
Liouville's theorem then implies that $h(w)$ must be a polynomial. 
But the fundamental theorem of algebra tells us that a 
polynomial can be a periodic function only when it 
is a constant.  
Hence $h(w)$ must be a constant, which is nonzero as 
$h(w)$ is nontrivial. 
Finally we show that $m = n$. 
We already know that $m \le n$. 
If $m < n$ then the right-hand side of estimate \eqref{eqn:bound3} 
would tend to zero as $|\mathrm{Im}(w)| \to \infty$. 
But this contradicts the fact that $h(w)$ is a nonzero constant. 
Thus we must have $m =n$. 
The proof is complete. \hfill $\Box$ \par\medskip              
Recall that we can multiply the rational function $S(w)$ by any nonzero 
constant without changing the form of expression \eqref{eqn:canonical}.  
Thus after multiplying $S(w)$ by a suitable constant if necessary, 
we may conclude that $h(w) \equiv 1$ in Lemma \ref{lem:const} and 
so definitions \eqref{eqn:g} and \eqref{eqn:h} yield   
\begin{equation} \label{eqn:f=g}
f(w) \equiv g(w) := S(w) \cdot d^w \cdot 
\frac{\varGamma(w+u_1) \cdots \varGamma(w+u_m)}{\varGamma(w+v_1) \cdots \varGamma(w+v_m)}.   
\end{equation}
\begin{proposition} \label{prop:f=g} 
If $f(w)$ is a solution to Problem $\ref{prob:ocf}$ in region 
\eqref{eqn:h-strip}, then we have an identity \eqref{eqn:f=g} with 
$m = n$, $S_0 = A$, $d = B$ and $v = u + s_0$ in formula 
\eqref{eqn:canonical}, where $S_0 \in \mathbb{R}^{\times}$ and $s_0 \in \mathbb{Z}$ are 
the numbers defined by the asymptotic condition $S(w) \sim S_0 w^{s_0}$ 
as $w \to \infty$; $A$ and $B$ are the constants defined 
in Proposition $\ref{prop:asympt}$; and $u := u_1+\cdots+u_m$ and 
$v:= v_1+\cdots+v_n$.   
\end{proposition}
{\it Proof}. 
We have only to evaluate $S_0$, $d$ and $v$. 
Now that we have $m = n$ by Lemma \ref{lem:const}, the asymptotic 
formula in Lemma \ref{lem:gamma2} reads $g(w) \sim S_0 \cdot 
d^w \cdot w^{u-v+s_0}$ as $w \to \infty$ on the half-plane 
$\mathrm{Re}(w) \ge R_2$. 
Compare this with the asymptotic formula $f(w) \sim 
A \cdot B^w$ in Proposition \ref{prop:asympt}. 
The identity $f(w) \equiv g(w)$ then implies that 
$S_0 = A$, $d = B$ and $u-v+s_0 = 0$.  
\hfill $\Box$ 
\section{Rational Functions in Canonical Form} \label{sec:rational}
We make a general discussion about rational functions in order 
to put representation \eqref{eqn:canonical} in a canonical form.  
Given a rational function $R(w) \in \mathbb{C}(w)$, consider an 
expression of 
the form  
\begin{equation} \label{eqn:rational1} 
R(w) = \frac{S(w+1)}{S(w)} \cdot d \cdot \frac{P(w)}{Q(w)},  
\end{equation}
where $S(w)$ is a rational function, $d$ is a nonzero constant, and 
$P(w)$ and $Q(w)$ are monic polynomials. 
$P(w)$ and $Q(w)$ are said to be {\sl strongly coprime} if $P(w)$ and 
$Q(w+j)$ are coprime over $\mathbb{C}$ for every integer $j$, in which 
case representation \eqref{eqn:rational1} is said to be {\sl canonical}.   
We remark that Gosper \cite{Gosper} considered expression 
\eqref{eqn:rational1} in a similar but somewhat different situation 
where $P(w)$ and $Q(w+j)$ were coprime for every nonnegative integer 
$j$ with $S(w)$ being a polynomial.      
\begin{lemma} \label{lem:rational1}
Any rational function $R(w) \in \mathbb{C}(w)$ admits a canonical 
representation \eqref{eqn:rational1}. 
If $R(w) \in \mathbb{R}(w)$ then $d$, $P(w)$, $Q(w)$ and $S(w)$ can be 
taken to be real.       
\end{lemma}
{\it Proof}. Start with the reduced expression $R(w) = d \cdot 
P_0(w)/Q_0(w)$, where $P_0(w)$ and $Q_0(w)$ are coprime monic 
polynomials. 
If they are strongly coprime then we are done with 
$P(w) = P_0(w)$, $Q(w) = Q_0(w)$ and $S(w) = 1$. 
Otherwise the argument proceeds as follows. 
Suppose that there is a representation \eqref{eqn:rational1} 
in which $P(w)$ and $Q(w)$ are coprime but {\sl not} strongly coprime. 
Then either there exists a positive integer $i$ such that $P(w+i)$ and 
$Q(w)$ have a common factor or there exists a positive integer $j$ such 
that $P(w)$ and $Q(w+j)$ have a common factor. 
In the former case there is a number $\alpha \in \mathbb{C}$ such that 
$(w+\alpha)|P(w)$ and $(w+\alpha+i)|Q(w)$.  
Put 
\begin{equation} \label{eqn:step1}
P_1(w) := \frac{P(w)}{w+\alpha}, \qquad 
Q_1(w) := \frac{Q(w)}{w+\alpha+i}, \qquad 
S_1(w) := \frac{S(w)}{(w+\alpha)_i}. 
\end{equation}
In the latter case there is a number $\beta \in \mathbb{C}$ such that 
$(w+\beta+j)|P(w)$ and $(w+\beta)|Q(w)$.  
Put 
\begin{equation} \label{eqn:step2} 
P_1(w) := \frac{P(w)}{w+\beta+j}, \qquad 
Q_1(w) := \frac{Q(w)}{w+\beta},  \qquad 
S_1(w) := S(w) \cdot (w+\beta)_j. 
\end{equation}
In either case it is easy to see that $P_1(w)$ and $Q_1(w)$ are 
coprime monic polynomials such that  
\[
 \frac{S_1(w+1)}{S_1(w)} \cdot d \cdot \frac{P_1(w)}{Q_1(w)} 
= \frac{S(w+1)}{S(w)} \cdot d \cdot \frac{P(w)}{Q(w)} = R(w),    
\]
with $\deg P_1(w) = \deg P(w)-1$ and $\deg Q_1(w) = \deg Q(w)-1$.  
If $P_1(w)$ and $Q_1(w)$ are strongly coprime then we are done. 
Otherwise, repeat the same procedure. 
This process must terminate in finite steps because the degrees of 
$P(w)$ and $Q(w)$ decrease by one in each step.  
\par
The proof of the second assertion requires a slight modification 
of the above argument.   
Suppose that $R(z)$ is real. 
Then $d$, $P_0(w)$ and $Q_0(w)$ above can be taken to be real. 
The induction procedure $(P,Q,S) \mapsto (P_1,Q_1,S_1)$ in 
formula \eqref{eqn:step1} resp. \eqref{eqn:step2} carries 
over if $\alpha$ resp. $\beta$ is real. 
But if $\alpha$ resp. $\beta$ is not real then 
formula \eqref{eqn:step1} resp. \eqref{eqn:step2} should be 
replaced by  
\begin{align*} 
P_1(w) &:= \textstyle \frac{P(w)}{(w+\alpha)(w+\bar{\alpha})}, &
Q_1(w) &:= \textstyle \frac{Q(w)}{(w+\alpha+i)(w+\bar{\alpha}+i)}, &  
S_1(w) &:= \textstyle \frac{S(w)}{(w+\alpha)_i \cdot (w+\bar{\alpha})_i}, \\
P_1(w) &:= \textstyle \frac{P(w)}{(w+\beta+j)(w+\bar{\beta}+j)}, &
Q_1(w) &:=  \textstyle \frac{Q(w)}{(w+\beta)(w+\bar{\beta})}, &
S_1(w) &:= S(w) \cdot (w+\beta)_j \cdot (w+\bar{\beta})_j,  
\end{align*}
respectively.  
This is well defined since if a real polynomial has a non-real 
root then its complex conjugate is also a root of the same 
polynomial.  
The modified procedure keeps realness, so that the real initial 
data $(d, P_0, Q_0)$ leads to a final real output 
$(S, d, P, Q)$. \hfill $\Box$ 
\begin{proposition} \label{prop:pole2} 
Let $f(w)$ be a non-elementary solution to Problem $\ref{prob:ocf}$ in 
region \eqref{eqn:h-strip}.  
If representation \eqref{eqn:canonical} is canonical, then we have 
$m=n \ge 1$, the number $r$ must be a positive integer with $r \ge m$ 
and there exist integers $s_1, \dots, s_m$ mutually distinct modulo $r$ 
such that \eqref{eqn:ui} holds. 
\end{proposition}
{\it Proof}. 
Let $W_{\mathrm{pole}}$ denote the set of all poles of $f(w)$, 
which is an infinite set since $f(w)$ is assumed to be non-elementary. 
In formula \eqref{eqn:f=g} the poles of $\varGamma(w+u_1) \cdots \varGamma(w+u_m)$ 
and those of $\varGamma(w+v_1) \cdots \varGamma(w+v_m)$ constitute two families of 
arithmetic progressions 
\begin{align} 
U_i &:= \{-u_i-k\}_{k = 0}^{\infty} \qquad \,\, (i = 1,\dots,m), \label{eqn:a-p1} \\ 
V_j &:= \{\, -v_j-k \, \}_{k = 0}^{\infty} \qquad (j = 1,\dots,m),  \label{eqn:a-p2} 
\end{align} 
respectively.    
Since representation \eqref{eqn:canonical} is canonical, $U_i$ and 
$V_j$ are disjoint for every $i, j = 1, \dots, m$, so that 
$W_{\mathrm{pole}}$ is commensurable to the union $\bigcup_{i=1}^m U_i$, 
where this union is disjoint because all poles of $f(w)$ are simple 
so that $u_i-u_j$ is not an integer for every $i \neq j$, that is,  
\begin{equation} \label{eqn:Wpole}
W_{\mathrm{pole}} \stackrel{\circ}{=} \coprod_{i=1}^m U_i. 
\end{equation}
Thus when expression \eqref{eqn:canonical} is canonical, $f(w)$ is 
non-elementary if and only if $m \ge 1$.  
\par
Take $i = 1$ and $k$ sufficiently large in the arithmetic 
progression \eqref{eqn:a-p1}. 
Equation \eqref{eqn:f=g} then shows that $w = -u_1-k$ and 
$w = -u_1-k-1$ are poles of $f(w) = g(w)$, so that they must lie 
in the arithmetic progression \eqref{eqn:a-p0}. 
Thus there exist nonnegative integers $j_1$ and $j_2$ with 
$j_1 < j_2$ such that $-u_1-k = -j_1/r$ and $-u_1-k-1 = -j_2/r$. 
Taking their difference gives $1 = (j_2-j_1)/r$, which shows 
that $r = j_2-j_1$ must be a positive integer.  
Similarly, for each $i = 1,\dots,m$ there exists an integer $k_i$ 
such that $w = -u_i-k_i$ is a pole of $g(w) = f(w)$ and so it 
must lie in the arithmetic progression \eqref{eqn:a-p0}, namely, 
it can be written $-u_i-k_i = -j_i/r$ for some integer $j_i$.  
If we put $s_i := j_i+r k_i$ then formula \eqref{eqn:ui} holds true. 
Note that $s_1, \dots,s_m$ are mutually distinct modulo $r$, 
because $U_i, \dots, U_m$ are mutually disjoint. 
\hfill $\Box$ \par\medskip
Note that putting Propositions \ref{prop:f=g} and \ref{prop:pole2} 
together yields Theorem \ref{thm:gpf}. 
\section{Asymptotics of the Residues} \label{sec:a-r}
Throughout this section let $\lambda = (p,q,r;a,b,x)$ and the 
associated $f(w) = f(w;\lambda)$ be a non-elementary solution to 
Problem \ref{prob:ocf} in region \eqref{eqn:h-strip} with formula 
\eqref{eqn:canonical} in a canonical form. 
The poles of $f(w)$ are commensurable to the disjoint 
union of $m$ arithmetic progressions $U_i$ $(i = 1, \dots, m)$ as 
in formula \eqref{eqn:Wpole}.   
In view of formulas \eqref{eqn:a-p0} and \eqref{eqn:ui} the general 
term of $U_i$ is expressed as $-u_i - k = -(r k + s_i)/r = w_j$ 
with $j = r k + s_i$. 
In this situation, if we put 
\[
\mathrm{Res}_k^{(i)} := \underset{\scriptstyle w = w_j}{\mathrm{Res}} \, f(w), \quad  
C_k^{(i)} := C_j, \quad F_k^{(i)} :=  {}_2F_1(a_j, \, b_j; \, j+2; \, x), 
\]
using the notation of Lemma \ref{lem:residue}, then 
formula \eqref{eqn:residue} reads  
\begin{equation} \label{eqn:residue2} 
\mathrm{Res}_k^{(i)} = C_k^{(i)} \cdot F_k^{(i)}.   
\end{equation} 
We study the asymptotic behavior of $\mathrm{Res}_k^{(i)}$ as $k \to \infty$ 
for a fixed $i = 1, \dots, m$.      
\begin{lemma} \label{lem:F_k^i} 
Let $B$ and $t_0$ be the same constants as in definition \eqref{eqn:AB} 
and put  
\begin{align}
\xi(t) &:= t^{-2 p/r-1} (1-t)^{2 p/r-3} (1-x t)^{2 q/r-1}, \label{eqn:xi} \\
\widetilde{A} &:= \dfrac{p^{a+2 p/r-1/2}}{(r-p)^{a+2 p/r-3/2} \cdot r^{1/2}} 
\cdot \dfrac{\xi(t_0)}{\eta(t_0) \sqrt{\phi''(t_0)}}, \label{eqn:tilde-A}  
\end{align}
where $\phi(t)$ and $\eta(t)$ are defined by formula \eqref{eqn:Phi-eta}. 
Then for each $i = 1, \dots, m$, we have  
\begin{equation} \label{eqn:F_k^i}
F_k^{(i)} \sim \{(1-x)^{(p+q-r)u_i-a-b} \widetilde{A} B^{u_i +2/r}\} 
\cdot \{(1-x)^{p+q-r} B \}^k \qquad \mbox{as} \quad k \to \infty.  
\end{equation} 
\end{lemma} 
{\it Proof}. Euler's transformation \eqref{eqn:transf1} and the 
definitions of $a_j$ and $b_j$ in Lemma \ref{lem:residue} yield 
\[ 
\begin{split}
F_k^{(i)} &= (1-x)^{j+2-a_j-b_j} \, 
{}_2F_1(j+2-a_j, \, j+2-b_j; \, j+2; \, x) \quad \mbox{with} \quad 
j = r k + s_i \\
&= (1-x)^{(p+q-r)(k+u_i)-a-b} \, 
{}_2F_1(p \, w_k^{(i)} + \widetilde{a}, \, 
q \, w_k^{(i)} + \widetilde{b}; \, r \, w_k^{(i)}; \, x),  
\end{split} 
\]
where $\widetilde{a} := 1-a-2p/r$, $\widetilde{b} := 1-b-2q/r$ and 
$w_k^{(i)} := k+u_i+2/r$. 
Asymptotic behavior of $F(p \, w_k^{(i)} + \widetilde{a}, 
\, q \, w_k^{(i)} + \widetilde{b}, \, r \, w_k^{(i)} ; \, x)$ can be extracted 
from that of $f(w) = F(pw+a, qw+b, rw; x)$ in formula \eqref{eqn:asympt2} 
by substitution $a \mapsto \widetilde{a}$, $b \mapsto \widetilde{b}$, 
$w \mapsto w_k^{(i)}$, where $p$, $q$, $r$, $x$ and so   
$\Phi(t)$, $\phi(t)$, $t_0$, $B$ in formula \eqref{eqn:AB} are 
left unchanged.  
This substitution replaces $\eta(t)$ with 
\[
\widetilde{\eta}(t) := t^{-a-2p/r}(1-t)^{a+2p/r-2}(1-xt)^{b+2q/r-1} 
= \dfrac{t^{-1-2p/r}(1-t)^{2p/r-3}(1-xt)^{2q/r-1}}{t^{a-1}(1-t)^{-a-1}(1-xt)^{-b}} 
= \dfrac{\xi(t)}{\eta(t)},   
\]
where $\xi(t)$ is defined by formula \eqref{eqn:xi}, which in turn induces 
the change of constant $A \mapsto \widetilde{A}$ in formula \eqref{eqn:tilde-A}.  
Lemma \ref{lem:asympt} then yields $F_k^{(i)} \sim (1-x)^{(p+q-r)(k+u_i)-a-b} 
\cdot \widetilde{A} \cdot B^{k+u_i+2/r}$ as $k \to \infty$. 
After a rearrangement, it just gives the desired formula \eqref{eqn:F_k^i}. 
\hfill $\Box$ \par\medskip 
We proceed to investigating $C_k^{(i)}$. 
Substituting $j = rk+s_i$ into formula \eqref{eqn:Cj} yields 
\begin{equation} \label{eqn:C_k^i}
C_k^{(i)} = \dfrac{(-1)^{r k+s_i}}{r} \cdot 
\dfrac{(-(pk+\alpha_i))_{rk+s_i+1} (-(qk+\beta_i))_{rk+s_i+1}}{(rk+s_i)! \, (rk+s_i+1)!} 
\, x^{rk+s_i+1},  
\end{equation} 
where $\alpha_i := p u_i-a$ and $\beta_i := q u_i - b$ as in definition 
\eqref{eqn:ai-bi}.  
Using this formula we study the asymptotic behavior of $C_k^{(i)}$ as 
$k \to \infty$ under the condition \eqref{eqn:h-strip} where 
$0 < p < r$ and $q < r$.    
\begin{lemma} \label{lem:factorial}  
For each $i = 1, \dots, m$, according to the value of $q$ we have      
\begin{subequations} \label{eqn:factorial1}
\begin{align}
C_k^{(i)} & \sim D_1^{(i)} \cdot E_1^k \cdot 
(-1)^{[p k+\alpha_i]+[q k+\beta_i]+r k+s_i} \cdot 
\sin\pi\{p k+\alpha_i\} \cdot \sin\pi\{q k + \beta_i\} & & (0 < q < r),  
\label{eqn:factorial1-1} \\ 
C_k^{(i)} & \sim D_2^{(i)} \cdot E_2^k \cdot k^{b-1/2} \cdot 
(-1)^{[p k+\alpha_i]+r k+s_i+1} \cdot \sin\pi\{p k+\alpha_i\} & & (q = 0 < r),  
\label{eqn:factorial1-2} \\
C_k^{(i)} &\sim D_3^{(i)} \cdot E_3^k \cdot 
(-1)^{[p k+\alpha_i]+r k+s_i+1} \cdot \sin\pi\{p k+\alpha_i\} & & (q < 0 < r), 
\label{eqn:factorial1-3} 
\end{align}
\end{subequations} 
as $k \to \infty$, where $D_{\nu}^{(i)}$ and $E_{\nu}$, $\nu = 1,2,3$, are 
constants defined by  
\begin{subequations} \label{eqn:DE123}
\begin{align}
D_1^{(i)} &:= \textstyle
\frac{2 \, p^{\alpha_i+1/2} (r-p)^{s_i-\alpha_i+1/2} q^{\beta_i+1/2} 
(r-q)^{s_i-\beta_i+1/2}}{\pi \, r^{2 s_i+3}} \, x^{s_i+1}, & 
E_1 &:= \textstyle 
\frac{p^p(r-p)^{r-p} q^q (r-q)^{r-q}}{r^{2 r}} \, x^r, 
\label{eqn:DE1} \\ 
D_2^{(i)} &:= \textstyle 
\sqrt{\frac{2}{\pi}} \cdot 
\frac{p^{\alpha_i+1/2} (r-p)^{s_i-\alpha_i+1/2}}{\varGamma(b) \, r^{s_i-b+5/2}} 
\, x^{s_i+1}, & 
E_2 &:= \textstyle  
\frac{p^p(r-p)^{r-p}}{r^r} \, x^r, 
\label{eqn:DE2} \\
D_3^{(i)} &:= \textstyle 
\frac{p^{\alpha_i+1/2} (r-p)^{s_i-\alpha_i+1/2} |q|^{\beta_i+1/2} 
(r-q)^{s_i-\beta_i+1/2}}{\pi \, r^{2 s_i+3}} \, x^{s_i+1}, & 
E_3 &:= \textstyle
\frac{p^p(r-p)^{r-p}|q|^q(r-q)^{r-q}}{r^{2 r}} \, x^r. 
\label{eqn:DE3}
\end{align}
\end{subequations}
\end{lemma}
{\it Proof}. 
It follows from $0 < p < r$ that $[pk+\alpha_i]+1$ and 
$rk+s_i-[pk+\alpha_i]$ are positive integers for every sufficiently 
large integer $k$. 
Since $pk+\alpha_i = [pk+\alpha_i] + \{pk+\alpha_i\}$, we have  
\[ 
\begin{split} 
(-(p k + \alpha_i))_{r k+s_i+1} &= \overbrace{(-(p k+\alpha_i))(1-(p k+\alpha_i)) 
\cdots (-\{ p k+\alpha_i \})}^{[p k+\alpha_i]+1} \\
& \quad \times \overbrace{(1-\{p k+\alpha_i\})(2-\{p k+\alpha_i\}) \cdots 
(r k+s_i - (p k+\alpha_i))}^{r k+s_i-[p k+\alpha_i]} \\
&= (-1)^{[p k+\alpha_i]+1} (\{p k+\alpha_i\})_{[p k+\alpha_i]+1} \, 
(1-\{p k+\alpha_i\})_{r k+s_i-[p k+\alpha_i]} \\
&= (-1)^{[p k+\alpha_i]+1} \frac{\varGamma(p k+\alpha_i+1)}{\varGamma(\{p k+\alpha_i\})} 
\cdot \frac{\varGamma((r-p)k+s_i-\alpha_i+1)}{\varGamma(1-\{p k+\alpha_i\})} \\
&= (-1)^{[p k+\alpha_i]+1} \frac{\sin\pi\{p k+\alpha_i\}}{\pi} \,  
\varGamma(p k+\alpha_i+1) \, \varGamma((r-p)k+s_i-\alpha_i+1),  
\end{split} 
\]
by the recursion and reflection formulas for the gamma function. 
By Stirling's formula \eqref{eqn:stirling}, 
\[
\frac{\varGamma(pk+\alpha_i+1) \varGamma((r-p)k+s_i-\alpha_i+1)}{\varGamma(r k+s_i+3/2)} 
\sim \sqrt{2\pi} \, \frac{p^{\alpha_i+1/2}(r-p)^{s_i-\alpha_i+1/2}}{r^{s_i+1}} 
\left(\frac{p^p(r-p)^{r-p}}{r^r} \right)^k 
\]
as $k \to \infty$. 
Using this asymptotic formula in the above equation we have   
\begin{equation} \label{eqn:factorial-p}
\frac{(-(p k + \alpha_i))_{r k+s_i+1}}{\varGamma(r k+s_i+3/2)} \sim \textstyle 
\sqrt{\frac{2}{\pi}} \cdot 
\frac{p^{\alpha_i+1/2}(r-p)^{s_i-\alpha_i+1/2}}{r^{s_i+1}} 
\left(\frac{p^p(r-p)^{r-p}}{r^r} \right)^k \!  
(-1)^{[p k+\alpha_i]+1} \sin\pi\{p k +\alpha_i\}   
\end{equation}
as $k \to \infty$. 
Exactly in the same manner, if $0 < q < r$ then we have as $k \to \infty$, 
\begin{equation} \label{eqn:factorial-q1}
\frac{(-(qk + \beta_i))_{rk+s_i+1}}{\varGamma(rk+s_i+3/2)} \sim \textstyle 
\sqrt{\frac{2}{\pi}} \cdot 
\frac{q^{\beta_i+1/2}(r-q)^{s_i-\beta_i+1/2}}{r^{s_i+1}} 
\left(\frac{q^q(r-q)^{r-q}}{r^r} \right)^k \!
(-1)^{[qk+\beta_i]+1} \sin\pi\{q k +\beta_i\}.    
\end{equation}
\par
Next we consider the case $q \le 0 < r$. 
For every sufficiently large integer $k$,  
\[
\frac{(-(qk+\beta_i))_{rk+s_i+1}}{\varGamma(rk+s_i+3/2)} = 
\begin{cases} 
\dfrac{\varGamma(b+r k+s_i+1)}{\varGamma(b) \varGamma(rk+s_i+3/2)} \quad & (q = 0 < r), \\[3mm]
\dfrac{\varGamma((r-q)k+s_i-\beta_i+1)}{\varGamma(|q|k-\beta_i)\varGamma(rk+s_i+3/2)}  
\quad & (q < 0 < r).
\end{cases} 
\]
Applying Stirling's formula \eqref{eqn:stirling} to the right-hand side 
above we have as $k \to \infty$,    
\begin{equation} \label{eqn:factorial-q2}
\frac{(-(qk + \beta_i))_{rk+s_i+1}}{\varGamma(rk+s_i+3/2)} \sim \textstyle
\begin{cases} 
(rk)^{b-1/2}/\varGamma(b) & (q = 0 < r), \\[1mm]  
\frac{1}{\sqrt{2\pi}} \cdot 
\frac{|q|^{\beta_i+1/2}(r-q)^{s_i-\beta_i+1/2}}{r^{s_i+1}} 
\left(\frac{|q|^q (r-q)^{r-q}}{r^r} \right)^k & (q < 0 < r),       
\end{cases} 
\end{equation}
\par
Notice that $(rk+s_i)! \, (rk+s_i+1)! \sim \varGamma(rk+s_i+3/2)^2$ as 
$k \to \infty$ by Stirling's formula \eqref{eqn:stirling}. 
Thus substituting formulas \eqref{eqn:factorial-p} and 
\eqref{eqn:factorial-q1} into \eqref{eqn:C_k^i} yields 
formula \eqref{eqn:factorial1-1}. 
Similarly substituting formulas \eqref{eqn:factorial-p} and 
\eqref{eqn:factorial-q2} into \eqref{eqn:C_k^i} yields 
formulas \eqref{eqn:factorial1-2} and \eqref{eqn:factorial1-3}.  
\hfill $\Box$ 
\begin{remark} \label{rem:factorial} 
When $q = 0$, the number $b$ must not be a non-positive integer, since 
elementary solutions \eqref{eqn:elementary1} are excluded from our 
consideration. 
Thus the constant $D_2^{(i)}$ in \eqref{eqn:DE2} is nonzero.  
\end{remark}
\begin{proposition} \label{prop:asympt-Res1} 
For each $i = 1, \dots, m$, according to the value of $q$ we have    
\begin{subequations} \label{eqn:factorial2}
\begin{align}
\mathrm{Res}_k^{(i)} &\sim D_4^{(i)} \cdot E_4^k \cdot 
(-1)^{[p k+\alpha_i]+[q k+\beta_i]+r k+s_i}  
\sin\pi\{p k+\alpha_i\} \cdot \sin\pi\{q k + \beta_i\} && (0 < q < r), 
\label{eqn:factorial2-1} \\ 
\mathrm{Res}_k^{(i)} &\sim D_5^{(i)} \cdot E_5^k \cdot k^{b-1/2} \cdot 
(-1)^{[p k+\alpha_i]+r k+s_i+1} \sin\pi\{pk+\alpha_i\} && (q = 0 < r), 
\label{eqn:factorial2-2} \\  
\mathrm{Res}_k^{(i)} &\sim D_6^{(i)} \cdot E_6^k \cdot 
(-1)^{[p k+\alpha_i]+r k+s_i+1} \sin\pi\{p k+\alpha_i\} && (q < 0 < r), 
\label{eqn:factorial2-3} 
\end{align}  
\end{subequations}
as $k \to \infty$, where $D_{\nu}^{(i)}$ and $E_{\nu}$, $\nu = 4, 5, 6$, are 
constants defined by  
\begin{equation} \label{eqn:DE456}
D_{\nu}^{(i)} := (1-x)^{(p+q-r)u_i-a-b} \widetilde{A} B^{u_i+2/r} D_{\nu-3}^{(i)}, 
\quad  
E_{\nu} := (1-x)^{p+q-r} B \cdot E_{\nu-3} \quad (\nu = 4,5,6).   
\end{equation}
\end{proposition}
{\it Proof}. 
This proposition is proved by putting Lemmas \ref{lem:F_k^i} and 
\ref{lem:factorial} together. \hfill $\Box$ 
\begin{lemma} \label{lem:asympt-Res2}
In the circumstances of Proposition $\ref{prop:f=g}$, 
we have for each $i = 1,\dots,m$,  
\begin{equation} \label{eqn:asympt-Res2}
\mathrm{Res}_k^{(i)} \sim K^{(i)} \cdot B^{-k} \quad \mbox{as} \quad  
k \to \infty, \qquad 
K^{(i)} := \frac{(-1)^{s_0}A}{\pi B^{u_i}} \cdot 
\frac{\prod_{j=1}^m \sin\pi(v_j-u_i)}{\prod_{j=1}^{* \, m} \sin\pi(u_j-u_i)},    
\end{equation}
where $\prod_{j=1}^{* \, m}$ denotes the product taken over all 
$j = 1, \dots, m$ but $j = i$. 
\end{lemma}
{\it Proof}. 
Applying the reflection formula for the gamma function to  
the identity \eqref{eqn:f=g} yields  
\[
f(w) = S(w) \cdot B^w \cdot 
\frac{\prod_{j=1}^m \varGamma(1-v_j-w)}{\prod_{j=1}^m \varGamma(1-u_j-w)} 
\cdot \frac{\prod_{j=1}^m \sin\pi(w+v_j)}{\prod_{j=1}^m \sin\pi(w+u_j)},  
\]
where $d = B$ is used in view of Proposition \ref{prop:f=g}. 
Taking its residue at $w = -k-u_i$ gives 
\[
\begin{split}
\mathrm{Res}_k^{(i)} &= S(-(k+u_i)) \cdot B^{-k-u_i} 
\frac{\prod_{j=1}^m \varGamma(k+u_i-v_j+1)}{\prod_{j=1}^m \varGamma(k+u_i-u_j+1)} \\
& \quad \times  
\frac{\prod_{j=1}^m \sin\pi(v_j-u_i-k)}{\prod_{j=1}^{* \, m} \sin\pi(u_j-u_i-k)} 
\lim_{w \to -k-u_i} \frac{w+k+u_i}{\sin\pi(w+u_i)}  \\
&= \frac{1}{\pi} S(-(k+u_i)) B^{-k-u_i} 
\frac{\prod_{j=1}^m \varGamma(k+u_i-v_j+1)}{\prod_{j=1}^m \varGamma(k+u_i-u_j+1)} 
\cdot \frac{\prod_{j=1}^m \sin\pi(v_j-u_i)}{\prod_{j=1}^{* \, m} \sin\pi(u_j-u_i)}.  
\end{split} 
\]
By Proposition \ref{prop:f=g}, $S(-(k+u_i)) \sim S_0 (-k)^{s_0} = A 
(-1)^{s_0} k^{s_0}$ and by Stirling's formula \eqref{eqn:stirling},  
\[
\frac{\prod_{j=1}^m \varGamma(k+u_i-v_j+1)}{\prod_{j=1}^m \varGamma(k+u_i-u_j+1)} 
\sim k^{u-v} = k^{-s_0} \qquad \mbox{as} \quad k \to \infty. 
\]
Substituting these asymptotic formulas into the above equation 
we get formula \eqref{eqn:asympt-Res2}. \hfill $\Box$ \par\medskip
It follows from definition \eqref{eqn:asympt-Res2} that $K^{(i)} \neq 0$ for 
$i = 1, \dots, m$, because $A$ and $B$ are positive by formula \eqref{eqn:AB} 
and $v_j -u_i \not\in \mathbb{Z}$ for $i,j = 1, \dots, m$, by the assumption 
that \eqref{eqn:canonical} be canonical. 
\section{Applying Kronecker's Theorem} \label{sec:kronecker}
Let $\lambda = (p,q,r;a,b;x)$ be the same as in \S\ref{sec:a-r}. 
In this section we shall describe how the asymptotic results in Proposition 
\ref{prop:asympt-Res1} and Lemma \ref{lem:asympt-Res2} are combined 
with Kronecker's theorem on Diophantine approximations 
to obtain an arithmetic result on $p$ and $q$, where we have already known   
that $r$ must be a positive integer by Proposition $\ref{prop:pole2}$.   
Since asymptotic representations \eqref{eqn:factorial2} and 
\eqref{eqn:asympt-Res2} must be equivalent, taking the ratio of them gives           
\begin{subequations} \label{eqn:factorial3} 
\begin{align}
D_7^{(i)} \cdot E_7^k \cdot(-1)^{[p k+\alpha_i]+[q k+\beta_i]+r k+s_i} \cdot 
\sin\pi\{p k+\alpha_i\} \cdot \sin\pi\{q k + \beta_i\} \to 1 \quad & 
(0 < q < r), \label{eqn:factorial3-1} \\     
D_8^{(i)} \cdot E_8^k \cdot k^{b-1/2} \cdot 
(-1)^{[p k+\alpha_i]+r k+s_i+1} \cdot \sin\pi\{p k+\alpha_i\} \to 1 \quad & 
(q = 0 < r),  \label{eqn:factorial3-2} \\    
D_9^{(i)} \cdot E_9^k \cdot 
(-1)^{[p k+\alpha_i]+r k+s_i+1} \cdot \sin\pi\{p k+\alpha_i\} \to 1 \quad & 
(q < 0 < r),  \label{eqn:factorial3-3}    
\end{align}
\end{subequations}
as $k \to \infty$, where $D_{\nu}^{(i)}$ and $E_{\nu}$, $\nu = 7, 8, 9$, 
$i = 1, \dots, m$, are constants defined by  
\begin{equation} \label{eqn:DE789}
D_{\nu}^{(i)} := D_{\nu-3}^{(i)}/K^{(i)}, \qquad E_{\nu} := B \cdot E_{\nu-3}. 
\qquad (\nu = 7,8,9). 
\end{equation}
\par 
Taking the absolute values of formulas \eqref{eqn:factorial3} yields       
\begin{subequations} \label{eqn:factorial4} 
\begin{align}
|D_7^{(i)}| \cdot E_7^k \cdot 
\sin\pi\{p k+\alpha_i\} \cdot \sin\pi\{q k + \beta_i\} \to 1 \qquad & 
(0 < q < r), \label{eqn:factorial4-1} \\     
|D_8^{(i)}| \cdot E_8^k \cdot k^{b-1/2} \cdot 
\sin\pi\{p k+\alpha_i\} \to 1 \qquad & (q = 0 < r),  \label{eqn:factorial4-2} \\    
|D_9^{(i)}| \cdot E_9^k \cdot \sin\pi\{p k+\alpha_i\} \to 1 \qquad & 
(q < 0 < r),  \label{eqn:factorial4-3}    
\end{align}
\end{subequations}
as $k \to \infty$. 
We study these formulas using Kronecker's approximation 
theorem \cite{Kronecker}.  
\begin{proposition} \label{prop:kronecker1} 
For each $i = 1, \dots, m$, formulas \eqref{eqn:factorial4-1}, 
\eqref{eqn:factorial4-2} and \eqref{eqn:factorial4-3} lead to 
\begin{subequations} \label{eqn:kronecker1}
\begin{alignat}{4}
p, q &\in \mathbb{Q}, \quad & 
 E_7 &= 1, \quad    &   
|D_7^{(i)}| \cdot \sin\pi\{p k+\alpha_i\} \cdot \sin\pi\{q k + \beta_i\} &= 1 
\qquad & (0 < q < r), \label{eqn:kronecker1-1} \\     
p    &\in \mathbb{Q}, \quad & 
 E_8 &= 1, \quad    & 
   b  = 1/2, \qquad  
|D_8^{(i)}| \cdot \sin\pi\{p k+\alpha_i\} &= 1 
\qquad & (q = 0 < r),  \label{eqn:kronecker1-2} \\    
p    &\in \mathbb{Q}, \quad & 
 E_9 &= 1, \quad    &  
|D_9^{(i)}| \cdot \sin\pi\{p k+\alpha_i\} &= 1 
\qquad & (q < 0 < r),  \label{eqn:kronecker1-3}       
\end{alignat}
\end{subequations}
for every integer $k \in \mathbb{Z}$, respectively.   
\end{proposition} 
{\it Proof}. 
We shall prove the implication \eqref{eqn:factorial4-1} $\Rightarrow$ 
\eqref{eqn:kronecker1-1} by using two-dimensional as well as 
one-dimensional versions of Kronecker's theorem.    
Proofs of the remaining implications \eqref{eqn:factorial4-2} $\Rightarrow$ 
\eqref{eqn:kronecker1-2} and \eqref{eqn:factorial4-3} $\Rightarrow$ 
\eqref{eqn:kronecker1-3} are left to the reader, because they use 
only  one-dimensional version and so less intricate.  
In what follows we fix an index $i = 1, \dots, m$. 
\par
First we claim that $p$, $q$ and $1$ must be linearly dependent over $\mathbb{Q}$. 
Suppose the contrary that they are linearly independent over $\mathbb{Q}$. 
By two-dimensional version of Kronecker's theorem, the sequence  
$(\{p k + \alpha_i\}, \, \{q k + \beta_i\})$ are dense in the square 
$[0, \, 1) \times [0, \, 1)$ as $k \to \infty$. 
In particular there exists a subsequence of the $k$'s along which 
$\{p k + \alpha_i\} \to 1/2$ and $\{q k + \beta_i\} \to 1/2$ so that 
$\sin\pi\{pk+\alpha_i\} \cdot \sin\pi\{qk + \beta_i\} \to 1$ as 
$k \to \infty$. 
Formula \eqref{eqn:factorial4-1} then says that along this 
subsequence $|D_7^{(i)}| \cdot E_7^k \to 1$ as $k \to \infty$, which forces 
$E_7 = 1$ and $|D_7^{(i)}| = 1$. 
But there exists another subsequence along which $\{p k + \alpha_i\} \to 0$ 
and $\{q k + \beta_i\} \to 0$ so that $\sin\pi\{p k+\alpha_i\} \cdot 
\sin\pi\{qk + \beta_i\} \to 0$ as $k \to \infty$. 
Formula \eqref{eqn:factorial4-1} now yields an absurd conclusion 
$0 \sim 1$ as $k \to \infty$ along the latter subsequence. 
Thus the claim is proved. 
\par
Next we shall show that $p$ and $q$ are rational (the proof will be 
completed at the end of next paragraph).  
Suppose the contrary that either $p$ or $q$ is irrational, where   
we may assume without loss of generality that $p$ is irrational. 
Since $p$, $q$ and $1$ are linearly dependent over $\mathbb{Q}$, there exist 
$\lambda, \, \mu \in \mathbb{Q}$ such that $q = \lambda p + \mu$.  
Let $\nu$ be the denominator of the reduced representation 
of $\mu$; by convention let $\nu = 1$ when $\mu = 0$. 
If we put $x_k = p\nu k + \alpha_i$ and $y_k = q \nu k + \beta_i$ for 
$k \in \mathbb{Z}$, then formula \eqref{eqn:factorial4-1} with $k$ replaced 
by $\nu k$ reads  
\begin{equation} \label{eqn:factorial5}
|D_7^{(i)}| \cdot E_7^{\nu k} \cdot 
\sin\pi\{x_k\} \cdot \sin\pi\{y_k \} \to 1 \qquad \mbox{as} \quad 
k \to \infty.  
\end{equation}
Observe that $y_k = \lambda x_k + \gamma_i + \mu \nu k$ with 
$\gamma_i := \beta_i-\lambda \alpha_i$ and $\mu \nu k \in \mathbb{Z}$, 
so that $\{y_k\} = \{\lambda x_k + \gamma_i\}$. 
Since $p\nu$ is irrational in $x_k = p\nu k + \alpha_i$, 
the limit set of the sequence $\{x_k\}$ as $k \to \infty$ is the 
whole unit interval $[0, \, 1]$ by one-dimensional version of 
Kronecker's theorem.   
With this in mind we describe the limit set of the sequence 
$(\{x_k\}, \, \{y_k\})$ as $k \to \infty$.   
If it is thought of as a sequence on the torus $\mathbb{T}^2 = 
\mathbb{R}^2/\mathbb{Z}^2$, its limit set is the torus line 
coming down from a line $y = \lambda x + \gamma_i$ in the universal 
covering $\mathbb{R}^2_{(x,y)}$.   
If $\lambda = \lambda_2/\lambda_1$ is the reduced representation 
of $\lambda$ (by convention put $\lambda_1 = 1$ and $\lambda_2 = 0$ 
when $\lambda = 0$), then the limit set is a 
$(\lambda_1, \lambda_2)$-torus knot as in Figure \ref{fig:torus}. 
Viewed on the square $[0, \, 1) \times [0, \, 1)$, it is a finite 
union of parallel line segments as in Figure \ref{fig:sin-sin-line}; 
when $\lambda = 0$, it is a single line segment parallel to the $x$-axis. 
\par
Consider the function $\varphi(x,y) = \sin \pi x \cdot \sin \pi y$ 
defined for $(x,y) \in [0, \, 1) \times [0, \, 1)$.    
For each $0 < h < 1$ the $h$-level set $\varphi(x,y) = h$ is a simple closed 
curve whose interior is a convex bounded domain; the $0$-level set is 
the union  of two lines $x = 0$ and $y = 0$; while the $1$-level set is 
a single point $(1/2, 1/2)$ (see Figure \ref{fig:sin-sin-line}). 
Thus it is clear that if $\lambda$ is nonzero then any single level 
set of the function $\varphi(x,y)$ cannot contain the limit set of the 
sequence $(\{x_k\}, \, \{y_k\})$.  
This is also the case when $\lambda = 0$, since the limit set 
is away from the $x$-axis. 
Indeed, if $\lambda = 0$ then $\{y_k\} = {\gamma_i}$ for every $k \in \mathbb{Z}$ 
so that the limit set is the line $y =\{\gamma_i\}$, where $\{\gamma_i\}$ 
must be nonzero for otherwise formula \eqref{eqn:factorial5} would imply 
an absurd conclusion $0 \sim 1$ as $k \to \infty$. 
Therefore there exist two subsequences of the $k$'s and two numbers 
$\sigma_1 \neq \sigma_2$ such that $\sin \pi\{x_k\} \cdot \sin \pi \{y_k\} 
\to \sigma_1$ along the first subsequence while $\sin \pi\{x_k\} \cdot 
\sin \pi \{y_k\} \to \sigma_2$ along the second one, where one 
may assume $\sigma_1 \neq 0$. 
Formula \eqref{eqn:factorial5} then implies $|D_7^{(i)}| \cdot E_7^{\nu k} 
\to 1/\sigma_1$ along the first subsequence, which forces 
$E_7^{\nu} = 1$ and $|D_7^{(i)}| = 1/\sigma_1$.     
But by taking the limit along the second subsequence, this formula   
yields $\sigma_2/\sigma_1 = 1$, which contradicts 
$\sigma_1 \neq \sigma_2$, showing that $p$ and $q$ must be rational. 
\par
Now that $p$ and $q$ are rational, there is a positive integer $\nu$ 
such that $\nu p$ and $\nu q$ are integers. 
If we put $c_k := |D_7^{(i)}| \cdot \sin \pi\{p k + \alpha_i\} 
\cdot \sin \pi\{q k + \beta_i\}$ for $k \in \mathbb{Z}$, then $\{c_k\}$ is a 
periodic sequence of period $\nu$, while formula \eqref{eqn:factorial4-1} 
reads $c_k \cdot E_7^k \to 1$ as $k \to \infty$. 
If $k$ is replaced by $\nu k$, one has $c_0 \cdot E_7^{\nu k} = 
c_{\nu k} \cdot E_7^{\nu k} \to 1$ as $k \to \infty$, which forces 
$E_7^{\nu} = 1$ and so $E_7 = 1$ since $E_7 > 0$. 
Now that $E_7 = 1$, formula \eqref{eqn:factorial4-1} 
becomes $c_k \to 1$ as $k \to \infty$. 
But this occurs with the periodic sequence $\{c_k\}$ only when 
$c_k = 1$ for every $k \in \mathbb{Z}$. 
Thus all the assertions in \eqref{eqn:kronecker1-1} are established.   
\hfill $\Box$   
\begin{proposition} \label{prop:kronecker2} 
For each $i = 1, \dots, m$, formulas \eqref{eqn:factorial3} and 
\eqref{eqn:kronecker1} imply that 
\begin{enumerate} 
\item if $0 < q < r$ then the parity of $[p k+\alpha_i]+[q k+\beta_i]+rk$ is  
independent of $k \in \mathbb{Z}$,  
\item if $q \le 0 < r$ then the parity of $[p k+\alpha_i]+rk$ is  
independent of $k \in \mathbb{Z}$,         
\end{enumerate} 
where the parity of an integer refers to whether it is odd or even. 
\end{proposition}
{\it Proof}. We prove assertion (1).  
It follows from formulas \eqref{eqn:factorial3-1} and \eqref{eqn:kronecker1-1}     
that $D_7^{(i)} \in \mathbb{R}^{\times}$ and $\mathrm{sgn} \, D_7^{(i)} \cdot 
(-1)^{[p k+\alpha_i]+[q k+\beta_i]+r k+s_i} \to 1$ as $k \to \infty$, where 
$\mathrm{sgn} \, D_7^{(i)} = \pm1$ is the signature of $D_7^{(i)}$. 
Thus the integer sequence $[p k+\alpha_i]+[q k+\beta_i]+r k$ has a stable parity 
as $k \to \infty$. 
Since $p$ and $q$ are rational by \eqref{eqn:kronecker1-1}, the sequence 
is periodic so that it must have a constant parity.  
In a similar manner assertion (2) follows from formulas \eqref{eqn:factorial3-2} 
and \eqref{eqn:kronecker1-2} when $q = 0 < r$ and from formulas 
\eqref{eqn:factorial3-3} and \eqref{eqn:kronecker1-3} when 
$q < 0 < r$, respectively.    
\hfill $\Box$
\begin{lemma} \label{lem:B} 
Everywhere in region \eqref{eqn:h-strip} the dilation constant is given by 
\begin{equation} \label{eqn:B}
d = B = \dfrac{r^r}{\sqrt{p^p \cdot |q|^q \cdot (r-p)^{r-p} \cdot 
(r-q)^{r-q} \cdot x^r \cdot (1-x)^{p+q-r}}},  
\end{equation}
where $|q|^q = 1$ when $q = 0$; this convention is reasonable since 
$|q|^q \to 1$ as $q \to 0$.   
\end{lemma}
{\it Proof}. 
We have $d = B$ by Proposition \ref{prop:f=g} and $B > 0$ by 
formula \eqref{eqn:AB}. 
If $0 < q < r$ then we use formulas \eqref{eqn:kronecker1-1}, \eqref{eqn:DE789}, 
\eqref{eqn:DE456} and \eqref{eqn:DE1} in this order to obtain
\[
1 = E_7 = B \cdot E_4 = B^2 \cdot (1-x)^{p+q-r} E_1 
= B^2 \cdot r^{-2r} p^p q^q (r-p)^{r-p} (r-q)^{r-q} x^x (1-x)^{p+q-r},   
\]
which gives formula \eqref{eqn:B} since $B > 0$. 
In a similar manner we use  \eqref{eqn:kronecker1-2}, \eqref{eqn:DE789}, 
\eqref{eqn:DE456}, \eqref{eqn:DE2} when $q = 0 < r$; and 
\eqref{eqn:kronecker1-3}, \eqref{eqn:DE789}, 
\eqref{eqn:DE456}, \eqref{eqn:DE3}  when $q < 0 < r$. 
In either case formula \eqref{eqn:B} follows. 
\hfill $\Box$
\section{Level Curves of the Sine-Sine} \label{sec:level}
Let $p$, $q \in \mathbb{Q}^{\times}$ and $\alpha$, $\beta \in \mathbb{R}$. 
Proposition \ref{prop:kronecker1} leads us to consider when 
$\sin \pi \{p j + \alpha \} \cdot \sin \pi\{q j + \beta \}$ 
or $\sin \pi\{p j + \alpha \}$ is a constant sequence, that is, 
independent of $j \in \mathbb{Z}$. 
(Hereafter index is denoted by $j$ instead of $k$.)   
We begin with the latter case which is more tractable than the former.     
\begin{lemma} \label{lem:sin}
Let $p \in \mathbb{Q}^{\times}$ and $\alpha \in \mathbb{R}$. 
The sequence $\sin \pi\{p j+ \alpha \}$ is independent of $j \in \mathbb{Z}$ 
if and only if either $(1)$ $p \in \mathbb{Z}$ or $(2)$ $p \in 1/2 + \mathbb{Z}$ 
and $\alpha \equiv \pm 1/4 \mod 1$.  
\end{lemma}
\begin{figure}[t]
\begin{center}
\includegraphics[width=60mm,clip]{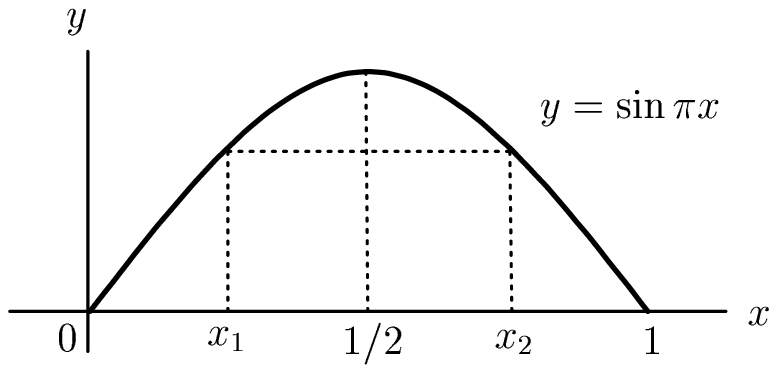}
\end{center}
\caption{Graph of the function $y = \sin \pi x$ $(0 \le x \le 1)$.}
\label{fig:sin}
\end{figure}
{\it Proof}. If $p \in \mathbb{Z}$ then $\{pj+ \alpha \} = \{\alpha \}$ 
and so $\sin \pi\{pj+\alpha \} = \sin \pi\{\alpha \}$ for every 
$j \in \mathbb{Z}$, being independent of $j$.   
This is case (1) of the lemma.  
Suppose that $p \in \mathbb{Q} \setminus \mathbb{Z}$ and let $p = l/k$ be its 
reduced representation with $k \ge 2$ and $(k, \, l) = 1$. 
Note that $\{p i + \alpha \} = \{p j + \alpha \}$ if and only if 
$i \equiv j \mod k$. 
Thus the sequence $\{pj+ \alpha \}$ takes exactly $k$ distinct values 
$\{pj + \alpha \}$ $(j = 0,\dots,k-1)$ as $j$ varies in $\mathbb{Z}$. 
In order for the sequence to be constant one must have 
\begin{equation}  \label{eqn:sin1} 
\sin \pi\{\alpha \} = \sin \pi\{p j + \alpha \} \qquad (j = 1, \dots, k-1). 
\end{equation}
Observe that the function $y = \sin \pi x$ $(0 \le x \le 1)$ 
is symmetric around $x = 1/2$ and strictly increasing 
in $0 \le x \le 1/2$ (see Figure \ref{fig:sin}). 
Thus $\sin \pi x$ cannot take a common value at distinct 
three points in $0 \le x \le 1$, so that one must have $k = 2$ in 
equation \eqref{eqn:sin1}.   
Since $(2, \, l) = (k, \, l) = 1$ the integer $l$ must be odd 
and hence $p = l/2 \in 1/2 + \mathbb{Z}$.  
Condition \eqref{eqn:sin1} now becomes $\sin\pi\{\alpha \} = 
\sin\pi\{p+ \alpha \} = \sin \pi \{\alpha + 1/2\}$, which is equivalent to 
$\{\alpha \}+\{\alpha + 1/2\} = 1$, because $\sin \pi x_1 = \sin \pi x_2$ for  
$0 \le x_1 < x_2 \le 1$ precisely when $x_1 + x_2 =1$ 
(see Figure \ref{fig:sin}). 
If $\alpha \equiv \alpha' \mod 1$ with $0 \le \alpha' < 1/2$ then 
$\{\alpha \} = \alpha'$ and $\{\alpha + 1/2\} = \alpha' + 1/2$ and so 
$\alpha' + (\alpha' + 1/2) = 1$ yields $\alpha' = 1/4$, namely, 
$\alpha \equiv 1/4 \mod 1$. 
If $\alpha \equiv \alpha' \mod 1$ with $1/2 \le \alpha' < 1$ then 
$\{\alpha \} = \alpha'$ and $\{\alpha + 1/2\} = \alpha' - 1/2$ and so 
$\alpha' + (\alpha' - 1/2) = 1$ yields $\alpha' = 3/4$, 
namely, $\alpha \equiv -1/4 \mod 1$.  \hfill $\Box$
\begin{proposition} \label{prop:sin2}  
Let $p$, $q \in \mathbb{Q}^{\times}$ and $\alpha$, $\beta \in \mathbb{R}$. 
The sequence $\sigma_j := \sin \pi\{p j+ \alpha \} \cdot \sin \pi\{q j + \beta \}$ 
is independent of $j \in \mathbb{Z}$ if and only if one of the following 
conditions is satisfied:   
\begin{enumerate}
\item either $(p, \alpha) \in \mathbb{Z}^2$ or $(q, \beta) \in \mathbb{Z}^2$,  
\item $p$, $q \in \mathbb{Z}$ and $\alpha$, $\beta \in \mathbb{R} \setminus \mathbb{Z}$, 
\item $p \in \mathbb{Z}$, $q \in 1/2 + \mathbb{Z}$, 
$\alpha \in \mathbb{R} \setminus \mathbb{Z}$ and 
$\beta \equiv \pm 1/4 \mod 1$, 
\item $p \in 1/2 + \mathbb{Z}$, $q \in \mathbb{Z}$, $\alpha \equiv \pm 1/4 \mod 1$ 
and $\beta \in \mathbb{R} \setminus \mathbb{Z}$, 
\item $p$, $q \in 1/2 + \mathbb{Z}$ and either $\alpha + \beta \in 1/2 + \mathbb{Z}$ or 
$\alpha - \beta \in 1/2 + \mathbb{Z}$,  
\item $p \in \varepsilon_p/3 + \mathbb{Z}$, $q \in \varepsilon_q/3 + \mathbb{Z}$, 
$\alpha \in \varepsilon_p(\delta/3 + \varepsilon \kappa) + \mathbb{Z}$, 
$\beta \in \varepsilon_q(\delta/3 - \varepsilon \kappa) + \mathbb{Z}$, 
where $\varepsilon_p, \, \varepsilon_q, \, \varepsilon = \pm 1$, 
$\delta = 0, \, \pm1$ and $\kappa$ is an irrational number defined by 
\begin{equation} \label{eqn:kappa}
\kappa = \frac{1}{\pi}\arctan \sqrt{3/5} = 0.209785\cdots, 
\end{equation}
\item $p \in \varepsilon_p/4 + \mathbb{Z}$, $q \in \varepsilon_q/4 + \mathbb{Z}$, 
$\alpha \in \varepsilon_p(\delta/8 + \varepsilon/4) + \mathbb{Z}$, 
$\beta \in \varepsilon_q(\delta/8 - \varepsilon/4) + \mathbb{Z}$, 
where $\varepsilon_p, \, \varepsilon_q, \, \varepsilon = \pm 1$, 
$\delta = \pm1, \, \pm3$.    
\end{enumerate}     
\end{proposition} 
{\it Proof}. 
There are two cases: (I) either $p$ or $q$ is an integer; 
(II) neither $p$ nor $q$ is an integer.  
\par 
Case (I) is divided into four subcases: (I-1) either $(p, \alpha) \in \mathbb{Z}^2$    
or $(q, \beta) \in \mathbb{Z}^2$; (I-2) $p$, $q \in \mathbb{Z}$ 
and $\alpha$, $\beta \in \mathbb{R} \setminus \mathbb{Z}$; (I-3) $p \in \mathbb{Z}$, 
$q \in \mathbb{Q} \setminus \mathbb{Z}$, 
$\alpha \in \mathbb{R} \setminus \mathbb{Z}$; (I-4) $p \in \mathbb{Q} \setminus \mathbb{Z}$, 
$q \in \mathbb{Z}$, $\beta \in \mathbb{R} \setminus \mathbb{Z}$. 
In subcase (I-1) one has $\sigma_j = 0$ for every $j \in \mathbb{Z}$, which is 
just case (1) of the lemma.  
In subcase (I-2) one has $\sigma_j = \sin \pi\{\alpha \} \cdot \sin \pi\{\beta \}$ 
for every $j \in \mathbb{Z}$, which falls into case (2) of the lemma.   
In subcase (I-3) one has $\sigma_j = \sin \pi\{\alpha \} \cdot \sin \pi\{q j+ \beta \}$   
with nonzero $\sin \pi\{\alpha \}$ so that $\sigma_j$ is independent of 
$j$ if and only if $\sin\pi\{q j+ \beta \}$ is independent of $j$.   
Since $q \in \mathbb{Q} \setminus \mathbb{Z}$, it follows from Lemma \ref{lem:sin} that    
$q \in 1/2 + \mathbb{Z}$ and $\beta \equiv \pm 1/4 \mod 1$, which leads to 
case (3) of the lemma. 
In a similar manner subcase (I-4) leads to case (4) of the lemma. 
Thus case (I) is completed.  
\par
We proceed to case (II). 
Let $p = l_p/k_p$ and $q = l_q/k_q$ be the reduced 
representations of $p$ and $q$, where $k_p$, $k_q \ge 2$ and 
$(k_p, \, l_p) = (k_q, \, l_q) = 1$. 
If $k$ is the least common multiple of $k_p$ and $k_q$, then 
$(\{pi+ \alpha \}, \, \{qi+ \beta \}) = (\{pj+ \alpha \}, \, \{qj+ \beta \})$ 
precisely when $i \equiv j \mod k$. 
Thus the sequence $(\{pj+ \alpha \}, \, \{qj+ \beta \})$ takes exactly $k$ 
distinct values $(\{pj+ \alpha \}, \, \{qj+ \beta \})$ 
$(j = 0, \dots, k-1)$ as $j$ varies in $\mathbb{Z}$.  
In order for $\sigma_j$ to be independent of $j$ one must have 
\begin{equation} \label{eqn:sin2}
\sin \pi\{p j + \alpha \} \cdot \sin \pi\{q j + \beta \} = 
\mbox{a constant} \qquad (j = 0, \dots, k-1). 
\end{equation}
Case (II) is divided into two subcases: (II-1) 
$k_p = k_q = k = 2$; (II-2) either $k_p \ge 3$ or $k_q \ge 3$.  
Note that in subcase (II-2) one must have either $k \ge 4$ or  
$k_p = k_q = k = 3$.  
\par
First we consider subcase (II-1). 
Since $(2, \, l_p) = (2, \, l_q) = 1$, the integers $l_p$ and $l_q$ 
must be odd, so that $p$, $q \in 1/2 + \mathbb{Z}$ and condition 
\eqref{eqn:sin2} becomes 
\begin{equation} \label{eqn:sin3} 
\sin\pi\{\alpha \} \cdot \sin\pi\{\beta \} = 
\sin\pi\{\alpha +1/2\} \cdot \sin\pi\{\beta + 1/2 \}. 
\end{equation}
Let $\alpha' = \{\alpha \}$, $\beta' = \{\beta \} \in [0, \, 1)$.  
If $0 \le \alpha', \, \beta' < 1/2$ then $\{\alpha +1/2\} = \alpha' + 1/2$ 
and $\{\beta + 1/2\} = \beta' + 1/2$, so \eqref{eqn:sin3} reads 
$\sin \pi \alpha' \cdot \sin \pi \beta' = \cos \pi \alpha' \cdot \cos \pi \beta'$, 
that is, $\cos \pi(\alpha'+\beta') = 0$ and hence $\alpha' + \beta' = 1/2$. 
If $0 \le \alpha' < 1/2 \le \beta' < 1$ then $\{\alpha +1/2\} = \alpha' + 1/2$ 
and $\{\beta + 1/2\} = \beta' - 1/2$, so \eqref{eqn:sin3} reads 
$\sin \pi \alpha' \cdot \sin \pi \beta' = - \cos \pi \alpha' \cdot \cos \pi \beta'$, 
that is, $\cos \pi(\beta'-\alpha') = 0$ and hence $\beta' - \alpha' = 1/2$. 
Similar reasoning shows that $\alpha'-\beta'=1/2$ if 
$0 \le \beta' < 1/2 \le \alpha' < 1$; and  $\alpha' + \beta' = 3/2$ 
if $1/2 \le \alpha', \, \beta' < 1$.  
Summing up, one has $\alpha + \beta \in 1/2 + \mathbb{Z}$ or 
$\alpha- \beta \in 1/2 + \mathbb{Z}$, 
which falls into case (5) of the lemma. 
\par
Secondly we consider subcase (II-2). 
We begin by the following claim. 
\par\medskip
{\bf Claim 1}. In subcase (II-2) the constant in 
condition \eqref{eqn:sin2} must be nonzero.    
\par\medskip 
Assume the contrary that it is zero. 
First suppose that $k \ge 4$. 
Putting $j = 0$ in condition \eqref{eqn:sin2} yields 
$\sin\pi\{\alpha \}\cdot\sin\pi\{\beta \} = 0$, 
which means that either $\alpha$ or $\beta$ is an integer.   
By symmetry we may assume that $\alpha$ is an integer, in which case 
$p + \alpha$ is not an integer and so $\sin\pi\{p+ \alpha \}$ is nonzero.  
Putting $j = 1$ in \eqref{eqn:sin2} yields $\sin\pi\{p+ \alpha \} \cdot 
\sin\pi\{q+ \beta \} = 0$ and so $\sin\pi\{q+ \beta \} = 0$, that is, 
$q+ \beta \in \mathbb{Z}$. 
Then $2q+ \beta$ is not an integer and so $\sin\pi\{2q+ \beta \}$ is nonzero. 
Putting $j = 2$ in \eqref{eqn:sin2} yields $\sin\pi\{2p+ \alpha \} \cdot 
\sin\pi\{2q+ \beta \} = 0$ and so $\sin\pi\{2p+ \alpha \} = 0$, that is, 
$2p+ \alpha \in \mathbb{Z}$. 
Then $3p+ \alpha$ is not an integer and so $\sin\pi\{3p+ \alpha \}$ is nonzero. 
Putting $j = 3$ in \eqref{eqn:sin2} yields $\sin\pi\{3p+ \alpha \} \cdot 
\sin\pi\{3q+ \beta \} = 0$ and so $\sin\pi\{3q+ \beta \} = 0$, that is, 
$3q+ \beta \in \mathbb{Z}$.  
Therefore one has $\alpha$, $2p+ \alpha$, $q+ \beta$, $3q+ \beta \in \mathbb{Z}$, 
and thus $2p$, $2q \in \mathbb{Z}$, that is, $p$, $q \in 1/2 + \mathbb{Z}$.  
This means that $k_p = k_q = 2$, which is a contradiction. 
Next suppose that $k_p = k_q = k = 3$. 
In this case the same argument as above with $j = 0,1,2$ remains 
true and shows that $\alpha$ is an integer but $3p+ \alpha$ is not. 
This is a contradiction because if $\alpha$ is an 
integer then so is $3p+ \alpha$ by $k_p = 3$.    
\hfill $\Box$ \par\medskip
{\bf Claim 2}. In subcase (II-2) one must have 
$k_p = k_q = k \ge 3$.   
\par\medskip 
Let $m_p = k/k_p$ and $m_q = k/k_q$, that is, $k = k_p m_p = k_q m_q$.  
Note that $(m_p, \, m_q) = 1$. 
Putting $j = k_p \mu$ with $\mu = 0, \dots, m_p-1$, we observe that 
$\sin\pi\{p j + \alpha \} = \sin\pi\{l_p \mu + \alpha \} = 
\sin\pi\{ \alpha \}$ is independent of $\mu$.  
Thus putting $j = k_p \mu$ with $\mu = 0, \dots, m_p-1$ in condition 
\eqref{eqn:sin2} implies that $\sin \pi\{ \alpha \} \cdot 
\sin \pi\{(q k_p) \mu + \beta \}$ is independent of $\mu = 0, \dots, m_p-1$. 
Since $\sin \pi\{\alpha \}$ is nonzero by Claim 1, 
$\sin \pi\{(q k_p) \mu + \beta \}$ is also independent of $\mu = 0, \dots, m_p-1$.   
This forces $m_p = 1$ or $m_p = 2$ as in the proof of Lemma \ref{lem:sin}. 
Similarly one must have $m_q = 1$ or $m_q = 2$. 
So there are at most three possibilities: (a) $m_p = m_q = 1$; 
(b) $m_p = 2$, $m_q = 1$; (c) $m_p = 1$, $m_q = 2$, because 
$m_p = m_q = 2$ is forbidden by $(m_p, \, m_q) = 1$. 
We show that case (b) is impossible.  
Indeed, in this case, $k = 2 k_p = k_q$ and hence $(k_q, \, l_q) 
= (2 k_p, \, l_q) = 1$ implies that $(k_p, \, l_q) = 1$ and $l_q$ is odd. 
For each $j = 0, \dots, k-1$, write $j = k_p \mu + \nu$ with 
$\mu = 0, 1$ and $\nu = 0, \dots, k_p-1$. 
Condition \eqref{eqn:sin2} says that $\sin\pi\{p i + \alpha \} \cdot
\sin\pi\{q i + \beta \} = \sin\pi\{l_p \mu + p \nu + \alpha \}\cdot
\sin\pi\{l_q \mu/2 + q \nu + \beta \} = \sin\pi\{p \nu + \alpha \} 
\cdot \sin\pi\{\mu/2 + q \nu + \beta \}$ is nonzero and independent of 
$\mu = 0, 1$ and $\nu = 0, \dots, k_p-1$. 
This in particular yields $\sin\pi\{p \nu + \alpha \} \cdot 
\sin\pi\{q \nu + \beta \} = \sin\pi\{p \nu + \alpha \} \cdot 
\sin\pi\{1/2 + q \nu + \beta \}$ and so  
$\sin\pi\{q \nu + \beta \} = \sin\pi\{1/2 + q \nu + \beta \}$ and hence 
$q \nu + \beta \equiv \pm 1/4 \mod 1$ for every $\nu = 0, \dots, k_p-1$. 
Putting $\nu = 0, 1$, one has $q = (q+ \beta) - \beta \equiv 
(\pm 1/4)-(\pm 1/4) \equiv 0, \, 1/2 \mod 1$. 
But this is impossible because $q$ has the reduced representation 
$q = l_q/k_q$ with $k_q = 2 k_p \ge 4$. 
Similarly case (c) is also impossible.  
Thus case (II-2) must fall into case (a) where $k_p = k_q = k \ge 3$. 
\hfill $\Box$ \par\medskip
By Claim 2 the rational numbers $p$ and $q$ have reduced representations 
$p = l_p/k$ and $q = l_q/k$ with common denominator $k \ge 3$, where 
$(k, \, l_p) = (k, \, l_q ) = 1$. 
The analysis of this final case requires several steps and occupies 
all the rest of the proof. 
Consider the $k$ points 
\begin{equation} \label{eqn:k-points}  
(\{p j + \alpha \}, \, \{q j + \beta \}) \in [0, \, 1) \times [0, \, 1) 
\qquad (j = 0, \dots, k-1). 
\end{equation} 
They lie on a level set of $\varphi(x, y) := \sin(\pi x) \cdot \sin(\pi y)$, 
whose height $h$ must be nonzero by Claim 1. 
\begin{figure}[t]
\begin{minipage}{0.5\hsize}
\begin{center}
\includegraphics[width=60mm,clip]{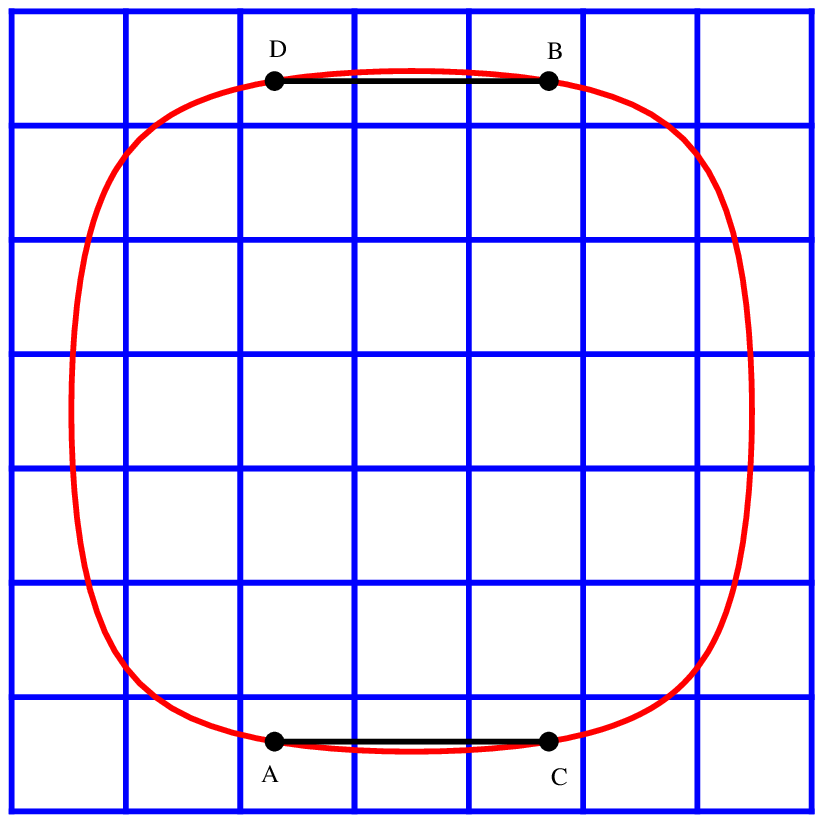} 
\end{center}
\vspace{-3mm} 
\caption{Level curve $\varphi(x, y) = h$.}
\label{fig:sinsin7}
\end{minipage}
\begin{minipage}{0.47\hsize}
\begin{center}
\includegraphics[width=60mm,clip]{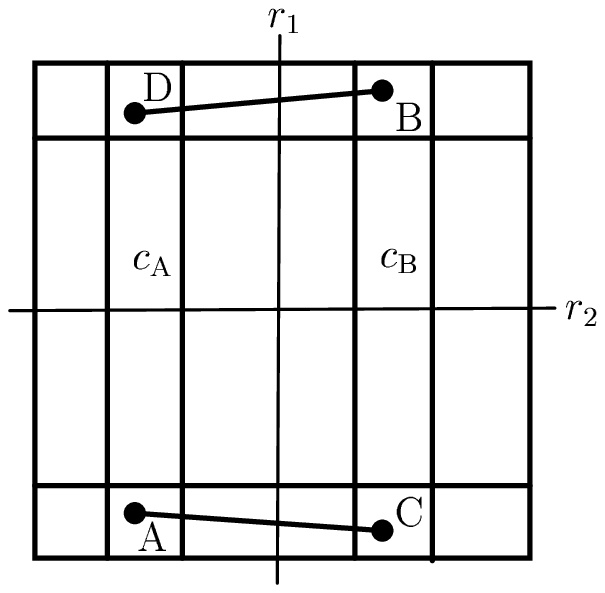} 
\vspace{6mm} 
\end{center}
\vspace{-7mm} 
\caption{Reflections of the square.}
\label{fig:reflections}
\end{minipage}
\end{figure}
The level set $\varphi(x,y) = h$ is a simple closed curve as in 
Figure \ref{fig:sinsin7} with its interior $\varphi(x,y) > h$ being  
a strictly convex domain. 
Since $(k, \, l_p) = (k, \, l_q) = 1$, there exist unique integers 
$L$, $M \in \{1, \dots, k-1\}$ coprime to $k$ such that 
$l_q \equiv L l_p$ and $l_p \equiv M l_q \mod k$. 
We can write  
\begin{equation} \label{eqn:ab}
\alpha = l_{\alpha}/k + \alpha', \qquad \beta = l_{\beta}/k + \beta', 
\qquad l_{\alpha}, \, l_{\beta} \in \mathbb{Z}, \qquad 
0 < \alpha', \, \beta' < 1/k,    
\end{equation}
where a priori the last condition should be $0 \le \alpha', \, \beta' < 1/k$ 
but equality is ruled out by Claim 1. 
Let $i_{\alpha} \in \{0,\dots,k-1\}$ be the unique integer such that 
$i_{\alpha} \equiv l_{\alpha} - M l_{\beta} \mod k$. 
\par\medskip
{\bf Claim 3}. The $k$ points in \eqref{eqn:k-points} can be rearranged as      
\begin{equation} \label{eqn:L} 
(i/k + \alpha', \, \{ L \,(i-i_{\alpha})/k + \beta' \}) \qquad  
(i = 0, \dots, k-1). 
\end{equation}
\par
Observe that 
$(\{p j + \alpha \}, \, \{q j + \beta \}) = (\{(l_p j + l_{\alpha})/k + \alpha' \}, \, 
\{(l_q j + l_{\beta})/k + \beta' \})$. 
We work with the quotient ring $\mathbb{Z}/k\mathbb{Z}$ and its unit group 
$(\mathbb{Z}/k\mathbb{Z})^{\times}$. 
Note that $L = l_q l_p^{-1}$ and $M = L^{-1} = l_p l_q^{-1}$ in 
$(\mathbb{Z}/k\mathbb{Z})^{\times}$.  
If $i = l_p j + l_{\alpha}$ in $\mathbb{Z}/k\mathbb{Z}$ then 
$j = l_p^{-1} (i-l_{\alpha})$ and so $l_q j + l_{\beta} = 
l_q  l_p^{-1} (i-l_{\alpha}) + l_{\beta} = 
(l_q  l_p^{-1})(i - l_{\alpha} + (l_p l_q^{-1})l_{\beta}) = 
L(i-l_{\alpha}+ M l_{\beta}) = L(i-i_{\alpha})$ 
in $\mathbb{Z}/k\mathbb{Z}$. 
Thus $(\{p j + \alpha \}, \, \{q j + \beta \})$ $(j = 0,\dots,k-1)$ can be 
rearranged as $(\{i/k + \alpha' \}, \, \{L(i-i_{\alpha})/k + \beta' \})$ 
$(i = 0,\dots,k-1)$. 
Since $\{i/k + \alpha'\} = i/k + \alpha'$ by $i = 0, \dots, k-1$ and 
$0 < \alpha' < 1/k$, we have formula \eqref{eqn:L}. \hfill $\Box$ 
\par\medskip
If the square $[0, \, 1) \times [0, \, 1)$ is equi-partitioned 
into $k$ columns and $k$ rows as in Figure \ref{fig:sinsin7}, then 
rearrangement \eqref{eqn:L} tells us that each column contains exactly one 
point from \eqref{eqn:k-points}, where the index $i$ in \eqref{eqn:L} 
corresponds to the $(i+1)$-th column of the square.   
Similarly, a row version of \eqref{eqn:L} implies that each 
row also contains exactly one point from \eqref{eqn:k-points}.  
By Claim 1 each point of \eqref{eqn:k-points} must be in the interior 
of a small square (or a box) created by the $k$-by-$k$ partition.   
\par
Notice that the function $\varphi(x,y)$ is invariant under two reflections 
$r_1 : (x, \, y) \mapsto (1-x, \, y)$ and $r_2 : (x, \, y) \mapsto (x, \, 1-y)$ 
(see Figure \ref{fig:reflections}).   
On the other hand we have $1 - \{ t \} = \{ -t \}$ whenever $t \in \mathbb{R} \setminus \mathbb{Z}$, 
which yields $1-\{p j+ \alpha \} = \{-p j- \alpha \}$ and $1-\{q j+ \beta \} = \{-q j- \beta \}$, 
since neither $p j + \alpha$ nor $q j + \beta$ is an integer by Claim 1. 
Thus these reflections induce two symmetries  
\begin{equation} \label{eqn:symmetry}
r_1 : (p,q; \alpha, \beta) \mapsto (-p,q; -\alpha, \beta), \qquad 
r_2 : (p,q; \alpha, \beta) \mapsto (p,-q; \alpha, -\beta).   
\end{equation} 
among solutions to our current problem, which may be settled 
only up to these symmetries.  
\par
Let $\mathrm{A}$ be the point from \eqref{eqn:k-points} that lie in 
the bottom row and $c_{\mathrm{A}}$ the column that contains $\mathrm{A}$. 
Similarly let $\mathrm{B}$ be the point from \eqref{eqn:k-points} 
that lie in the top row and $c_{\mathrm{B}}$ the column that contains 
$\mathrm{B}$ (see Figure \ref{fig:reflections}). 
For two columns $c_1$ and $c_2$ we write $c_1 \prec c_2$ 
if $c_1$ is located to the left of $c_2$. 
We may assume $c_{\mathrm{A}} \prec c_{\mathrm{B}}$, for otherwise $c_{\mathrm{A}}$ and $c_{\mathrm{B}}$ 
can be exchanged by the reflection $r_1$. 
\par\medskip
{\bf Claim 4}. If $c_{\mathrm{A}} \prec c_{\mathrm{B}}$ then $c_{\mathrm{B}}$ must be 
the right neighbor of $c_{\mathrm{A}}$, $L = k-1$ in \eqref{eqn:L}, and  
\begin{equation} \label{eqn:ell}
l_p + l_q \equiv 0, \qquad l_{\alpha} + l_{\beta} 
\equiv i_{\alpha} \quad \mod k. 
\end{equation}    
\par
Suppose the contrary that there is an intermediate column 
$c$ between $c_{\mathrm{A}}$ and $c_{\mathrm{B}}$.   
Let $\mathrm{C}$ resp. $\mathrm{D}$ be the $r_2$-reflection of the point $\mathrm{B}$ resp. 
$\mathrm{A}$ (see Figure \ref{fig:reflections}). 
Since the interior $\varphi(x,y) > h$ of the level curve $\varphi(x,y) 
= h$ is strictly convex, the level arc $\mathrm{A}\mathrm{C}$ is 
lower than  the line segment $\mathrm{A}\mathrm{C}$ while the level arc $\mathrm{D}\mathrm{B}$ 
is upper than the line segment $\mathrm{D}\mathrm{B}$ (see Figure \ref{fig:sinsin7}). 
Thus the level curve $\varphi(x,y) = h$ meets the column $c$  
in its top and bottom boxes only. 
So the unique point from \eqref{eqn:k-points} lying in the 
column $c$ must be either in the top box or in the 
bottom box.  
But this is impossible because the top and bottom rows are already 
occupied by the points $\mathrm{B}$ and $\mathrm{A}$ respectively. 
Therefore the columns $c_{\mathrm{A}}$ and $c_{\mathrm{B}}$ must be consecutive, that is, 
$c_{\mathrm{B}}$ must be the right neighbor of $c_{\mathrm{A}}$ since 
$c_{\mathrm{A}} \prec c_{\mathrm{B}}$. 
Next we show that $L = k-1$.  
In view of formula \eqref{eqn:L} with $i = i_{\alpha}$ the point $\mathrm{A}$ 
has coordinates $(i_{\alpha}/k + \alpha', \, \beta')$, where 
$i_{\alpha} \in \{0,\dots,k-2\}$ because $c_{\mathrm{A}}$ is not the rightmost column. 
The point $\mathrm{B}$ is then given by formula \eqref{eqn:L} with 
$i = i_{\alpha} + 1$, that is, by $((i_{\alpha}+1)/k + \alpha', \, \{L/k + \beta'\})$, 
where $\{L/k + \beta'\} = L/k + \beta'$ since $L \in \{1,\dots,k-1\}$ 
and $0 < \beta' < 1/k$. 
In order for $\mathrm{B}$ to be in the top row, one must have $L = k-1$. 
Accordingly, $L \equiv M \equiv -1 \mod k$ so that 
$l_q \equiv L l_p \equiv -l_p$ and $i_{\alpha} \equiv l_{\alpha} - M l_{\beta} 
\equiv l_{\alpha} + l_{\beta} \mod k$.  \hfill $\Box$ 
\par\medskip
{\bf Claim 5}. The case $k \ge 5$ is impossible and so one has 
at most $k = 3$, $4$. 
\par
\begin{figure}[t]
\begin{center}
\includegraphics[width=60mm,clip]{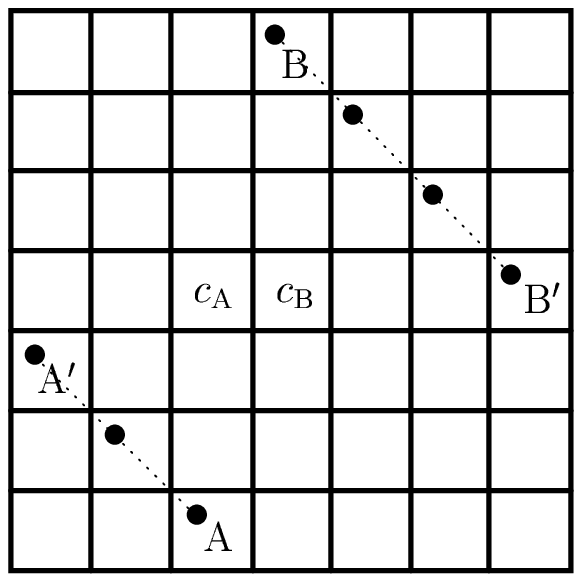} 
\end{center}
\caption{A configuration of the $k$ points (when $k = 7$).}
\label{fig:config}
\end{figure}

\medskip
Claim 4 and formula \eqref{eqn:L} show that if $c_{\mathrm{A}} \prec c_{\mathrm{B}}$  
then the $k$ points in \eqref{eqn:k-points} are given by 
\begin{equation} \label{eqn:config}
\begin{cases}
(i/k + \alpha', \, (i_{\alpha}-i)/k + \beta') & (i = 0, \dots, i_{\alpha}), \\
(i/k + \alpha', \, (i_{\alpha}-i)/k + \beta' + 1) & (i = i_{\alpha}+1, \dots, k-1),
\end{cases}
\end{equation}
which comprises two straight strings $\mathrm{A}'\mathrm{A}$ and 
$\mathrm{B}\mathrm{B}'$ of slope $-1$ as in Figure \ref{fig:config}. 
Each string cannot contain more than two points, because the 
interior $\varphi(x,y) > h$ of the level curve $\varphi(x,y) = h$ 
is strictly convex so that three collinear points cannot lie 
on the level curve.  
Thus if $k_{\mathrm{A}}$ resp. $k_{\mathrm{B}}$ stands for the number of points 
on the string $\mathrm{A}'\mathrm{A}$ resp. $\mathrm{B}\mathrm{B}'$, then one must have 
$k_{\mathrm{A}} \le 2$ and $k_{\mathrm{B}} \le 2$ so that $k = k_{\mathrm{A}} + k_{\mathrm{B}} \le 4$. 
(Figure \ref{fig:config} exhibits a configuration that is impossible.)   
\hfill $\Box$
\par\medskip
Note that $k_{\mathrm{A}} = i_{\alpha} + 1$ since the point $\mathrm{A}$ belongs to  
the $(i_{\alpha}+1)$-th column.  
Observe that the transformation $r_1 \circ r_2$ brings 
$\alpha' \mapsto 1/k-\alpha'$, $\beta' \mapsto 1/k-\beta'$ and 
$i_{\alpha} \mapsto k-2-i_{\alpha}$ without violating the condition 
$c_{\mathrm{A}} \prec c_{\mathrm{B}}$. 
So the symmetry allows us to assume $i_{\alpha} \le k-2-i_{\alpha}$, that is,  
$0 \le i_{\alpha} \le [k/2]-1$ in addition to $c_{\mathrm{A}} \prec c_{\mathrm{B}}$; 
when $k = 3$ this means that we may assume $i_{\alpha} = 0$.     
\par\medskip
{\bf Claim 6}. When $k =3$, if $c_{\mathrm{A}} \prec c_{\mathrm{B}}$ and 
$i_{\alpha} = 0$ then 
\begin{equation} \label{eqn:k=3}
p \in \frac{\varepsilon}{3} + \mathbb{Z}, \quad 
p \in -\frac{\varepsilon}{3} + \mathbb{Z}, \quad 
\alpha \in \frac{\delta}{3} + \kappa + \mathbb{Z}, \quad 
\beta \in -\frac{\delta}{3} + \kappa + \mathbb{Z}, 
\end{equation}
where $\varepsilon = \pm 1$, $\delta = 0, \, \pm 1$ and $\kappa$ is 
the irrational constant defined by formula \eqref{eqn:kappa}. 
\par
\begin{figure}[t]
\begin{minipage}{0.45\hsize}
\begin{center}
\includegraphics[width=45mm,clip]{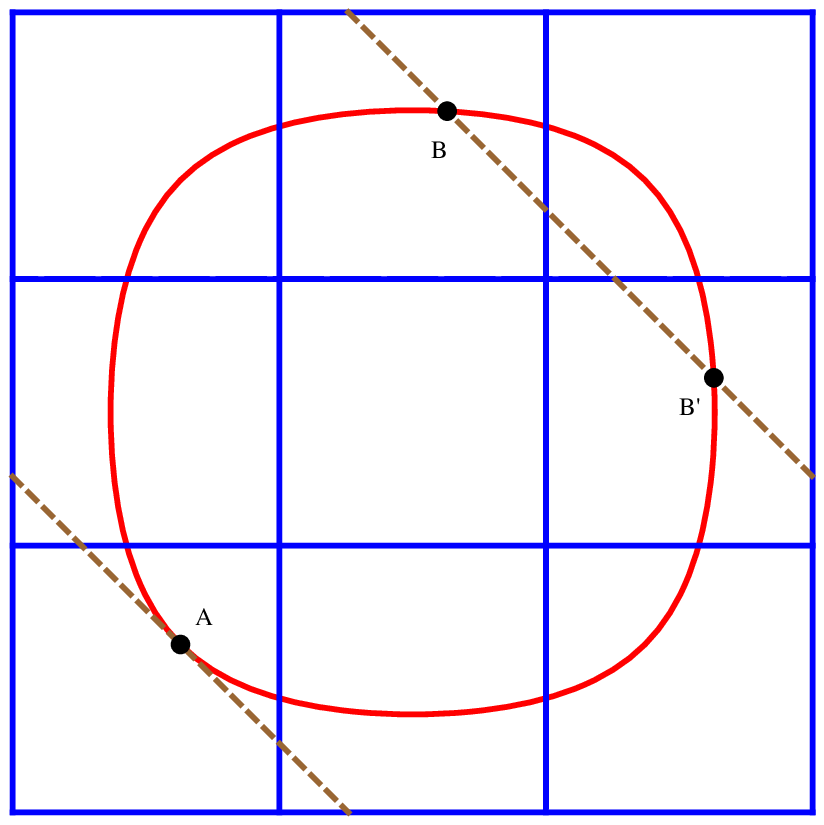} 
\end{center}
\vspace{-7mm} 
\caption{$k = 3$.}
\label{fig:sinsin3}
\end{minipage}
\begin{minipage}{0.45\hsize}
\begin{center}
\includegraphics[width=45mm,clip]{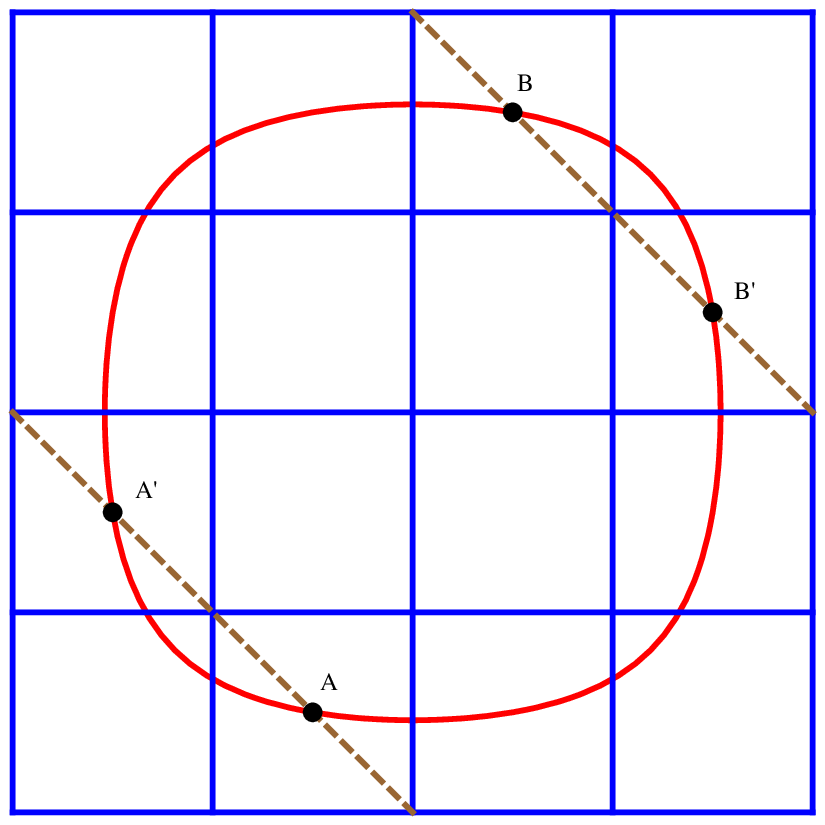} 
\end{center}
\vspace{-7mm} 
\caption{$k = 4$.}
\label{fig:sinsin4}
\end{minipage}
\end{figure}

\medskip
If $k = 3$ and $i_{\alpha} = 0$, the three points in formula 
\eqref{eqn:config} are given by $\mathrm{A} (\alpha', \, \beta')$, 
$\mathrm{B} (\alpha'+1/3, \, \beta'+2/3)$ and $\mathrm{B}' (\alpha'+2/3, \, \beta'+1/3)$, 
where $0< \alpha', \, \beta' < 1/3$.  
Condition \eqref{eqn:sin2} now reads  
\[
\sin\pi\alpha' \cdot \sin\pi\beta' = 
\sin\pi(\alpha'+1/3)\cdot \sin\pi(\beta'+2/3) = 
\sin\pi(\alpha'+2/3)\cdot \sin\pi(\beta'+1/3).  
\] 
Addition formula for sine recasts the second equality to  
$\tan \pi \alpha' = \tan \pi \beta'$ and so $\alpha' = \beta'$. 
The first equality then becomes $\tan^2\pi\alpha' = 3/5$ yielding 
$\alpha' = \beta' = \frac{1}{\pi} \arctan \sqrt{3/5}$, which is 
the constant $\kappa$ defined by \eqref{eqn:kappa}.  
The configuration of the three points is as in 
Figure \ref{fig:sinsin3}.  
Formula \eqref{eqn:k=3} follows from $p = l_p/3$, $q = l_q/3$ with 
$(3, \, l_p) = (3, \, l_q) = 1$ and formulas \eqref{eqn:ab} and 
\eqref{eqn:ell} with $i_{\alpha} = 0$.    
\hfill $\Box$
\par\medskip
{\bf Claim 7}. When $k =4$, if $c_{\mathrm{A}} \prec c_{\mathrm{B}}$ then 
\begin{equation} \label{eqn:k=4}
p \in \frac{\varepsilon}{4} + \mathbb{Z}, \quad 
q \in -\frac{\varepsilon}{4} + \mathbb{Z}, \quad 
\alpha \in \frac{2+\delta}{8} + \mathbb{Z}, \quad 
\beta \in \frac{2-\delta}{8} + \mathbb{Z} \quad 
(\varepsilon = \pm 1, \, \delta = \pm 1, \, \pm 3). 
\end{equation}
\par
We must have $k_{\mathrm{A}} = k_{\mathrm{B}} = 2$ and so $i_{\alpha} = 1$, because 
$k_{\mathrm{A}} \le 2$, $k_{\mathrm{B}} \le 2$ and $k_{\mathrm{A}} + k_{\mathrm{B}} = k = 4$.  
If $k = 4$ and $i_{\alpha} = 1$, the four points in \eqref{eqn:config} are 
given by $\mathrm{A} (\alpha'+1/4, \, \beta')$, $\mathrm{B} (\alpha'+2/4, \, \beta'+3/4)$, 
$\mathrm{A}' (\alpha', \, \beta'+1/4)$ and $\mathrm{B}' (\alpha'+3/4, \, \beta'+2/4)$, 
where $0< \alpha', \, \beta' < 1/4$.  
Condition \eqref{eqn:sin2} now implies    
\[
\sin\pi\alpha' \cdot \sin\pi(\beta'+1/4) = 
\sin\pi(\alpha'+1/4) \cdot \sin\pi\beta' = 
\sin\pi(\alpha'+2/4)\cdot \sin\pi(\beta'+3/4),   
\] 
where the first equality comes from $\varphi(\mathrm{A}') = \varphi(\mathrm{A})$ 
while the second one from $\varphi(\mathrm{A}) = \varphi(\mathrm{B})$.   
Addition formula for sine recasts the first equality to  
$\tan \pi \alpha' = \tan \pi \beta'$ and so $\alpha' = \beta'$. 
The second equality then becomes $\tan 2\pi\alpha' = 1$ yielding 
$\alpha' = \beta' = 1/8$.  
The configuration of the four points is indicated in Figure \ref{fig:sinsin4}. 
Formula \eqref{eqn:k=4} follows from $p = l_p/4$, $q = l_q/4$ with 
$(4, \, l_p) = (4, \, l_q) = 1$ and formulas \eqref{eqn:ab} and 
\eqref{eqn:ell} with $i_{\alpha} = 1$. \hfill $\Box$ \par\medskip
We are now in a position to establish cases (6) and (7) of Proposition 
\ref{prop:sin2}.  
Replace $\delta$ by $\varepsilon \delta$ in the third and fourth 
formulas in \eqref{eqn:k=3}.  
Using $\mathbb{Z}_2 \times \mathbb{Z}_2$ symmetries in \eqref{eqn:symmetry}, multiply 
the first and third formulas by $\varepsilon_p \varepsilon$, while the 
second and fourth formulas by $-\varepsilon_q \varepsilon$, where 
$\varepsilon_p, \varepsilon_q = \pm1$.  
We then arrive at case (6). 
Similarly, replace $\delta$ by $\varepsilon \delta$ in the third and 
fourth formulas in \eqref{eqn:k=4}.   
Multiply the first and third formulas by $\varepsilon_p \varepsilon$, 
while the second and fourth formulas by $-\varepsilon_q \varepsilon$, 
where $\varepsilon_p$, $\varepsilon_q = \pm1$.  
We then arrive at case (7). 
The proof of Proposition \ref{prop:sin2} is complete. \hfill $\Box$     
\section{Parity} \label{sec:parity}
Proposition \ref{prop:kronecker1} was able to confine the 
possibilities of $(p,q; \alpha, \beta)$ substantially as in Proposition 
\ref{prop:sin2}. 
We shall discuss how further Proposition \ref{prop:kronecker2} 
can reduce those possibilities.  
\begin{lemma} \label{lem:parity1} 
Let $r \in \mathbb{N}$, $p \in \mathbb{Q}^{\times}$ and $\alpha \in \mathbb{R}$. 
The sequence $\sin \pi\{pj+ \alpha \}$ takes a nonzero constant value 
and the integer sequence $[pj+\alpha]+rj$ has a constant parity independent 
of $j \in \mathbb{Z}$, if and only if $p$ is an integer with the same parity 
as $r$ and $\alpha \in \mathbb{R} \setminus \mathbb{Z}$.  
\end{lemma}
{\it Proof}. We are either in case (1) or in case (2) of 
Lemma \ref{lem:sin}. 
We begin by case (1).  
We must have $\alpha \in \mathbb{R} \setminus \mathbb{Z}$ because 
$\sin\pi\{\alpha \}$ should be nonzero. 
Since $p \in \mathbb{Z}$ we have $[pj+\alpha]+rj = (p+r)j + [\alpha]$, 
which has a constant parity if and only if $p+r$ is even, that is, 
$p$ and $r$ have the same parity. 
We proceed to case (2).  
Note that $2p$ is an odd integer since $p \in 1/2 + \mathbb{Z}$.  
Replacing $j$ by $2j$ we have  $[2pj+\alpha]+2rj = (2p+2r)j + [\alpha] 
\equiv j + [\alpha] \mod 2$ since $2p$ is odd and $2r$ is even. 
Clearly $j +[\alpha]$ cannot have a constant parity, so that case 
(2) must be ruled out.  \hfill $\Box$  
\begin{proposition} \label{prop:parity2}
Let $r \in \mathbb{N}$, $p, \, q \in \mathbb{Q}^{\times}$ and 
$\alpha, \, \beta \in \mathbb{R}$. 
The sequence $\sin \pi\{p j+ \alpha \} \cdot \sin \pi\{q j+ \beta \}$ 
takes a nonzero constant value and the integer sequence 
$[p j+\alpha] + [q j+ \beta] +r j$ 
has a constant parity independent of $j \in \mathbb{Z}$, 
if and only if either condition $(\mathrm{A})$ or $(\mathrm{B})$ 
below is satisfied: 
\begin{enumerate}
\item[$(\mathrm{A})$] \ $p, \, q \in \mathbb{Z}$, $p+q+r$ is even, and 
$\alpha,\, \beta \in \mathbb{R} \setminus \mathbb{Z}$.  
\item[$(\mathrm{B})$] \ $p, \, q \in 1/2 + \mathbb{Z}$, 
$\alpha, \, \beta \in \mathbb{R} \setminus \mathbb{Z}$, 
and $\alpha -(-1)^{p+q+r} \beta \in 1/2 + \mathbb{Z}$.   
\end{enumerate}
\end{proposition}
{\it Proof}. We are in one of seven cases (1)--(7) of 
Proposition \ref{prop:sin2}, but case (1) is ruled out because the 
constant $\sin\pi\{\alpha\} \cdot \sin\pi\{\beta\}$ vanishes in this case. 
We shall make a case-by-case check as to whether the sequence 
$[p j+\alpha]+[q j+\beta]+r j$ can have a constant parity. 
In case (2), since $p$ and $q$ are integers, we have $[p j+\alpha] + 
[q j+\beta]+r j = (p+q+r)j + [\alpha] + [\beta]$, which has a constant parity 
if and only if $p+q+r$ is even. 
This is just the case of type $(\mathrm{A})$ in Proposition \ref{prop:parity2}.   
Next we consider case (3) where $p \in \mathbb{Z}$ and $q \in 1/2 + \mathbb{Z}$. 
Replacing $j$ by $2j$ in $[p j+ \alpha]+[q j+ \beta]+r j$ we have 
$[2pj+\alpha]+[2q j+\beta]+ 2r j = (2p+2q+2r)j + [\alpha] + [\beta] 
\equiv j + [\alpha] + [\beta] \mod 2$, since $2 p$ and $2 r$ are even while 
$2 q$ is odd. 
Clearly $j+[\alpha]+[\beta]$ cannot have a constant parity, so that case (3) 
must be ruled out. 
For the same reason case (4) is also ruled out. 
\par
We turn our attention to case (5). 
Since $p, \, q \in 1/2+\mathbb{Z}$, observe that $2p$ and $2q$ are odd integers 
and so $2p+2q$ is an even integer.  
Putting $j = 2k+l$ with $k \in \mathbb{Z}$ and $l = 0, 1$, we have    
$[p j+\alpha]+[q j+\beta]+r j = [2p k+p l+\alpha]+[2q k+q l+\beta] + 2r k +r l= 
(2 p+2 q+2 r)k + [p l+\alpha]+[q l+\beta] + r l \equiv [p l+\alpha]+[q l+\beta] 
+ r l \mod 2$ for every $k \in \mathbb{Z}$. 
So the sequence $[p j+\alpha]+[q j+\beta]+r j$ has a constant parity if and 
only if $[\alpha]+[\beta] \equiv [p+\alpha]+[q+\beta] + r \mod 2$, 
which is equivalent to
\begin{equation} \label{eqn:parity2}
[p]+[q] + r \equiv [1/2+\{\alpha \}]+[1/2+ \{\beta \}] \mod 2,  
\end{equation} 
because $[p+\alpha]+[q+\beta] = 
[p]+[\alpha]+[1/2+\{\alpha\}]+[q]+[\beta]+[1/2+\{\beta\}]$. 
In case (5) we have either $\alpha +\beta \in 1/2 + \mathbb{Z}$ or 
$\alpha - \beta \in 1/2 + \mathbb{Z}$, 
where the former condition implies (i) $\{\alpha \}+\{\beta \} = 1/2$ or 
(ii) $\{\alpha \}+\{\beta \} = 3/2$, while the latter implies (iii) 
$\{\alpha \}-\{\beta \} = 1/2$ or (iv) $\{\beta \}-\{\alpha \} = 1/2$.     
Note that $\{\alpha \}$ and $\{\beta \}$ must be positive, that is, 
$\alpha$, $\beta \in \mathbb{R} \setminus \mathbb{Z}$, since 
$\sin\pi\{\alpha \} \cdot \sin\pi\{\beta \}$ is nonzero. 
So $[1/2 + \{\alpha \}] = [1/2 + \{\beta \}] = 0$ in subcase (i);  
$[1/2 + \{\alpha \}] = [1/2 + \{\beta \}] = 1$ in subcase (ii); 
$[1/2 + \{\alpha \}] = 1$, $[1/2 + \{\beta \}] = 0$ in subcase (iii); 
$[1/2 + \{\alpha \}] = 0$, $[1/2 + \{\beta \}] = 1$ in subcase (iv). 
Thus condition \eqref{eqn:parity2} becomes 
$[p] + [q] + r \equiv 0 \mod 2$ in subcases (i) and (ii); and 
$[p] + [q] + r \equiv 1 \mod 2$ in subcases (iii) and (iv).   
This shows that $p+q+r$ is odd if $\alpha + \beta \in 1/2 + \mathbb{Z}$; 
and even if $\alpha - \beta \in 1/2 + \mathbb{Z}$, since $p+q = (1/2+[p]) 
+ (1/2+[q]) = [p] + [q] + 1$. 
In either case we have $\alpha -(-1)^{p+q+r} \beta \in 1/2 + \mathbb{Z}$, 
so that we are led to the case of type $(\mathrm{B})$ 
in Proposition \ref{prop:parity2}.   
\par
We proceed to case (6). 
In this case $p = \varepsilon_p/3 + p'$ and $q = \varepsilon_q/3 + q'$ 
with some $\varepsilon_p, \, \varepsilon_q = \pm 1$ and $p', \, q' \in \mathbb{Z}$. 
Note that $\varepsilon_p + \varepsilon_q$ is even in every choice of 
$\varepsilon_p$ and $\varepsilon_q$.  
Putting $j = 3k + l$ with $k \in \mathbb{Z}$ and $l = 0, \pm1$, we have 
$[p j+ \alpha]+[q j+ \beta]+r j = (\varepsilon_p + \varepsilon_q + 3p'+3q'+3r)k 
+ [l p+\alpha] + [l q+\beta] + l r \equiv (p'+q'+r)k + [p l+\alpha] + 
[q l+\beta] + r l \mod 2$. 
In order for this to have a common parity for every $k \in \mathbb{Z}$ 
the integer $p'+q'+r$ must be even. 
Under this condition, $[p+\alpha]+[q+\beta]+r = p'+q'+r+
[\varepsilon_p/3 + \alpha] + [\varepsilon_q/3 + \beta] \equiv 
[\varepsilon_p/3 + \alpha] + [\varepsilon_q/3 + \beta] \mod 2$ and 
$[-p+\alpha] + [-q+\beta]-r = -(p'+q'+r) +[-\varepsilon_p/3 + \alpha] + 
[-\varepsilon_q/3 + \beta] \equiv [-\varepsilon_p/3 + \alpha] + 
[-\varepsilon_q/3 + \beta] \mod 2$. 
Thus the parity of $[p l+ \alpha] + [q l+ \beta] + r l$ 
being constant for $l = 0, \pm1$ means 
\[
[\alpha]+[\beta] \equiv [\varepsilon_p/3 + \alpha] + 
[\varepsilon_q/3 + \beta] \equiv 
[-\varepsilon_p/3 + \alpha] + [-\varepsilon_q/3 + \beta] \mod 2. 
\]
Since $\alpha \in \varepsilon_p(\delta/3 + \varepsilon \kappa) + \mathbb{Z}$ 
and $\beta \in \varepsilon_q(\delta/3 - \varepsilon \kappa) + \mathbb{Z}$ with 
$\varepsilon = \pm1$, $\delta = 0, \pm1$, this becomes  
\[
\textstyle 
[\varepsilon_p(\frac{\delta}{3}+\varepsilon\kappa)]+
[\varepsilon_q(\frac{\delta}{3}-\varepsilon\kappa)] \equiv  
[\varepsilon_p(\frac{\delta+1}{3}+\varepsilon\kappa)]+
[\varepsilon_q(\frac{\delta+1}{3}-\varepsilon\kappa)] \equiv  
[\varepsilon_p(\frac{\delta-1}{3}+\varepsilon\kappa)]+
[\varepsilon_q(\frac{\delta-1}{3}-\varepsilon\kappa)]. 
\]
Note that if $t$ is not an integer then $[-t] = -[t]-1 \equiv [t] +1 
\mod 2$ and hence $[\epsilon t] \equiv [t] + \frac{1-\epsilon}{2} 
\mod 2$ for $\epsilon = \pm1$. 
This rule allows us to remove $\varepsilon_p$ and $\varepsilon_q$ 
from the above condition, that is,   
\begin{equation} \label{eqn:parity3}
\textstyle 
[\frac{\delta}{3}+\kappa]+
[\frac{\delta}{3}-\kappa] \equiv 
[\frac{\delta+1}{3}+\kappa]+
[\frac{\delta+1}{3}-\kappa] \equiv 
[\frac{\delta-1}{3}+\kappa]+
[\frac{\delta-1}{3}-\kappa] \mod 2, 
\end{equation}
where $\varepsilon$ can also be removed by symmetry. 
Since $\kappa = 0.209785 \cdots$ in formula \eqref{eqn:kappa}, 
congruences \eqref{eqn:parity3} become 
$-1 \equiv 0 \equiv -2$ if $\delta = 0$; $0 \equiv 0 \equiv -1$ if  
$\delta = 1$; and $-2 \equiv -1 \equiv -2$ if $\delta = -1$, 
respectively.  
In any case we have a contradiction and thus case (6) is ruled out. 
\par
Finally we consider case (7).  
In this case $p = \varepsilon_p/4 + p'$ and $q = \varepsilon_q/4 + q'$ 
with some $\varepsilon_p, \, \varepsilon_q = \pm 1$ and $p', \, q' \in \mathbb{Z}$. 
The parity of $[p j+\alpha] + [q j+\beta] + r j$ being equal for $j = 0, 2$ yields 
$[\alpha]+[\beta] \equiv [2p+\alpha] + [2q+\beta] + 2r = 
[\varepsilon_p/2 + \alpha] + [\varepsilon_q/2 + \beta] + 2(p'+q'+r) \equiv 
[\varepsilon_p/2 + \alpha] + [\varepsilon_q/2 + \beta] \mod 2$. 
Since $\alpha \in \varepsilon_p (\delta/8 + \varepsilon/4) + \mathbb{Z}$ and 
$\beta \in \varepsilon_q (\delta/8 - \varepsilon/4) + \mathbb{Z}$ with 
$\varepsilon = \pm1$ and $\delta = \pm1, \pm3$,  
\[
\textstyle
[\varepsilon_p(\frac{\delta}{8}+\frac{\varepsilon}{4})] + 
[\varepsilon_q(\frac{\delta}{8}-\frac{\varepsilon}{4})] \equiv
[\varepsilon_p(\frac{\delta}{8}+\frac{\varepsilon}{4} + \frac{1}{2})] + 
[\varepsilon_q(\frac{\delta}{8}-\frac{\varepsilon}{4} + \frac{1}{2})] 
\mod 2. 
\]
For the same reason as in the last paragraph, $\varepsilon_p$ and 
$\varepsilon_q$ can be removed from the above condition: 
\[
\textstyle
[\frac{\delta}{8}+\frac{1}{4}] + 
[\frac{\delta}{8}-\frac{1}{4}] \equiv
[\frac{\delta}{8}+\frac{1}{4} + \frac{1}{2}] + 
[\frac{\delta}{8}-\frac{1}{4} + \frac{1}{2}] 
\mod 2,  
\]
where $\varepsilon$ can also be removed by symmetry.  
But the right-hand side above is $[\frac{\delta}{8}+\frac{3}{4}] + 
[\frac{\delta}{8}+\frac{1}{4}] = [\frac{\delta}{8}+\frac{1}{4}] + 
[\frac{\delta}{8}-\frac{1}{4}]+1$, which leads to a contradiction 
$0 \equiv 1 \mod 2$. 
Thus case (7) is ruled out.   
\hfill $\Box$ \par\medskip 
Note that Theorem \ref{thm:integer} follows from Propositions 
\ref{prop:kronecker1}, \ref{prop:kronecker2}, \ref{prop:parity2} 
and Lemma \ref{lem:B}. 
\section{The Method of Contiguous Relations} \label{sec:contig}
If $p$, $q$, $r$ are integers then contiguous relations of Gauss 
lead to a general three-term relation \eqref{eqn:contig-F} and 
a specialization \eqref{eqn:contig-f} of it evaluated at 
$(\alpha,\beta;\gamma;z) = (p w+a, q w+b; r w; x)$. 
In this section we shall see how these formulas contribute to the 
discussions of Problems \ref{prob:contiguous1} and \ref{prob:necessary}.    
\subsection{Rational Independence} \label{subsec:indep}
To settle Problem \ref{prob:contiguous1} in region \eqref{eqn:cross} 
and prove Theorem \ref{thm:contiguous1}, we begin by the following. 
\begin{lemma} \label{lem:indep} 
If $\lambda = (p,q,r;a,b;x)$ is a non-elementary solution to Problem 
$\ref{prob:ocf}$ in region \eqref{eqn:cross}, then $f(w)$ and 
$\tilde{f}(w)$ in \eqref{eqn:f} and \eqref{eqn:f-tilde}are 
linearly independent over the rational function field $\mathbb{C}(w)$. 
\end{lemma}
{\it Proof}. 
If $\gamma$ is not an integer then the Gauss hypergeometric equation 
admits  
\[
u_1 := {}_2F_1(\alpha, \beta; \gamma; z), \quad  
u_2 := z^{1-\gamma} {}_2F_1(\alpha-\gamma+1, \beta-\gamma+1; 2-\gamma; z)
\]
as a fundamental set of solutions, whose Wronskian is given by 
\[
W := u_1 u_2' - u_1' u_2 = (1-\gamma) z^{-\gamma}(1-z)^{\gamma-\alpha-\beta-1}, 
\]
where $u' = du/dz$.  
From formulas (20) and (22) of Erd\'elyi \cite[Chapter II, \S 2.8]{Erdelyi},    
\[
u_1' = \frac{\alpha \beta}{\gamma} \, {}_2F_1(\alpha+1, \beta+1; \gamma+1; z), 
\quad 
u_2' = (1-\gamma) z^{-\gamma} \, {}_2F_1(\alpha-\gamma+1, \beta-\gamma+1; 1-\gamma; z).   
\]
Substituting these into the Wronskian formula above one has 
\begin{equation} \label{eqn:quad} 
\begin{split}
& (1-z)^{\gamma-\alpha-\beta-1} = {}_2F_1(\alpha, \beta; \gamma; z) \, 
{}_2F_1(\alpha-\gamma+1, \beta-\gamma+1; 1-\gamma; z) + \\ 
& \frac{\alpha \beta z}{\gamma(\gamma-1)} \, {}_2F_1(\alpha+1, \beta+1; \gamma+1; z) \, 
{}_2F_1(\alpha-\gamma+1, \beta-\gamma+1; 2-\gamma; z).    
\end{split}   
\end{equation}
\par
Since $\lambda$ is a non-elementary solution to Problem \ref{prob:ocf}, 
$f(w)$ has a GPF \eqref{eqn:f=g} with $m \ge 1$ by Propositions 
\ref{prop:f=g} and \ref{prop:pole2}. 
If $f(w)$ and $\tilde{f}(w)$ were linearly dependent over $\mathbb{C}(w)$, then   
there would be a rational function $T(w)$ such that $\tilde{f}(w) = T(w) f(w)$.  
Putting $\alpha = p w+a$, $\beta = q w+b$, $\gamma = r w$ and $z = x$ into 
formula \eqref{eqn:quad} yields $f(w) f_1(w) = (1-x)^{(r-p-q)w-a-b-1}$, where  
\begin{align*} 
f_1(w) &:= {}_2F_1((p-r)w+a+1, \, (q-r)w+b+1; \, 1-r w; \, x) \\ 
       & \quad + \frac{(p w+a)(q w+b) x}{r w(r w-1)} \cdot T(w) \cdot 
{}_2F_1((p-r)w+a+1, \, (q-r)w+b+1; \, 2-r w; \, x). 
\end{align*} 
Take a negative real number $R_5$ so that all poles of $T(w)$ and 
$S(w)$ are in the right half-plane $\mathrm{Re}(z) > R_5$, where 
$S(w)$ is the rational function in formula \eqref{eqn:f=g}.  
Since $r$ is positive, $f_1(w)$ is holomorphic on the left 
half-plane $\mathrm{Re}(z) < R_5$.  
Choose a positive integer $j$ so that $-j-v_1 < R_5$. 
Then $f(w)$ has a zero at $w = -j-v_1$ while $f_1(w)$ is 
holomorphic at this point. 
Therefore, $0 = f(-j-v_1) f_1(-j-v_1) = (1-x)^{-(r-p-q)(j+v_1)-a-b-1} 
\neq 0$, which is a contradiction.    
\hfill $\Box$ 
\begin{proposition} \label{prop:ocr}
If $\lambda = (p,q,r;a,b,x)$ is a non-elementary solution to Problem 
$\ref{prob:ocf}$ in region \eqref{eqn:cross} with $p,q,r \in \mathbb{Z}$,  
then it comes from contiguous relations. 
In particular Theorem $\ref{thm:contiguous1}$ holds.    
\end{proposition}
{\it Proof}. 
Since $\lambda$ is a solution to Problem \ref{prob:ocf}, there exists  
a rational function $\widetilde{R}(w)$ such that $f(w+1) = 
\widetilde{R}(w) f(w)$. 
Subtracting this from three-term relation \eqref{eqn:contig-f} 
yields a linear relation 
$\{R(w)-\widetilde{R}(w)\} f(w) + Q(w) \tilde{f}(w) = 0$ 
over $\mathbb{C}(w)$. 
By Lemma \ref{lem:indep} one must have $R(w)-\widetilde{R}(w) = Q(w) = 0$ 
in $\mathbb{C}(w)$, so that three-term relation \eqref{eqn:contig-f} degenerates 
to a two-term one \eqref{eqn:ocf}. \hfill $\Box$ 
\subsection{Contiguous Relations in Matrix Form} \label{subsec:c-r}
It is convenient to rewrite contiguous relations in a matrix form by 
putting  
\[
\mbox{\boldmath $F$}(\mbox{\boldmath $a$}) := 
{}^t({}_2F_1(\mbox{\boldmath $a$};z), \, 
{}_2F_1(\mbox{\boldmath $a$} + \mbox{\boldmath $1$};z)), \qquad 
\mbox{\boldmath $a$} := (a_1, a_2; a_3) = (\alpha, \beta; \gamma), 
\quad \mbox{\boldmath $1$} := (1, 1; 1),   
\]
and $\mbox{\boldmath $e$}_1 = (1,0;0)$, $\mbox{\boldmath $e$}_2 = 
(0,1;0)$, $\mbox{\boldmath $e$}_3 = (0,0;1)$.  
From formulas in Erd\'elyi \cite[Chapter II, \S2.8]{Erdelyi}, we 
observe that the contiguous relation raising parameter $a_i$ by one can 
be written
\begin{equation} \label{eqn:c-r}   
\mbox{\boldmath $F$}(\mbox{\boldmath $a$} + \mbox{\boldmath $e$}_i) = A_i(\mbox{\boldmath $a$}) \, 
\mbox{\boldmath $F$}(\mbox{\boldmath $a$}) \qquad (i = 1, 2, 3), 
\end{equation}
where the matrix $A_i(\mbox{\boldmath $a$})$ is given in Table \ref{tab:contiguity}, 
together with its determinant $\det A_i(\mbox{\boldmath $a$})$. 
\begin{table}[t]
\begin{align*}
A_1(\mbox{\boldmath $a$}) &:= 
\begin{pmatrix}
1 & \dfrac{\beta z}{\gamma} \\[4mm] 
-\dfrac{\gamma}{(\alpha+1)(z-1)} & 
\dfrac{\gamma-\alpha-1-\beta z}{(\alpha+1)(z-1)}
\end{pmatrix}, & 
\det A_1(\mbox{\boldmath $a$}) &= \dfrac{\gamma-\alpha-1}{(\alpha+1)(z-1)}, 
\\[4mm] 
 A_2(\mbox{\boldmath $a$}) &:= 
\begin{pmatrix}
1 & \dfrac{\alpha z}{\gamma} \\[4mm] 
-\dfrac{\gamma}{(\beta+1)(z-1)} & 
\dfrac{\gamma-\beta-1-\alpha z}{(\beta+1)(z-1)}
\end{pmatrix}, &
\det A_2(\mbox{\boldmath $a$}) &= \dfrac{\gamma-\beta-1}{(\beta+1)(z-1)}, 
\\[4mm]
A_3(\mbox{\boldmath $a$}) &:= 
\begin{pmatrix}
\dfrac{\gamma(\gamma-\alpha-\beta)}{(\gamma-\alpha)(\gamma-\beta)} & 
-\dfrac{\alpha \beta (z-1)}{(\gamma-\alpha)(\gamma-\beta)} \\[5mm] 
\dfrac{\gamma(\gamma+1)}{(\gamma-\alpha)(\gamma-\beta) z} & 
\dfrac{\gamma(\gamma+1)(z-1)}{(\gamma-\alpha)(\gamma-\beta)z}
\end{pmatrix}, &
\det A_3(\mbox{\boldmath $a$}) &= \dfrac{\gamma(\gamma+1)(z-1)}{(\gamma-\alpha)(\gamma-\beta)z}. 
\end{align*}
\caption{Contiguous matrices and their determinants.} 
\label{tab:contiguity}
\end{table}
As the compatibility conditions for three relations \eqref{eqn:c-r} 
one has the commutation relations:   
\begin{equation} \label{eqn:compatibility}
A_i(\mbox{\boldmath $a$}+\mbox{\boldmath $e$}_j) \, A_j(\mbox{\boldmath $a$}) 
= A_j(\mbox{\boldmath $a$}+\mbox{\boldmath $e$}_i) \,
 A_i(\mbox{\boldmath $a$}) \qquad 
(i,j = 1, 2, 3).    
\end{equation} 
\par
Given a lattice point $\mbox{\boldmath $p$} = (p_1, p_2; p_3) = 
(p, q; r) \in \mathbb{Z}^3_{\ge 0}$, 
a lattice path in $\mathbb{Z}^3_{\ge0}$ from $\mbox{\boldmath $0$} 
= (0,0;0)$ to $\mbox{\boldmath $p$}$ can be 
represented by a sequence $i = (i_1,\dots,i_k)$ of indexes in 
$\{1, 2, 3\}$ such that $\mbox{\boldmath $p$} = 
\mbox{\boldmath $e$}_{i_1} + \cdots + \mbox{\boldmath $e$}_{i_k}$ 
where $k = p+q+r$.  
By compatibility conditions \eqref{eqn:compatibility} the matrix product       
\[
A(\mbox{\boldmath $a$};\mbox{\boldmath $p$}) := A_{i_k}(\mbox{\boldmath $a$}+ 
\mbox{\boldmath $e$}_{i_1}+\cdots+\mbox{\boldmath $e$}_{i_{k-1}}) \cdots 
A_{i_3}(\mbox{\boldmath $a$}+\mbox{\boldmath $e$}_{i_1}+
\mbox{\boldmath $e$}_{i_2}) \, A_{i_2}(\mbox{\boldmath $a$}+
\mbox{\boldmath $e$}_{i_1}) \, A_{i_1}(\mbox{\boldmath $a$}) 
\]
is independent of the path $i = (i_1,\dots,i_k)$, that is, depends 
only on the initial point $\mbox{\boldmath $a$}$ and the terminal point 
$\mbox{\boldmath $a$}+\mbox{\boldmath $p$}$. 
The matrix version of three-term relation \eqref{eqn:contig-F} is 
expressed as 
\begin{equation} \label{eqn:m-contig-F}
\mbox{\boldmath $F$}(\mbox{\boldmath $a$}+\mbox{\boldmath $p$}) = 
A(\mbox{\boldmath $a$};\mbox{\boldmath $p$}) \, \mbox{\boldmath $F$}(\mbox{\boldmath $a$}). 
\end{equation}
\begin{lemma} \label{lem:degree} 
If $1 \le p \le r$ and $1 \le q \le r$ then 
$A(\mbox{\boldmath $a$};\mbox{\boldmath $p$})$ in 
\eqref{eqn:m-contig-F} admits a representation  
\begin{equation} \label{eqn:A(a;p)}
A(\mbox{\boldmath $a$};\mbox{\boldmath $p$}) = \dfrac{1}{(\gamma-\alpha)_{r-p} \, 
(\gamma-\beta)_{r-q}}
\begin{pmatrix}
\dfrac{(\gamma)_r \, \phi_{11}^{(r-2)}}{(\alpha+1)_{p-1} \, (\beta+1)_{q-1}} & 
\dfrac{(\gamma+1)_{r-1} \, \phi_{12}^{(r-1)}}{(\alpha+1)_{p-1} \, (\beta+1)_{q-1}} \\[6mm]
\dfrac{(\gamma)_{r+1} \, \phi_{21}^{(r-1)}}{(\alpha+1)_p \, (\beta+1)_q} & 
\dfrac{(\gamma+1)_r \, \phi_{22}^{(r)}}{(\alpha+1)_p \, (\beta+1)_q}
\end{pmatrix}, 
\end{equation}
where $\phi_{11}^{(-1)} = 0$ and $\phi_{ij}^{(k)} = 
\phi_{ij}^{(k)}(\mbox{\boldmath $a$}; \mbox{\boldmath $p$})$ is a 
polynomial of degree at most $k$ in 
$\mbox{\boldmath $a$} = (\alpha,\beta;\gamma)$ with coefficients in 
the ring $\mathbb{Z}[z^{\pm1}, 1/(z-1)]$.  
Moreover the determinant of $A(\mbox{\boldmath $a$};\mbox{\boldmath $p$})$ 
is given by  
\begin{equation} \label{eqn:detA(a;p)}
\det A(\mbox{\boldmath $a$};\mbox{\boldmath $p$}) = 
\dfrac{z^{-r}(z-1)^{r-p-q} \cdot (\gamma)_r \, (\gamma+1)_r}{(\alpha+1)_p \, 
(\beta+1)_q \, (\gamma-\alpha)_{r-p} \, (\gamma-\beta)_{r-q}}.  
\end{equation}
\end{lemma}
{\it Proof}. Formula \eqref{eqn:A(a;p)} is proved by induction on $r$,   
where the main claim is the assertion about the degrees of 
$\phi_{ij}^{(r)}$, $i,j = 1,2$, in $\mbox{\boldmath $a$} = 
(\alpha,\beta;\gamma)$. 
A direct check shows that it is true for $r = 1$, that is, for 
$\mbox{\boldmath $p$} = \mbox{\boldmath $1$}$. 
Assuming the assertion is true for $r$ we show it for $r+1$, 
that is, for $(p,q;r+1)$ with $1 \le p \le r+1$ and $1 \le q \le r+1$, 
where symmetry allows us to assume $p \le q$.   
There are three cases to deal with: (i) $p \le q \le r$; 
(ii) $p < q = r+1$; and (iii) $p = q = r+1$. 
In case (i) the relation 
$A(\mbox{\boldmath $a$};\mbox{\boldmath $p$}+\mbox{\boldmath $e$}_3) 
= A_3(\mbox{\boldmath $a$}+\mbox{\boldmath $p$}) \, 
A(\mbox{\boldmath $a$};\mbox{\boldmath $p$})$ 
with $\mbox{\boldmath $p$} = (p,q;r)$ leads to the recurrence 
\[
\begin{pmatrix}
\phi_{11}^{(r-1)} & \phi_{12}^{(r)\phantom{+1}} \\[2mm]
\phi_{21}^{(r)\phantom{-1}} & \phi_{22}^{(r+1)} 
\end{pmatrix}
= 
\begin{pmatrix}
\gamma-\alpha-\beta+r-p-q & 1-z \\[2mm]
\dfrac{(\alpha+p)(\beta+q)}{z} & \dfrac{(\gamma+r)(z-1)}{z}  
\end{pmatrix}
\begin{pmatrix}
\phi_{11}^{(r-2)} & \phi_{12}^{(r-1)} \\[2mm]
\phi_{21}^{(r-1)} & \phi_{22}^{(r)\phantom{-1}} 
\end{pmatrix},  
\]
by which the assertion for $r+1$ follows from induction hypothesis.  
In case (ii) the relation 
$A(\mbox{\boldmath $a$};\mbox{\boldmath $p$}+\mbox{\boldmath $e$}_2+
\mbox{\boldmath $e$}_3) 
= A_2(\mbox{\boldmath $a$}+\mbox{\boldmath $p$}+\mbox{\boldmath $e$}_3) \, 
A_3(\mbox{\boldmath $a$}+\mbox{\boldmath $p$}) \, 
A(\mbox{\boldmath $a$};\mbox{\boldmath $p$})$ 
with $\mbox{\boldmath $p$} = (p,r;r)$ leads to  
\[
\begin{pmatrix}
\phi_{11}^{(r-1)} & \phi_{12}^{(r)\phantom{+1}} \\[2mm]
\phi_{21}^{(r)\phantom{-1}} & \phi_{22}^{(r+1)} 
\end{pmatrix}
= 
\begin{pmatrix}
\beta + r & z-1 \\[2mm]
-\dfrac{(\alpha+p)(\beta+r)}{z} & \dfrac{\gamma+r-(\alpha+p) z}{z} 
\end{pmatrix}
\begin{pmatrix}
\phi_{11}^{(r-2)} & \phi_{12}^{(r-1)} \\[2mm]
\phi_{21}^{(r-1)} & \phi_{22}^{(r)\phantom{-1}} 
\end{pmatrix},  
\]
by which the assertion follows from induction hypothesis. 
Finally, in case (iii) the relation 
$A(\mbox{\boldmath $a$};\mbox{\boldmath $p$}+\mbox{\boldmath $1$}) = 
A_1(\mbox{\boldmath $a$}+\mbox{\boldmath $p$}+\mbox{\boldmath $e$}_2+
\mbox{\boldmath $e$}_3) \, 
A_2(\mbox{\boldmath $a$}+\mbox{\boldmath $p$}+\mbox{\boldmath $e$}_3) \, 
A_3(\mbox{\boldmath $a$}+\mbox{\boldmath $p$}) \, 
A(\mbox{\boldmath $a$};\mbox{\boldmath $p$})$ 
with $\mbox{\boldmath $p$} = (r,r;r)$ yields   
\[
\begin{pmatrix}
\phi_{11}^{(r-1)} & \phi_{12}^{(r)\phantom{+1}} \\[2mm]
\phi_{21}^{(r)\phantom{-1}} & \phi_{22}^{(r+1)} 
\end{pmatrix}
= 
\begin{pmatrix}
0 & 1 \\[4mm]
\dfrac{(\alpha+r)(\beta+r)}{z(1-z)} & 
\dfrac{\gamma+r-(\alpha+\beta+2r+1)z }{z(z-1)}  
\end{pmatrix}
\begin{pmatrix}
\phi_{11}^{(r-2)} & \phi_{12}^{(r-1)} \\[2mm]
\phi_{21}^{(r-1)} & \phi_{22}^{(r)\phantom{-1}} 
\end{pmatrix},  
\]
by which the assertion follows and the induction completes itself. 
\par
Determinant formula \eqref{eqn:detA(a;p)} is obtained by taking the 
determinant of matrix products 
\begin{align}
A(\mbox{\boldmath $a$};\mbox{\boldmath $p$}) 
&= A_1(\alpha+p-1,\beta+q;\gamma+r) \cdots A_1(\alpha+1,\beta+q;\gamma+r) \, 
A_1(\alpha,\beta+q;\gamma+r) \nonumber \\
& \qquad \cdot A_2(\alpha,\beta+q-1;\gamma+r) \cdots A_2(\alpha,\beta+1;\gamma+r) \, 
A_2(\alpha,\beta;\gamma+r)  \label{eqn:A(a;p)-2} \\
& \qquad \cdot A_3(\alpha,\beta;\gamma+r-1) \cdots A_3(\alpha,\beta;\gamma+1) \, 
A_3(\alpha,\beta;\gamma), \nonumber
\end{align} 
and by using determinant formulas in Table \ref{tab:contiguity}. 
\hfill $\Box$ \par\medskip
Lemma \ref{lem:degree} readily leads to a matrix version of 
three-term relation \eqref{eqn:contig-f}:  
\begin{equation} \label{eqn:m-contig-f}
\mbox{\boldmath $f$}(w+1) = 
A(w) \mbox{\boldmath $f$}(w), \qquad \mbox{\boldmath $f$}(w) 
:= {}^t(f(w), \tilde{f}(w)),  
\end{equation}
where the matrix $A(w)$ is described by Corollary \ref{cor:degree} below.  
Note that the $(1,1)$-entry and $(1,2)$-entry of $A(w)$ are just  
$R(w)$ and $Q(w)$ in formula \eqref{eqn:contig-f} respectively. 
\begin{corollary} \label{cor:degree} 
If $1 \le p \le r$ and $1 \le q \le r$ then $A(w)$ in \eqref{eqn:m-contig-f} 
admits a representation  
\begin{equation} \label{eqn:A(w)}
\begin{split}
A(w) &= \dfrac{1}{((r-p)w-a)_{r-p} \, ((r-q)w-b)_{r-q}} \\[2mm]
& \times 
\begin{pmatrix}
\dfrac{(r w)_r \, \phi_{11}^{(r-2)}(w)}{(p w+a+1)_{p-1} \, (q w+b+1)_{q-1}} & 
\dfrac{(r w+1)_{r-1} \, \phi_{12}^{(r-1)}(w)}{(p w+a+1)_{p-1} \, (q w+b+1)_{q-1}} \\[6mm]
\dfrac{(r w)_{r+1} \, \phi_{21}^{(r-1)}(w)}{(p w+a+1)_p \, (q w+b+1)_q} & 
\dfrac{(r w+1)_r \, \phi_{22}^{(r)}(w)}{(p w+a+1)_p \, (q w+b+1)_q}
\end{pmatrix}, 
\end{split}
\end{equation}
where $\phi_{11}^{(-1)}(w) = 0$ and $\phi_{ij}^{(k)}(w)$ is a polynomial of 
degree at most $k$ in $w$.   
Moreover,   
\begin{equation} \label{eqn:detA(w)}
\det A(w) = 
\dfrac{x^{-r}(x-1)^{r-p-q} \cdot (r w)_r \, (r w+1)_r }{(p w+a+1)_p \, 
(q w+b+1)_q \, ((r-p)w-a)_{r-p} \, ((r-q)w-b)_{r-q}}. 
\end{equation}
\end{corollary}
{\it Proof}. 
In view of \eqref{eqn:m-contig-F}, substitute $\mbox{\boldmath $a$} = 
\mbox{\boldmath $\alpha$}(w) :=(p w+a, q w+b; r w)$ and $z = x$ 
into \eqref{eqn:A(a;p)}.  \hfill $\Box$
\begin{remark} \label{rem:degree} 
In three-term relation \eqref{eqn:contig-f} we have $Q(w) = A_{12}(w)$ 
and $R(w) = A_{11}(w)$ evaluated at $z = x$, where $A_{ij}(w)$ denotes  
the $(i,j)$-th entry of the matrix $A(w)$. 
Thus a solution $\lambda = (p,q,r;a,b;x)$ to Problem \ref{prob:ocf} comes 
from contiguous relations exactly when $A_{12}(w)$ or equivalently 
$\phi_{12}^{(r-1)}(w)$ vanishes in $\mathbb{C}(w)$ (upon putting $z = x$). 
If this is the case, then taking the determinant of formula 
\eqref{eqn:A(w)} and comparing the result with formula \eqref{eqn:detA(w)}, 
we find  
\begin{equation} \label{eqn:degree} 
\phi_{11}^{(r-2)}(w) \cdot \phi_{22}^{(r)}(w) 
= x^{-r}(x-1)^{r-p-q} \cdot  
(p w+a+1)_{p-1}(q w+b+1)_{q-1}((r-p)w-a)_{r-p}((r-q)w-b)_{r-q}.  
\end{equation}
This implies $\deg \phi_{11}^{(r-2)} = r-2$ and $\deg \phi_{22}^{(r)} = r$, 
since $\deg \phi_{11}^{(r-2)} \le r-2$ and $\deg \phi_{22}^{(r)} \le r$ 
while the right-hand side of \eqref{eqn:degree} is of degree $2r-2$.  
Using formula \eqref{eqn:degree} in $R(w) = A_{11}(w)$ yields   
\begin{equation} \label{eqn:R(w)-2}
R(w) = x^{-r} (x-1)^{r-p-q} \cdot \dfrac{(r w)_r}{\phi_{22}^{(r)}(w)}. 
\end{equation}
\end{remark}
\subsection{Principal Parts of Contiguous Matrices} \label{subsec:principal}
For each $i = 1,2,3$, the matrix $A_i(\mbox{\boldmath $a$})$ with 
$\mbox{\boldmath $a$} = \mbox{\boldmath $\alpha$}(w)$ admits a limit 
$B_i := \displaystyle \lim_{w \to \infty} A_i(\mbox{\boldmath $\alpha$}(w))$, 
the ``principal part"  of $A_i(\mbox{\boldmath $\alpha$}(w))$, whose explicit 
form is given in Table \ref{tab:contiguity2}. 
Compatibility condition \eqref{eqn:compatibility} or a direct check of 
formulas in Table \ref{tab:contiguity2} implies that $B_1$, $B_2$ and $B_3$ 
are mutually commutative. 
Taking the limit as $w \to \infty$ in formula \eqref{eqn:A(w)} enables us 
to extract some information about $x$ for a solution 
$\lambda = (p,q,r;a,b;x)$ to Problem \ref{prob:ocf} in region 
\eqref{eqn:apollo}.  
\begin{table}[t]
\begin{gather*}
B_1 = 
\begin{pmatrix}
1                  & \dfrac{q z}{r} \\[4mm] 
-\dfrac{r}{p(z-1)} & \dfrac{r-p-q z}{p(z-1)}
\end{pmatrix}, 
\qquad  
B_2 = 
\begin{pmatrix}
1                  & \dfrac{p z}{r} \\[4mm] 
-\dfrac{r}{q(z-1)} & \dfrac{r-q-p z}{q(z-1)}
\end{pmatrix},  
\\[4mm] 
B_3 = \dfrac{1}{(r-p)(r-q)} 
\begin{pmatrix}
\, r(r-p-q) & - p q (z-1) \, \\[4mm] 
\dfrac{r^2}{z} & \dfrac{r^2(z-1)}{z}
\end{pmatrix}.   
\end{gather*}
\caption{Principal part $B_i$ of $A_i(\mbox{\boldmath $a$})$ with 
$\mbox{\boldmath $a$} = \mbox{\boldmath $\alpha$}(w) :=(p w+a, q w+b; r w)$.} 
\label{tab:contiguity2}
\end{table}
\begin{lemma} \label{lem:Bpqr} 
In formula \eqref{eqn:A(w)} we have  
\begin{equation} \label{eqn:Bpqr}
\lim_{w\to\infty} A(w) = B_1^p B_2^q B_3^r = c   
\begin{pmatrix}
X(z)-\{r-(p+q) z\} Y(z) & 2 (p q/r) \, z(z-1)\, Y(z) \\[2mm]
 -2 r \, Y(z)           & X(z) + \{ r-(p+q)z \} Y(z) 
\end{pmatrix}, 
\end{equation}
where $X(z)$ and $Y(z)$ are polynomials in formula \eqref{eqn:CpmAB} 
and $c$ is a constant defined by 
\[
c := \dfrac{r^r}{2^r p^p q^q (r-p)^{r-p} (r-q)^{r-q} z^r}. 
\]
\end{lemma}
{\it Proof}. 
Substitute $\mbox{\boldmath $a$} = \mbox{\boldmath $\alpha$}(w)$ into 
formula \eqref{eqn:A(a;p)-2} and take the limit $w \to \infty$ to get 
$\displaystyle \lim_{w\to\infty} A(w) = B_1^p B_2^q B_3^r$ 
by the commutativity of $B_1$, $B_2$ and $B_3$. 
From formulas in Table \ref{tab:contiguity2} we have  
\[
B_1 B_3 = 
\begin{pmatrix}
\dfrac{r}{r-q} & \dfrac{q (z-1)}{r-q} \\[4mm]
-\dfrac{r^2}{p (r-q) z} & \dfrac{r (r-q z)}{p (r-q) z}
\end{pmatrix}, 
\qquad 
B_2 B_3 = 
\begin{pmatrix}
\dfrac{r}{r-p} & \dfrac{p (z-1)}{r-p} \\[4mm]
-\dfrac{r^2}{q (r-p) z} & \dfrac{r (r-p z)}{q (r-p) z}
\end{pmatrix}.     
\]
Observe that $B_1B_3$, $B_2B_3$ and $B_3$ are simultaneously 
diagonalized as 
\begin{align*}
T^{-1} (B_1 B_3) T &= 
\dfrac{r}{2 p(r-q)z} \cdot 
\mathrm{diag}\left\{r+(p-q)z + \sqrt{\varDelta}, \,\, r+(p-q)z - \sqrt{\varDelta} \right\}, \\
T^{-1} (B_2 B_3) T &= 
\dfrac{r}{2 q(r-p)z} \cdot 
\mathrm{diag}\left\{r-(p-q)z + \sqrt{\varDelta}, \,\, r-(p-q)z - \sqrt{\varDelta} \right\}, \\
T^{-1} B_3 T &= 
\dfrac{r}{2(r-p)(r-q)z} \cdot 
\mathrm{diag}\left\{(2 r-p-q)z-r - \sqrt{\varDelta}, \,\, (2 r-p-q)z-r + \sqrt{\varDelta} \right\},    
\end{align*}
where the diagonalizing matrix $T$ is given by  
\[
T := 
\begin{pmatrix}
 \dfrac{r-(p+q)z-\sqrt{\varDelta}}{2 r} & \dfrac{r-(p+q) z+ \sqrt{\varDelta}}{2 r} \\[3mm]
 1 & 1
\end{pmatrix}.  
\]
In view of formulas \eqref{eqn:Cpm} and \eqref{eqn:CpmAB}, $B_1^p B_2^q B_3^r 
= (B_1 B_3)^p(B_2 B_3)^q B_3^{r-p-q}$ is diagonalized as   
\[
T^{-1} (B_1^p B_2^q B_3^r) T = c \cdot 
\mathrm{diag} \left\{ X(z)+Y(z)\sqrt{\varDelta}, \,\,  X(z)-Y(z) \sqrt{\varDelta} \right\}.  
\]
Then \eqref{eqn:Bpqr} follows from $B_1^p B_2^q B_3^r = c \cdot 
T \cdot \mathrm{diag} \left\{ X(z)+Y(z) \sqrt{\varDelta}, \, 
X(z)-Y(z) \sqrt{\varDelta} \right\} \cdot T^{-1}$. \hfill $\Box$
\begin{lemma} \label{lem:roots}
Polynomial $Y(w)$ has the following properties.
\begin{enumerate}
\item If $r-p-q =0$ then $Y(z)$ has no root in $0 \le z < 1$. 
\item If $r-p-q$ is a positive even integer then $Y(z)$ has 
at least one root in $0 < z < 1$.  
\end{enumerate}
\end{lemma}
{\it Proof}. 
To prove assertion (1) we assume $p+q=r$. 
Let $0 \le z < 1$ and put $s := p-q \ge 0$.  
Formulas \eqref{eqn:Delta} and \eqref{eqn:Cpm} yield  
$\varDelta = (1-z)(r^2-s^2z) > 0$ and $Z_{\pm}(z) = (r+s z \pm 
\sqrt{\varDelta})^p (r-s z \pm \sqrt{\varDelta})^q$. 
Since $(r-s z)^2 - \varDelta = 4 q^2 z \ge 0$, we have $r + s z \ge r-s z 
\ge \sqrt{\varDelta} > 0$, so $r + s z + \sqrt{\varDelta} > 
r + s z - \sqrt{\varDelta} \ge 0$ and $r - s z + \sqrt{\varDelta} > 
r - s z - \sqrt{\varDelta} \ge 0$. 
Thus formula \eqref{eqn:CpmAB} yields 
\[
2 \, Y(z) \, \sqrt{\varDelta} = (r+s z + \sqrt{\varDelta})^p(r-s z + \sqrt{\varDelta})^q
-(r+s z - \sqrt{\varDelta})^p(r-s z - \sqrt{\varDelta})^q > 0, 
\]  
which implies $Y(z) > 0$. 
Therefore $Y(z)$ has no root in $0 \le z < 1$. 
\par
To show assertion (2) we assume $r-p-q > 0$.   
Since $\sqrt{\varDelta(0)} = r > 0$, formula \eqref{eqn:Cpm} gives 
$Z_{+}(0) = (2r)^p(2r)^q (-2r)^{r-p-q} = (-1)^{r-p-q}(2r)^r$ and 
$Z_{-}(0) = 0^p \cdot 0^q \cdot 0^{r-p-q} = 0$, which are valid even if $r-p-q=0$. 
Similarly, since $\sqrt{\varDelta(1)} = r-p-q > 0$, formula \eqref{eqn:Cpm} yields  
$Z_{+}(1) = \{2(r-q)\}^p\{2(r-p)\}^q \cdot 0^{r-p-q} = 0$ and 
$Z_{-}(1) = (2p)^p(2q)^q \{2(r-p-q)\}^{r-p-q} = 2^r p^p q^q \, (r-p-q)^{r-p-q}$. 
Thus it follows from formula \eqref{eqn:CpmAB} that    
\[
Y(0) = (-1)^{r-p-q} \, (2r)^{r-1}, \qquad 
Y(1) = - 2^{r-1} p^p q^q \, (r-p-q)^{r-p-q-1} < 0.  
\]
Accordingly, if $r-p-q$ is positive and even, then we have $Y(0) > 0$ 
and $Y(1) < 0$ so that $Y(z)$ has at least one root in the interval 
$0 < z < 1$. \hfill $\Box$
\begin{proposition} \label{prop:roots}
Assertions $(1)$ and $(2)$ of Theorem $\ref{thm:contiguous3}$ hold true. 
\end{proposition}
{\it Proof}. 
Any non-elementary solution $\lambda = (p,q,r;a,b;x)$ of type $(\mathrm{A})$ 
to Problem \ref{prob:ocf} in region \eqref{eqn:apollo} comes from 
contiguous relations by Theorem \ref{thm:contiguous1}. 
Then $A_{12}(w)$ vanishes identically by Remark \ref{rem:degree}. 
Thus formula \eqref{eqn:Bpqr} evaluated at $z = x$ yields 
$0 = \displaystyle \lim_{w \to \infty} A_{12}(w) = 2 (p q/r) x(x-1) Y(x)$, namely,  
$Y(x) =0$, which shows assertion (1) of Theorem \ref{thm:contiguous3}. 
Note that $r-p-q \ge 0$ is a priori assumed since $\lambda$ is 
in region \eqref{eqn:apollo}. 
The possibility of $r-p-q=0$ is ruled out because in that case 
$Y(z)$ would have no root in $0 \le z < 1$ by assertion (1) of 
Lemma \ref{lem:roots}, while $r-p-q$ must be even by Theorem \ref{thm:integer}. 
Therefore assertion (2) of Theorem \ref{thm:contiguous3} follows. 
\hfill $\Box$
\subsection{Truncated Hypergeometric Products} \label{subsec:ebisu}
In our vectorial notation the three-term relation \eqref{eqn:contig-F} 
can be written 
\begin{equation} \label{eqn:three-term}
{}_2F_1(\mbox{\boldmath $a$}+\mbox{\boldmath $p$};z) = 
r(\mbox{\boldmath $a$};\mbox{\boldmath $p$};z) \, 
{}_2F_1(\mbox{\boldmath $a$};z) + 
q(\mbox{\boldmath $a$};\mbox{\boldmath $p$};z) \, 
{}_2F_1(\mbox{\boldmath $a$}+\mbox{\boldmath $1$};z),  
\end{equation}
where dependence upon $\mbox{\boldmath $p$} = (p,q;r)$ is also emphasized. 
Replacing $\mbox{\boldmath $p$}$ by $\mbox{\boldmath $p$}+\mbox{\boldmath $1$}$ in 
\eqref{eqn:three-term} gives another relation of a similar sort.   
Comparing these two relations with formula \eqref{eqn:m-contig-F}, 
we find 
\[
A(\mbox{\boldmath $a$};\mbox{\boldmath $p$}) = 
\begin{pmatrix}
r(\mbox{\boldmath $a$};\mbox{\boldmath $p$};z) & 
q(\mbox{\boldmath $a$};\mbox{\boldmath $p$};z) \\
r(\mbox{\boldmath $a$};\mbox{\boldmath $p$}+\mbox{\boldmath $1$};z) & 
q(\mbox{\boldmath $a$};\mbox{\boldmath $p$}+\mbox{\boldmath $1$};z) 
\end{pmatrix}. 
\]
Ebisu \cite{Ebisu1} made an extensive study of relation 
\eqref{eqn:three-term} and his result allows us to represent each 
entry of the matrix $A(\mbox{\boldmath $a$};\mbox{\boldmath $p$})$ in 
terms of a truncated hypergeometric product; we are especially 
interested in the second column of $A(\mbox{\boldmath $a$};\mbox{\boldmath $p$})$, 
that is, in $q(\mbox{\boldmath $a$};\mbox{\boldmath $p$};z)$ and 
$q(\mbox{\boldmath $a$};\mbox{\boldmath $p$}+\mbox{\boldmath $1$};z)$.     
Recall that $\langle \varphi(z) \rangle_k = \sum_{j=0}^k c_k z^j$ 
stands for the truncation at degree $k$ of a power series 
$\varphi(z) = \sum_{j=0}^{\infty} c_k z^j$.  
\begin{lemma} \label{lem:ebisu}
Let $\mbox{\boldmath $p$} = (p,q;r) \in \mathbb{Z}^3$. 
If $0 \le q \le p$ and $p+q \le r$ then
\begin{equation} \label{eqn:ebisu-q1} 
q(\mbox{\boldmath $a$};\mbox{\boldmath $p$};z) = z^{1-r}(z-1) 
\cdot C(\mbox{\boldmath $a$};\mbox{\boldmath $p$}) \cdot
\langle {}_2F_1(\mbox{\boldmath $a$}^*;z) \cdot 
{}_2F_1(\mbox{\boldmath $v$}-\mbox{\boldmath $a$}^*-\mbox{\boldmath $p$}^*;z) 
\rangle_{r-q-1},   
\end{equation}
while if $-1 \le q \le p$ and $p+q \le r-1$ then
\begin{equation} \label{eqn:ebisu-q2}
q(\mbox{\boldmath $a$};\mbox{\boldmath $p$}+\mbox{\boldmath $1$};z) = 
z^{-r}(z-1) \cdot C(\mbox{\boldmath $a$};\mbox{\boldmath $p$}+\mbox{\boldmath $1$}) 
\cdot \langle {}_2F_1(\mbox{\boldmath $a$}^*;z) \cdot 
{}_2F_1(\mbox{\boldmath $1$}-\mbox{\boldmath $a$}^*-\mbox{\boldmath $p$}^*;z) 
\rangle_{r-q-1},    
\end{equation}
where $\mbox{\boldmath $a$}^* := (\gamma-\alpha, \gamma-\beta ; \gamma)$, 
$\mbox{\boldmath $p$}^* := (r-p,r-q;r)$, $\mbox{\boldmath $v$} = (1,1;2)$ 
and    
\[ 
C(\mbox{\boldmath $a$};\mbox{\boldmath $p$}) := (-1)^{r-p-q} \, 
\dfrac{(\gamma)_{r-1} \, (\gamma+1)_{r-1}}{(\alpha+1)_{p-1}(\beta+1)_{q-1}
(\gamma-\alpha)_{r-p}(\gamma-\beta)_{r-q}}.  
\]
\end{lemma}
{\it Proof}.  
By case (ii) of 
Ebisu \cite[Proposition 3.4, Theorems 3.7, Remark 3.11]{Ebisu1}, 
if $p \ge q$ and $r \ge \max\{p+q, 0\}$ then 
$q(\mbox{\boldmath $a$};\mbox{\boldmath $p$};z) = 
z^{1-r} (z-1) \, q_0(\mbox{\boldmath $a$};\mbox{\boldmath $p$};z)$, 
where $q_0(\mbox{\boldmath $a$};\mbox{\boldmath $p$};z)$ is a polynomial 
of degree at most $r-q-1$ in $z$ and it is explicitly given by    
\[
q_0(\mbox{\boldmath $a$};\mbox{\boldmath $p$};z) = 
C(\mbox{\boldmath $a$};\mbox{\boldmath $p$}) \,
 {}_2F_1(\mbox{\boldmath $a$}^*; z) \, 
{}_2F_1 (\mbox{\boldmath $v$}-\mbox{\boldmath $a$}^*-\mbox{\boldmath $p$}^*; z) - 
\dfrac{\alpha \beta}{\gamma(\gamma-1)} \, 
z^r \, {}_2F_1(\mbox{\boldmath $v$}-\mbox{\boldmath $a$}; z) \, 
{}_2F_1(\mbox{\boldmath $a$}+\mbox{\boldmath $p$}; z). 
\]
In particular, if $p \ge q \ge 0$ and $p+q \le r$ then we have 
$r-q-1 < r$ so that 
\begin{equation} \label{eqn:ebisu-q3}
q_0(\mbox{\boldmath $a$};\mbox{\boldmath $p$};z) = 
C(\mbox{\boldmath $a$};\mbox{\boldmath $p$}) \, 
\langle  {}_2F_1(\mbox{\boldmath $a$}^*; z) \cdot 
{}_2F_1 (\mbox{\boldmath $v$}-\mbox{\boldmath $a$}^*-\mbox{\boldmath $p$}^*; z) 
\rangle_{r-q-1},   
\end{equation}
which yields formula \eqref{eqn:ebisu-q1}. 
Formula \eqref{eqn:ebisu-q2} follows from formula \eqref{eqn:ebisu-q3} 
with $\mbox{\boldmath $p$}$ replaced by $\mbox{\boldmath $p$} + 
\mbox{\boldmath $1$}$, which holds true provided 
$p+1 \ge q+1 \ge 0$ and $(p+1)+(q+1) \le r+1$, that is, $p \ge q \ge -1$ 
and $p+q \le r-1$. 
Here we used $\mbox{\boldmath $v$}-\mbox{\boldmath $a$}^*-
(\mbox{\boldmath $p$}+\mbox{\boldmath $1$})^* = 
\mbox{\boldmath $1$}-\mbox{\boldmath $a$}^*-\mbox{\boldmath $p$}^*$. 
\hfill $\Box$ 
\begin{lemma} \label{lem:ebisu2}
If $\mbox{\boldmath $p$} = (p,q;r) \in \mathbb{Z}^3$ satisfies 
$0 \le q \le p$ and $p+q \le r$,  then evaluated at $z =x$, 
\begin{equation} \label{eqn:phi12} 
\phi_{12}^{(r-1)}(w) = (-1)^{r-p-q} \cdot x^{1-r}(x-1) 
\cdot \Phi(w;\lambda), 
\end{equation}
where $\Phi(w;\lambda)$ is defined by formula \eqref{eqn:truncate}. 
In particular, $\Phi(w,\lambda)$ is a polynomial of degree at most 
$r-1$ in $w$ and condition $Y(x) = 0$ is equivalent to equation 
$\Phi_{r-1}(\lambda) = 0$ in system \eqref{eqn:Phi}.    
If moreover $\mbox{\boldmath $p$}$ satisfies $p + q \le r-1$, 
then evaluated at $z =x$, 
\begin{equation} \label{eqn:phi22} 
\phi_{22}^{(r)}(w) = (-1)^{r-p-q-1} \cdot x^{-r} (x-1) \cdot P(w),  
\end{equation}
where $P(w)$ is defined by formula \eqref{eqn:P}. 
In particular $P(w)$ is a polynomial of degree at most $r$.    
\end{lemma}
{\it Proof}. 
Substituting $\mbox{\boldmath $a$} = \mbox{\boldmath $\alpha$}(w) := 
(p w+a, q w+b; r w)$ and $z = x$ 
into formula \eqref{eqn:ebisu-q1} we find 
\begin{equation} \label{eqn:A12}
A_{12}(w) = 
\dfrac{(-1)^{r-p-q} \cdot x^{1-r}(x-1) \cdot (r w+1)_{r-1} \cdot 
\Phi(w;\lambda)}{(p w+a+1)_{p-1}(q w+b+1)_{q-1}((r-p)w-a)_{r-p} 
((r-q)w-b)_{r-q}}.   
\end{equation}
Comparing this with formula \eqref{eqn:A(w)} yields 
formula \eqref{eqn:phi12}, which shows that $\Phi(w;\lambda)$ is 
of degree at most $r-1$ in $w$ since so is $\phi_{12}^{(r-1)}(w)$ 
by Corollary \ref{cor:degree}.  
Moreover formula \eqref{eqn:A12} gives
\[
\lim_{w\to\infty} A_{12}(w) = 
\dfrac{(-1)^{r-p-q} \cdot x^{1-r}(x-1) \cdot r^{r-1} \cdot 
\Phi_{r-1}(\lambda)}{p^{p-1} q^{q-1} (r-p)^{r-p} (r-q)^{r-q}}.   
\]
This together with formula \eqref{eqn:Bpqr} at $z = x$ yields   
$Y(x) = (-1)^{r-p-q} 2^{r-1}\Phi_{r-1}(\lambda)$. 
Thus $Y(x) = 0$ is equivalent to $\Phi_{r-1}(\lambda) = 0$.  
Finally, if $P(w)$ is defined by \eqref{eqn:P} then 
formula \eqref{eqn:ebisu-q2} is compared with  \eqref{eqn:A(w)} to 
yield formula \eqref{eqn:phi22}, from which the assertion for 
$\deg P(w)$ also follows. \hfill $\Box$ 
\begin{proposition} \label{prop:ebisu} 
Theorem $\ref{thm:contiguous2}$ and assertion $(3)$ of Theorem 
$\ref{thm:contiguous3}$ hold true. 
\end{proposition}
{\it Proof}. 
By remark \ref{rem:degree}, a solution $\lambda$ comes from 
contiguous relations if and only if $\phi_{12}^{(r-1)}(w)$ evaluated 
at $z = x$ vanishes in $\mathbb{C}(w)$.  
By formula \eqref{eqn:phi12}, this is equivalent to saying 
that $\Phi(w;\lambda)$ vanishes in $\mathbb{C}(w)$, from which 
Theorem \ref{thm:contiguous2} follows. 
When $\lambda$ comes from contiguous relations, the degree of $P(w)$ 
is exactly $r$ since so is $\phi_{22}^{(r)}(w)$ by Remark \ref{rem:degree}. 
Formula \eqref{eqn:R(w)} in Theorem \ref{thm:contiguous3} is then  
obtained by substituting formula \eqref{eqn:phi22} into \eqref{eqn:R(w)-2}. 
Division relation \eqref{eqn:division} follows easily from formulas 
\eqref{eqn:degree} and \eqref{eqn:phi22}. 
Thus assertion (3) of Theorem \ref{thm:contiguous3} is established. 
\hfill $\Box$
\subsection{Terminating Hypergeometric Sums} \label{subsec:terminate}
The condition \eqref{eqn:Phi} leads to an algebraic system 
involving {\sl terminating} hypergeometric sums.  
To see this we employ a {\sl renormalized} terminating 
hypergeometric sum:  
\begin{equation} \label{eqn:thg}
\mathcal{F}_k(\beta;\gamma;z) := (\gamma)_k \cdot {}_2F_1(-k,\beta;\gamma;z) = 
\sum_{j=0}^k (-1)^j \binom{k}{j} (\beta)_j (\gamma+j)_{k-j} \, z^j 
\quad (k \in \mathbb{N} \cup \{0\}).   
\end{equation}
Note that $\mathcal{F}_k(\beta;\gamma;z)$ is a polynomial 
of $(\beta;\gamma;z)$. 
By evaluating $\Phi(w;\lambda)$ in definition \eqref{eqn:truncate} at 
\begin{equation} \label{eqn:wkwj}
w_k^* := \dfrac{a-k}{r-p} \quad (0 \le k \le r-p-1); 
\qquad 
w_j := -\dfrac{a+j}{p} \quad (0 \le j \le p-1),  
\end{equation}
where $k$ and $j$ are integers in the indicated intervals, 
Theorem \ref{thm:contiguous2} yields the following.  
\begin{proposition} \label{prop:terminate1} 
System \eqref{eqn:Phi} leads to a total of $r$ algebraic 
equations for $(a,b;x):$   
\begin{subequations} \label{eqn:terminate1}
\begin{align}
(\gamma_k^*+k)_p \cdot \mathcal{F}_k(\beta_k^*; \gamma_k^*; x) \cdot 
\mathcal{F}_{r-p-1-k}(\tilde{\beta}_k^*; \tilde{\gamma}_k^*;x) &=0 
&& (0 \le k \le r-p-1), \label{eqn:terminate1-1} \\
(\gamma_j+j)_{r-p} \cdot \mathcal{F}_j(\beta_j; \gamma_j; x) \cdot 
\mathcal{F}_{p-1-j}(\tilde{\beta}_j+1; \tilde{\gamma}_j+1;x) &=0 
&& (0 \le j \le p-1), \label{eqn:terminate1-2}
\end{align}
\end{subequations}
each of which consists of a factorial and two 
terminating hypergeometric factors, where 
\begin{align*}
\beta_k^* &:= (r-q) w_k^*-b, & 
\gamma_k^* &:= r w_k^*, & 
\tilde{\beta}_k^* &:= 1-(r-q)(w_k^*+1)+b, & 
\tilde{\gamma}_k^* &:= 2-r(w_k^*+1) \\ 
\beta_j &:= q w_j+b, & 
\gamma_j &:= r w_j, & 
\tilde{\beta}_j &:= -q(w_j+1)-b, &   
\tilde{\gamma}_j &:= 1-r(w_j+1). 
\end{align*}
Moreover system \eqref{eqn:terminate1} leads back to and 
hence is equivalent to the original system \eqref{eqn:Phi}, if 
\begin{equation} \label{eqn:distinct}
a \neq \frac{p_1}{r_p}(k+j) - j \quad 
(0 \le {}^{\forall} k \le r-p-1, \,\, 0 \le {}^{\forall} j \le p-1),  
\end{equation}
in particular, if $r_p a \not\in \mathbb{Z}$, where $p_1/r_p$ is the 
reduced expression of $p/r \in \mathbb{Q}$. 
\end{proposition}
{\it Proof}. 
If we substitute $w = w_k^*$ in the first formula of definition  
\eqref{eqn:truncate}, then the two hypergeometric series inside the 
bracket $\langle \cdots \rangle_{r-q-1}$ terminate at degrees $k$ 
and $r-p-1-k$ in $z$, so their product is of degree at most 
$k + (r-p-1-k) = r-p-1 \le r-q-1$ in $z$.  
Thus $\Phi(w;\lambda)$ can be evaluated at $w = w_k^*$ without 
taking truncation.  
A bit of calculation shows    
\[
\Phi(w_k^*;\lambda) = (-1)^{r-p-1-k} \cdot (\gamma_k^*+k)_p \cdot  
\mathcal{F}_k(\beta_k^*; \gamma_k^*; x) \cdot 
\mathcal{F}_{r-p-1-k}(\tilde{\beta}_k^*; \tilde{\gamma}_k^*;x).  
\]
Thus the vanishing of $\Phi(w;\lambda)$ stated in  
Theorem \ref{thm:contiguous2} yields the $r-p$ equations in 
\eqref{eqn:terminate1-1}. 
Similarly, if we substitute $w = w_j$ in the second formula of 
\eqref{eqn:truncate}, then the two hypergeometric series 
inside the bracket $\langle \cdots \rangle_{r-q-1}$ terminate at 
degrees $k$ and $p-1-k$ in $z$, so $(1-z)^{r-p-q}$ times their 
product is of degree at most $(r-p-q) + k + (p-1-k) = r-q-1$ in $z$.  
Thus $\Phi(w;\lambda)$ can also be evaluated at $w = w_j$ without 
taking truncation.  
After some calculations,    
\[
\Phi(w_j;\lambda) = (-1)^{p-1-k} \cdot (1-x)^{r-p-q} \cdot 
(\gamma_j+j)_{r-p} \cdot  
\mathcal{F}_j(\beta_j; \gamma_j; x) \cdot 
\mathcal{F}_{p-1-j}(\tilde{\beta}_j+1; \tilde{\gamma}_j+1;x),   
\]
which together with the vanishing of $\Phi(w_j;\lambda)$ leads to 
the $p$ equations in formula \eqref{eqn:terminate1-2}. 
Note that \eqref{eqn:distinct} is the condition that any pair of 
$w_k^*$ and $w_j$ in \eqref{eqn:wkwj} be distinct. 
If this is the case then equations \eqref{eqn:terminate1} imply that 
$\Phi(w;\lambda)$, which is a polynomial of degree at most $r-1$ 
in $w$, vanishes at distinct $r$ points and hence vanishes identically, 
leading to equations \eqref{eqn:Phi}. \hfill $\Box$ \par\medskip 
As in the proof of Proposition \ref{prop:terminate1}, $P(w)$ in 
definition \eqref{eqn:P} can be evaluated as    
\begin{subequations} \label{eqn:Pwkwj}
\begin{align}
P(w_k^*) &= (-1)^{r-p-1-k} \cdot (\gamma_k^*+k)_{p+1} \cdot    
\mathcal{F}_k(\beta_k^*; \gamma_k^*; x) \cdot 
\mathcal{F}_{r-p-1-k}(\tilde{\beta}_k^*; \tilde{\gamma}_k^*-1;x), 
\label{eqn:Pwk} \\
P(w_j) &= (-1)^{p-k} \cdot (1-x)^{r-p-q-1} \cdot 
(\gamma_j+j)_{r-p} \cdot \mathcal{F}_j(\beta_j;\gamma_j;x) \cdot 
\mathcal{F}_{p-j}(\tilde{\beta}_j; \tilde{\gamma}_j;x). 
\label{eqn:Pwj} 
\end{align} 
\end{subequations}
In item (2) of Remark \ref{rem:contiguous} we posed a question about 
the factors of $P(w)$. 
It can be discussed by comparing formulas \eqref{eqn:Pwkwj} with 
equations \eqref{eqn:terminate1} and by using the following.      
\begin{lemma} \label{lem:terminate} 
Let $k \in \mathbb{N}$, $\beta$, $\gamma \in \mathbb{C}$ and 
$x \in \mathbb{C} \setminus \{0,1\}$ be fixed, while $z$ be a symbolic variable.  
\begin{enumerate}
\item $\mathcal{F}_k(\beta;\gamma;z) = 0$ in $\mathbb{C}[z]$ if and only if 
$\beta, \gamma \in \mathbb{Z}$ and $0 \le -\beta \le -\gamma \le k-1$.
\item If $\mathcal{F}_k(\beta;\gamma;x) = \mathcal{F}_{k-1}(\beta+1;\gamma+1;x) = 0$ 
then $\mathcal{F}_k(\beta;\gamma;z) = 0$ in $\mathbb{C}[z]$. 
\item If $\mathcal{F}_k(\beta;\gamma;x) = \mathcal{F}_k(\beta;\gamma-1;x) = 0$ 
then $\mathcal{F}_k(\beta;\gamma;z) = 0$ in $\mathbb{C}[z]$.  
\end{enumerate}
\end{lemma}
{\it Proof}. 
By definition \eqref{eqn:thg}, $\mathcal{F}_k(\beta;\gamma;z) = 0$ in 
$\mathbb{C}[z]$ if and only if $(\beta)_j (\gamma+j)_{k-j} = 0$ for every 
$j = 0,\dots,k$.  
Putting $j=0$ there implies $\gamma = -j_0$ with 
some $j_0 \in \{0,\dots,k-1\}$.  
Putting $j = j_0+1$ then implies $\beta = -i_0$ with some 
$i_0 \in \{0, \dots, j_0\}$. 
These are sufficient to guarantee the condition 
$(\beta)_j (\gamma+j)_{k-j} = 0$ for every $j = 0, \dots,k$, 
and hence assertion (1) follows. 
\par
It follows from Andrews et al. \cite[formulas (2.5.1) and (2.5.7)]{AAR} 
and definition \eqref{eqn:thg} that 
\begin{subequations} \label{eqn:t-contig}
\begin{align}
\textstyle \frac{d}{dz} \mathcal{F}_k(\beta;\gamma;z) 
&= -k \beta \, \mathcal{F}_{k-1}(\beta+1;\gamma+1;z),  \label{eqn:t-contig1} \\
z \textstyle \frac{d}{dz} \mathcal{F}_k(\beta;\gamma;z) 
&= (\gamma+k-1) \, \mathcal{F}_k(\beta;\gamma-1;z) 
-(\gamma-1) \, \mathcal{F}_k(\beta;\gamma;z). \label{eqn:t-contig2} 
\end{align}
\end{subequations}
Assumption of assertion (2) and formula \eqref{eqn:t-contig1} 
yield a vanishing initial condition $\mathcal{F}_k(\beta;\gamma;z) = 
\frac{d}{dz}\mathcal{F}_k(\beta;\gamma;z) = 0$ at $z = x$. 
As a solution to a Gauss hypergeometric equation, which is regular 
at $z = x$ $(\neq 0, 1)$, the polynomial $\mathcal{F}_k(\beta;\gamma;z)$ 
vanishes identically in $\mathbb{C}[z]$. 
Thus assertion (2) is established.  
Assertion (3) is proved in a similar manner by using 
formula \eqref{eqn:t-contig2}.  
\hfill $\Box$ \par\medskip
Assertion (1) of Lemma \ref{lem:terminate} leads us to think of 
the following conditions:
\begin{subequations} \label{eqn:vanish}
\begin{align}
\tilde{\beta}_k^*, \,\, \tilde{\gamma}_k^* \in \mathbb{Z}, \qquad & 
0 \le -\tilde{\beta}_k^* \le -\tilde{\gamma}_k^* \le r-p-k-2, 
\label{eqn:vanish1} \\
\tilde{\beta}_j, \,\, \tilde{\gamma}_j \in \mathbb{Z}, \qquad & 
0 \le -\tilde{\beta}_j \le -\tilde{\gamma}_j \le p-j-1. 
\label{eqn:vanish2} 
\end{align} 
\end{subequations}
Each of them is an extremely restrictive condition which in particular 
implies $r_p a \in \mathbb{Z}$ and $r_q b \in \mathbb{Z}$, 
where $p_1/r_p$ and $q_1/r_q$ are the reduced expressions of 
$p/r$ and $q/r$, respectively. 
\begin{proposition} \label{prop:terminate2}
As to the question in item $(2)$ of Remark $\ref{rem:contiguous}$,  
\begin{enumerate}
\item $(w-w_k^*)|P(w)$ if and only if $(\gamma_k^*+k)_{p+1} 
\cdot \mathcal{F}_k(\beta_k^*; \gamma_k^*; x) = 0$, unless \eqref{eqn:vanish1} 
is satisfied;  
\item $(w-w_j)|P(w)$ if and only if $(\gamma_j+j)_{r-p}  
\cdot \mathcal{F}_j(\beta_j; \gamma_j; x) = 0$, unless \eqref{eqn:vanish2} 
is satisfied,  
\end{enumerate}
where the ``unless" phrase is not needed for stating the ``if" part.  
\end{proposition}
{\it Proof}. 
The ``if" part of assertion (1) follows immediately from formula 
\eqref{eqn:Pwk}. 
We shall show the ``only if" part of assertion (1).  
Assume that $(w-w_k^*)|P(w)$, that is, $P(w_k^*) = 0$, but 
$(\gamma_k^*+k)_{p+1} \cdot \mathcal{F}_k(\beta_k^*; \gamma_k^*; x) \neq 0$. 
Formulas \eqref{eqn:terminate1-1} and \eqref{eqn:Pwk} then imply 
$\mathcal{F}_{r-p-1-k}(\tilde{\beta}_k^*; \tilde{\gamma}_k^*;x) = 0$ and 
$\mathcal{F}_{r-p-1-k}(\tilde{\beta}_k^*; \tilde{\gamma}_k^*-1;x) = 0$. 
By assertion (3) of Lemma \ref{lem:terminate}, 
$\mathcal{F}_{r-p-1-k}(\tilde{\beta}_k^*; \tilde{\gamma}_k^*;z)$ must 
vanish identically in $\mathbb{C}[z]$. 
Assertion (1) of the same lemma then leads to 
condition \eqref{eqn:vanish1}. 
Assertion (2) can be proved in a similar manner by using formulas 
\eqref{eqn:terminate1-2} and \eqref{eqn:Pwj}. 
\hfill $\Box$ \par\medskip
In Propositions \ref{prop:terminate1} and \ref{prop:terminate2} 
one can replace $(p,a)$ by $(q,b)$, since definitions 
\eqref{eqn:truncate} and \eqref{eqn:P} are symmetric with 
respect to $(p,a)$ and $(q,b)$. 
Indeed, a priori the truncation there should be  
$\langle \cdots \rangle_{\max\{r-p-1, \, r-q-1\}}$, but it 
becomes $\langle \cdots \rangle_{r-q-1}$ since we are 
working in region \eqref{eqn:apollo}. 
The $(q,b)$-version of these propositions should equally be 
taken into account in our consideration. 
\par
As is mentioned in item (2) of Remark \ref{rem:contiguous}, each 
term in formula \eqref{eqn:factor} appears as a factor of $P(w)$. 
The reason for this statement is as follows:   
If we put $k = r-p-1$ in equation \eqref{eqn:terminate1-1}, 
then we get $(\gamma_k^*+k)_p \cdot \mathcal{F}_k(\beta_k^*;\gamma_k^*;x) 
= 0$, since $\mathcal{F}_{r-p-1-k}(\tilde{\beta}_k^*; \tilde{\gamma}_k^*;x) 
= \mathcal{F}_0(\tilde{\beta}_k^*; \tilde{\gamma}_k^*;x) = 1$, so the  
first term in \eqref{eqn:factor} must be a factor of $P(w)$ by the 
``if" part of assertion (1) of Proposition \ref{prop:terminate2}. 
To deal with the third term in formula \eqref{eqn:factor}, put 
$j = p-1$ in equation \eqref{eqn:terminate1-1} (if $p \ge 2$) 
and use assertion (2) of Proposition \ref{prop:terminate2}. 
As for the second and fourth terms in \eqref{eqn:factor}, 
proceed in a similar manner with the $(q,b)$-versions of 
Propositions \ref{prop:terminate1} and \ref{prop:terminate2}.    
\section{Concluding Discussions} \label{sec:closing}
We conclude this article by providing a further result and discussing 
some future directions. 
\begin{table}[t]
\begin{center}
\begin{tabular}{c|ccc}
\hline
     &      &      &      \\[-4mm]
case & (C1) & (C2) & (C3) \\[1mm]
\hline
 &   &   &     \\[-4mm]
I& F & F & --- \\[1mm]
II & F & T & --- \\[1mm]
III & T & F & --- \\[1mm]  
IV & T & T & F  \\[1mm]
V & T & T & T  \\[1mm]
\hline  
\end{tabular}
\end{center}
\caption{Division into five cases.}
\label{tab:five-cases}
\end{table}
Working in region \eqref{eqn:square}, we are interested in the (equal) 
number $m = n$ of gamma factors in the numerator or denominator of GPF 
\eqref{eqn:gpf} when it is written in a canonical form.  
Arithmetically, the difference $N := r-m$, which is referred to as 
the {\sl deficiency}, is more meaningful than $m$ itself. 
In region \eqref{eqn:square} and hence under condition (A) or (B) in 
Theorem \ref{thm:integer}, we set: \\[2mm]
\quad $\bullet$ $p/r = p_1/r_p$ : the reduced expression, that is, 
$\gcd\{p_1, \, r_p\} = 1$; \\[2mm]  
\quad $\bullet$ $q/r = q_1/r_q$ : the reduced expression, that is, 
$\gcd\{q_1, \, r_q\} = 1$; \\[2mm] 
\quad $\bullet$ $a' := p' r_p a$ $(\in \mathbb{R})$, where $p'$ is an integer 
such that $p' p_1 \equiv 1 \mod r_p$; \\[2mm]  
\quad $\bullet$ $b' := q' r_q b$ $(\in \mathbb{R})$, where $q'$ is an integer 
such that $q' q_1 \equiv 1 \mod r_q$. \\[2mm]
With this notation we introduce the following three conditions: 
\begin{center}
(C1) \, $r_p \, a \in \mathbb{Z}$; \qquad   
(C2) \, $r_q \, b \in \mathbb{Z}$; \qquad 
(C3) \, $a' \equiv b' \mod r_{pq} := \gcd\{ r_p, \, r_q \}$,    
\end{center}
where condition (C3) is well defined, that is, independent of the choice 
of $p'$ and $q'$. 
We divide non-elementary solutions into five cases as in 
Table \ref{tab:five-cases} according to whether these conditions are true 
T or false F, where (C3) makes sense only when both (C1) and (C2) are true. 
Then an amplification of the density argument in \S\ref{sec:p-r} yields 
the following result.       
\begin{result} \label{res:deficiency} 
Let $\lambda = (p,q,r;a,b;x)$ be a non-elementary solution in region 
\eqref{eqn:square}.  
\begin{enumerate}
\item If $\lambda$ is a solution of type $(\mathrm{A})$, 
then the deficiency $N$ is given as in Table $\ref{tab:deficiency1}$.  
\item If $\lambda$ is a solution of type $(\mathrm{B})$, 
then the deficiency $N$ is given as in Table $\ref{tab:deficiency2}$, 
where cases $\mathrm{II}$, $\mathrm{III}$, $\mathrm{V}$ cannot occur.  
In case $\mathrm{IV}$ we must have $\gcd\{2p, \, r\} = \gcd\{2q, \, r\}$ 
with this equal number giving the deficiency $N$ and upon putting 
$\rho := r/N$ $(\in \mathbb{N})$ we must also have 
\begin{center}
$(\mathrm{i})$ \quad $a' \not\equiv b' \mod 2\rho$, \hspace{20mm} 
$(\mathrm{ii})$ \quad $a' \equiv b' \mod \rho$.    
\end{center}
\end{enumerate}
\end{result}
\par
The proof of this result is not given here 
to keep this article in a moderate length. 
We remark that (i) is equivalent to the defining condition for case $\mathrm{IV}$ 
that (C3) should be false, while (ii) is a further necessary condition for 
this case to occur. 
Note that all solutions in Tables \ref{tab:solutions} and \ref{tab:solutions2} 
are in case $\mathrm{IV}$. 
So far we have known no solutions of any other cases. 
In particular we do not know if there is any solution with {\sl null} 
deficiency $N = 0$, that is, with gamma factors in full, in region 
\eqref{eqn:square}. 
\begin{table}[t]
\begin{minipage}{0.5\hsize}
\begin{center}
\begin{tabular}{c|l}
\hline
     &                \\[-4mm]
case & deficiency $N$ \\[1mm]
\hline
  &     \\[-4mm]
I & $0$ \\[1mm]
II & $\gcd\{q, \, r\}$ \\[1mm]
III & $\gcd\{p, \, r\}$ \\[1mm]
IV & $\gcd\{p, \, r\} + \gcd\{q, \, r\}$ \\[1mm]
V & $\gcd\{p, \, r\} + \gcd\{q, \, r\} - \gcd\{p, \, q, \, r\}$ \\[1mm]
\hline  
\end{tabular}
\end{center}
\caption{The case of type (A).}
\label{tab:deficiency1}
\end{minipage}
\begin{minipage}{0.45\hsize}
\begin{center}
\begin{tabular}{c|l}
\hline
     &                \\[-4mm]
case & deficiency $N$ \\[1mm]
\hline
  &     \\[-4mm]
I & $0$ \\[1mm]
II & cannot occur \\[1mm]
III & cannot occur \\[1mm]
IV & $\gcd\{2p, \, r\} = \gcd\{2q, \, r\}$ \\[1mm]
V & cannot occur \\[1mm]
\hline  
\end{tabular}
\end{center}
\caption{The case of type (B).}
\label{tab:deficiency2}
\end{minipage}
\end{table}
Elsewhere, however, such solutions certainly exists. 
\par
Indeed, for any positive integers $j$ and $k$ with $j > k$, if we put 
$p = -q = j-k$, $r = j+k$, $a = c$, $b = 1-c$ and $x = 1/2$, where $c$ 
is a free parameter, then there exists a GPF:   
\begin{equation} \label{eqn:bailey}
\begin{split}
& {}_2F_1(\, (j-k)w+c, \, -(j-k)w+1-c; \, (j+k)w; \, 1/2) \\[2mm] 
& \qquad = \dfrac{\sqrt{2} \, k^{c/2}}{j^{(c-1)/2} \, (j+k)^{1/2}} 
\left\{\dfrac{(j+k)^{j+k}}{2^{j+k} \, j^j \, k^k} \right\}^w 
\dfrac{\prod_{\nu=0}^{j+k-1} \varGamma\big(w+\frac{\nu}{j+k}\big)}{\prod_{\nu=0}^{j-1} 
\varGamma\big(w+\frac{c}{2j}+\frac{\nu}{j}\big) \prod_{\nu=0}^{k-1} 
\varGamma\big(w+\frac{1-c}{2k}+\frac{\nu}{k}\big)},     
\end{split}
\end{equation}
which can be derived from Bailey's formula with two free parameters 
\cite[\S2.4, formula (3)]{Bailey}: 
\begin{equation} \label{eqn:bailey2}
{}_2F_1(\alpha, 1-\alpha; \beta; 1/2) = 
\dfrac{\varGamma\big(\frac{\beta}{2}\big)
\varGamma\big(\frac{\beta+1}{2}\big)}{\varGamma\big(\frac{\beta+\alpha}{2}\big)
\varGamma\big(\frac{\beta-\alpha+1}{2}\big)}
= 
\dfrac{\sqrt{\pi} \, 2^{1-\beta} \, 
\varGamma(\beta)}{\varGamma\big(\frac{\beta+\alpha}{2}\big)
\varGamma\big(\frac{\beta-\alpha+1}{2}\big)} \qquad (\alpha, \beta \in \mathbb{C}), 
\end{equation} 
by putting $\alpha = (j-k)w+c$ and $\beta = (j+k)w$ and using Gauss's multiplication 
formula for the gamma function \cite[Theorem 1.5.2]{AAR}. 
Not lying in region \eqref{eqn:square}, solution \eqref{eqn:bailey} belong to that 
part of region \eqref{eqn:cross} whose $(p,q)$-component corresponds to $E_1$ in 
Figure \ref{fig:pq-plane}.  
It is easy to see that \eqref{eqn:bailey} is of null deficiency $N=0$ 
if and only if $c$ satisfies the following generic condition: 
\[
c \not\in \mathbb{Z} \chi \cup (1+\mathbb{Z}\chi), \qquad \mbox{where} \quad 
\chi := \dfrac{2 \cdot \gcd\{j, \, k\}}{j+k} \,\, (\in \mathbb{Q}). 
\]
\par
There are a variety of studies on hypergeometric identities, 
especially, on gamma product formulas, and a lot of interesting 
formulas have been obtained not only for the Gauss hypergeometric 
series ${}_2F_1$ but also for its various generalizations. 
However, the study of necessary constraints for the existence of 
such identities lags far behind the well-developed 
ideas for discovering and verifying them, even in the most classical 
case of ${}_2F_1$.    
With the results in this article, our understanding of the former 
direction has advanced to some extent in region \eqref{eqn:square} 
and to a smaller extent in region \eqref{eqn:cross}, but remains  
almost null outside region \eqref{eqn:cross}. 
Even in region \eqref{eqn:square} we do not know whether $a$ and $b$ 
are always rational, although various evidences tempt us to 
guess positively. 
Note that the answer is certainly negative in region 
\eqref{eqn:cross} because of solution \eqref{eqn:bailey} and the 
possible existence of such a solution in region \eqref{eqn:square} 
makes the question much hard.        
\par
This article ends with a few examples of solutions outside region 
\eqref{eqn:cross}.   
The first formula of Table \ref{tab:erdelyi} is a solution whose 
$(p,q)$-component lies in region $G_4$ of Figure \ref{fig:pq-plane}. 
On the other hand,       
\begin{equation} \label{eqn:vidunas1}
{}_2F_1(4w-3, -4w+7;2w;1/2) = 2\pi \cdot 3^{7/2 -3 w} 
\dfrac{\varGamma(w+1/2)}{\varGamma(w-2/3)\varGamma(w-1/3)\varGamma(-w+5/2)}   
\end{equation}
is a solution corresponding to region $G_1$. 
This follows from formula (32) of Vidunas \cite{Vidunas0} by putting 
$c = 2w$ and using Pfaff's transformation \eqref{eqn:transf3} and 
the multiplication formula for the gamma function. 
A gamma factor $\varGamma(\pm w+\mbox{const.})$ is said to be 
{\sl positive} or {\sl negative} according to the choice of a sign. 
Note that our archetypal formula \eqref{eqn:gpf} has positive gamma 
factors only, while the present formula \eqref{eqn:vidunas1} contains a 
negative factor $\varGamma(-w+5/2)$ in its denominator.      
In any case, provided $r > 0$, every gamma factor in the numerator must 
be positive, since the function $f(w)$ has no poles in $\mathrm{Re}(w) > 0$.  
In region \eqref{eqn:cross} this is also the case with gamma 
factors in the denominator, because asymptotic formula \eqref{eqn:asympt2} 
shows that $f(w)$ has no zeros in some right half-plane,  
hence the setup of formula \eqref{eqn:gpf} is legitimate.   
For solution \eqref{eqn:vidunas1}, however, there are infinitely many 
zeros in both positive and negative directions along the real line. 
\par
If one wants to avoid the negative gamma factor in formula 
\eqref{eqn:vidunas1}, then Euler's reflection formula for the 
gamma function \cite[Theorem 1.2.1]{AAR} can be used to eliminate 
it, giving   
\begin{equation} \label{eqn:vidunas2}
{}_2F_1(4w-3, -4w+7;2w;1/2) = 2 \cdot 3^{7/2 -3 w} \cdot (\cos \pi w) \cdot  
\dfrac{\varGamma(w-3/2) \varGamma(w+1/2)}{\varGamma(w-2/3)\varGamma(w-1/3)},  
\end{equation}
where a {\sl trigonometric factor} appears instead. 
Formula \eqref{eqn:vidunas1} or \eqref{eqn:vidunas2} suggests that the 
archetypal formula \eqref{eqn:gpf} in region \eqref{eqn:cross} should 
be revised, when discussed in some other regions.   
Vidunas's formula mentioned above is equivalent to formula (i) of 
Maier \cite[Theorem 4.1]{Maier}, and formula (iv) of Maier's 
theorem leads to another solution of the same sort: 
\begin{align*}
{}_2F_1(6 w+2,-2 w;2 w;1/2) 
&= \frac{8\pi}{3} \cdot \left(w+\frac{1}{4}\right) \cdot 
\left(\frac{4}{27}\right)^w \cdot \dfrac{\varGamma(w)}{\varGamma(w+5/6)
\varGamma(w+7/6)\varGamma(-w)} \\[2mm]
&= - \frac{8}{3} \cdot \dfrac{w(w+1/4)}{w+1/6} \cdot 
\left(\frac{4}{27}\right)^w \cdot (\sin \pi w) \cdot 
\dfrac{\varGamma(w)^2}{\varGamma(w+1/6)\varGamma(w+5/6)}.   
\end{align*}
\par
More examples are derived from formula \eqref{eqn:bailey2} by 
putting $\alpha = (j+k)w+c$ and $\beta = (j-k)w$: 
\begin{align*}
& {}_2F_1(\, (j+k)w+c, \, -(j+k)w+1-c; \, (j-k)w; \, 1/2) \\[2mm] 
& \qquad = C \cdot d^w \cdot 
\prod_{\nu=0}^{k-1} \sin\pi\left(w+\frac{1+c}{2 k}+\frac{\nu}{k}\right) \cdot 
\dfrac{\prod_{\nu=0}^{j-k-1} \varGamma\big(w+\frac{\nu}{j-k}\big) 
\prod_{\nu=0}^{k-1} \varGamma\big(w+\frac{1+c}{2 k}+\frac{\nu}{k}\big)}{\prod_{\nu=0}^{j-1} 
\varGamma\big(w+\frac{c}{2 j}+\frac{\nu}{j}\big)},     
\end{align*}
where $j$ and $k$ are positive integers such that $j > k$, with $c$ being 
a free parameter and
\[ 
C :=  \dfrac{2^{k+1/2} \, k^{c/2}}{j^{(c-1)/2} \, (j-k)^{1/2}}, \qquad  
d := \dfrac{k^k (j-k)^{j-k}}{2^{j-k} \, j^j}. 
\]
A unified approach to exact evaluations with trigonometric factors 
is an interesting problem.        

\end{document}